\documentclass[12pt]{article}
\usepackage{graphicx}

\usepackage{epsfig}
\usepackage{amsmath}
\usepackage{amssymb}
\usepackage{amsthm}
\usepackage{latexsym}
\usepackage{cite}
\usepackage{amsbsy,amsfonts,euscript,eepic,rotating}

\newcommand{\comment}[1]{}
\newtheorem{theorem}{Theorem}
\newtheorem{lemma}{Lemma}[section]
\newtheorem{remark}{Remark}[section]
\newtheorem{corollary}{Corollary}[section]

\newtheorem{definition}{Definition}[section]
\begin{document}

\title{\LARGE
{\bf The Full Scaling Limit of \\ Two-Dimensional Critical Percolation}
}

\author{
{\bf Federico Camia}
\thanks{Research partially supported by a Marie Curie Intra-European Fellowship
under contract MEIF-CT-2003-500740.}\,
\thanks{E-mail: camia@eurandom.tue.nl}\\
{\small \sl EURANDOM, P.O. Box 513, 5600 MB Eindhoven, The Netherlands}\\
\and
{\bf Charles M.~Newman}
\thanks{Research partially supported by the
U.S. NSF under grant DMS-01-04278.}\,
\thanks{E-mail: newman@courant.nyu.edu}\\
{\small \sl Courant Inst.~of Mathematical Sciences,
New York University, New York, NY 10012, USA}
}

\date{}

\maketitle

\begin{abstract}
We use $SLE_6$ paths to construct a process of continuum nonsimple loops
in the plane and prove that this process coincides with the full continuum
scaling limit of 2D critical site percolation on the triangular lattice --
that is, the scaling limit of the set of all interfaces between different
clusters.
Some properties of the loop process, including conformal invariance,
are also proved.
In the main body of the paper these results are proved while assuming,
as argued by Schramm and Smirnov, that the percolation exploration path
converges in distribution to the trace of chordal $SLE_6$.
Then, in a lengthy appendix, a detailed proof is provided for this
convergence to $SLE_6$, which itself relies on Smirnov's result that
crossing probabilities converge to Cardy's formula.
\end{abstract}



\section{Introduction and Motivation} \label{intro}
In the theory of critical phenomena it is usually assumed
that a physical system near a continuous phase transition
is characterized by a single length scale (the ``correlation
length'') in terms of which all other lengths should be measured.
When combined with the experimental observation that the
correlation length diverges at the phase transition, this
simple but strong assumption, known as the scaling hypothesis,
leads to the belief that at criticality the system has
no characteristic length, and is therefore invariant under
scale transformations.
This suggests that all thermodynamic functions at criticality
are homogeneous functions, and predicts the appearance of
power laws.
It also means that it should be possible to rescale a
critical system appropriately and obtain a continuum model
(the ``continuum scaling limit'') which may have more
symmetries and be easier to study than the original discrete
model defined on a lattice.

Indeed, thanks to the work of Polyakov~\cite{polyakov}
and others~\cite{bpz1,bpz2}, it was understood by physicists
since the early seventies that critical statistical mechanical
models should possess continuum scaling limits with a global
conformal invariance that goes beyond simple scale invariance,
as long as the discrete models have ``enough'' rotation invariance.
This property gives important information, enabling the determination
of two- and three-point functions at criticality, when they are
nonvanishing.
Because the conformal group is in general a finite dimensional Lie
group, the resulting constraints are limited in number; however,
the situation becomes particularly interesting in two dimensions,
since there every analytic function $\omega=f(z)$ defines a conformal
transformation, at least at points where $f'(z) \neq 0$.
As a consequence, the conformal group in two dimensions is
infinite-dimensional.

After this observation was made, a large number of critical
problems in two dimensions were analyzed using conformal
methods, which were applied, among others, to Ising and Potts
models, Brownian motion, Self-Avoiding Walk (SAW), percolation,
and Diffusion Limited Aggregation (DLA).
The large body of knowledge and techniques that resulted, starting
with the work of Belavin, Polyakov and Zamolodchikov~\cite{bpz1,bpz2}
in the early eighties, goes under the name of Conformal Field Theory
(CFT).
In two dimensions, one of the main goals of CFT and its most important
application to statistical mechanics is a complete classification
of all universality classes via irreducible representations of the
infinite-dimensional Virasoro algebra.

Partly because of the success of CFT, work in recent
years on critical phenomena seemed to slow down somewhat,
probably due to the feeling that most of the leading
problems had been resolved.
Nonetheless, however powerful and successful it may be,
CFT has some limitations and leaves various open problems.
First of all, the theory deals primarily with correlation
functions of {\it local} (or quasi-local) operators, and is
therefore not always the best tool to investigate other
quantities.
Secondly, given some critical lattice model, there is no
way, within the theory itself, of deciding to which CFT
it corresponds.
A third limitation, of a different nature, is due to the fact
that the methods of CFT, although very powerful, are generally
speaking not completely rigorous from a mathematical point
of view.

In a somewhat surprising twist, the most recent developments
in the area of two-dimensional critical phenomena have
emerged in the mathematics literature and have followed
a new direction, which has provided new tools and a way of
coping with at least some of the limitations of CFT.
The new approach may even provide a reinterpretation of CFT,
and seems to be complementary to the traditional one in the sense
that questions that are difficult to pose and/or answer within
CFT are easy and natural in this new approach and vice versa.

The main tool of this radically new approach is the
Stochastic Loewner Evolution ($SLE$), or Schramm Loewner
Evolution, as it is also known, introduced by
Schramm~\cite{schramm}.
The new approach, which is probabilistic in nature,
focuses directly on non-local structures that characterize
a given system, such as cluster boundaries in Ising, Potts
and percolation models, or loops in the $O(n)$ model.
At criticality, these non-local objects become, in the
continuum limit, random curves whose distributions can be
uniquely identified thanks to their conformal invariance
and a certain ``Markovian" property.
There is a one-parameter family of $SLE$s, indexed by
a positive real number $\kappa$, and they appear to be
the only possible candidates for the scaling limits of
interfaces of two-dimensional critical systems that are
believed to be conformally invariant.

In particular, substantial progress has been made in
recent years, thanks to $SLE$, in understanding the
fractal and conformally invariant nature of (the scaling
limit of) large percolation clusters, which has attracted
much attention and is of interest both for intrinsic
reasons, given the many applications of percolation,
and as a paradigm for the behavior of other systems.
The work of Schramm~\cite{schramm} and Smirnov~\cite{smirnov}
has identified the scaling limit of a certain percolation
interface with $SLE_6$, providing, along with the work of
Lawler-Schramm-Werner~\cite{lsw1,lsw5} and Smirnov-Werner~\cite{sw},
a confirmation of many results in the physics literature,
as well as some new results.

However, $SLE_6$ describes a single interface, which
can be obtained by imposing special boundary conditions,
and is not in itself sufficient to immediately describe
the full scaling limit of the system.
In fact, not only the nature and properties, but the very
existence of the full scaling limit remained an open question.
This is true of all models, such as Ising and Potts models,
that are represented in terms of clusters.
Werner~\cite{werner3} considered this problem in the context
of $SLE_{\kappa}$ for values of $\kappa$ between $8/3$ and $4$.
For percolation (corresponding to $\kappa = 6$), the same
problem was addressed in~\cite{cn}, where $SLE_6$ was used
to construct a random process of continuous loops in the
plane, which was identified with the full scaling limit of
critical two-dimensional percolation, but without detailed proofs.

In this paper, we complete the analysis of~\cite{cn}, making
rigorous the connection between the construction given there
and the full scaling limit of percolation, and we prove some
properties of the full scaling limit, the Continuum Nonsimple
Loop process, including (one version of) conformal invariance.
We do this in two parts.
First, we give proofs in which we assume the validity of what
we will call statement (S) (see Section~\ref{main-results}),
which is a specific version of the results of Schramm and of
Smirnov~\cite{schramm,smirnov,smirnov-long,smirnov1,smirnov-private-comm}
concerning convergence of percolation exploration paths to $SLE_6$
(see the discussion towards the end of Section~\ref{explo}).
Since no detailed proof of statement (S) (or indeed, any version of
convergence to $SLE_6$) has been available, in Appendix~\ref{convergence-to-sle}
we give a proof based only on that part of Smirnov's results about
the convergence of crossing probabilities to Cardy's formula~\cite{smirnov}
(see Theorem~\ref{cardy-smirnov} in Appendix~\ref{convergence-to-sle}).
We note that statement (S) is restricted to Jordan domains while no
such restriction is indicated in~\cite{smirnov,smirnov-long}.



The rest of the paper is organized as follows.
In Section~\ref{defs}, we give necessary definitions and
introduce $SLE_6$.
Section~\ref{construct} is devoted to the construction of
the Continuum Nonsimple Loop process.
In Section~\ref{lapa}, we introduce the discrete model and
a discrete construction analogous to the continuum one presented
in Section~\ref{construct}.
Most of the main results of this paper are stated in Section~\ref{main-results},
while Section~\ref{proofs} contains the proofs of those results, which use (S).
The long Appendix~\ref{convergence-to-sle} contains the proof of
statement (S) (it is a consequence of Corollary~\ref{conv-to-sle-1} there)
and the short Appendix~\ref{rado} contains convergence
results for sequences of conformal maps which are used throughout the paper.

We remark that although our proof in Appendix~\ref{convergence-to-sle}
of convergence of exploration paths to $SLE_6$ roughly follows
Smirnov's outline~\cite{smirnov,smirnov-long}, based on his
proof~\cite{smirnov,smirnov-long} of convergence of crossing probabilities
to Cardy's formula and on the Markovian properties of hulls and tips,
there are at least two technically significant modifications.
The first is that we use a different sequence of stopping times to
obtain a Markov chain approximation to $SLE_6$, which results in a
different geometry for the approximation (see Remark~\ref{remark-difference}).
The second is that the control of ``close encounters" by the exploration
path to the domain boundary is not handled by general results for
``three-arms" events at the boundary of a half-plane, but rather by
an argument based on continuity of crossing probabilities with respect
to domain boundaries (see Lemmas~\ref{double-crossing}, \ref{touching}, \ref{equal}
and~\ref{mushroom}).
Moreover, we cannot use directly Smirnov's result on convergence of
crossing probabilities (see Theorem~\ref{cardy-smirnov}), but need
an extended version which is given in Theorem~\ref{strong-cardy} of
Appendix~\ref{convergence-to-sle}.

We conclude by noting that the convergence results of Appendix~\ref{convergence-to-sle}
are sufficient not only for our purposes of obtaining the full scaling
limit, but also for obtaining the critical exponents (see~\cite{sw}).

\section{Preliminary Definitions} \label{defs}

We will find it convenient to identify the
real plane ${\mathbb R}^2$ and the complex plane $\mathbb C$.
We will also refer to the Riemann sphere ${\mathbb C} \cup \infty$
and the open upper half-plane $\mathbb H = \{ x+iy : y>0 \}$
(and its closure $\overline{\mathbb H}$),
where chordal $SLE$ will be defined (see Section~\ref{sle1}).
$\mathbb D$ will denote the open unit disc
${\mathbb D} = \{ z \in {\mathbb C} : |z|<1 \}$.

A domain $D$ of the complex plane $\mathbb C$ is a
nonempty, connected, open subset of $\mathbb C$; a
simply connected domain $D$ is said to be a Jordan
domain if its (topological) boundary $\partial D$
is a Jordan curve (i.e., a simple continuous loop).

We will make repeated use of Riemann's mapping theorem,
which states that if $D$ is any simply connected domain
other than the entire plane $\mathbb C$ and $z_0 \in D$,
then there is a unique conformal
map $f$ of $D$ onto $\mathbb D$ such that $f(z_0)=0$
and $f'(z_0)>0$.

\subsection{Compactification of ${\mathbb R}^2$}

When taking the scaling limit $\delta \to 0$ one can focus
on fixed finite regions, $\Lambda \subset {\mathbb R}^2$,
or consider the whole ${\mathbb R}^2$ at once.
The second option avoids dealing with boundary conditions,
but requires an appropriate choice of metric.

A convenient way of dealing with the whole ${\mathbb R}^2$
is to replace the Euclidean metric with a distance function
$\Delta(\cdot,\cdot)$ defined on
${\mathbb R}^2 \times {\mathbb R}^2$ by
\begin{equation}
\Delta(u,v) = \inf_{\varphi} \int (1 + | {\varphi} |^2)^{-1} \, ds,
\end{equation}
where the infimum is over all smooth curves $\varphi(s)$
joining $u$ with $v$, parametrized by arclength $s$, and
where $|\cdot|$ denotes the Euclidean norm.
This metric is equivalent to the Euclidean metric in bounded
regions, but it has the advantage of making ${\mathbb R}^2$
precompact.
Adding a single point at infinity yields the compact space
$\dot{\mathbb R}^2$ which is isometric, via stereographic
projection, to the two-dimensional sphere.

\subsection{The Space of Curves} \label{space}

In dealing with the scaling limit we use the approach of
Aizenman-Burchard~\cite{ab}.
Denote by ${\cal S}_R$ the complete separable metric space
of continuous curves in the closure $\overline{\mathbb D}_R$
of the disc ${\mathbb D}_R$ of radius $R$ with the metric~(\ref{distance})
defined below.
Curves are regarded as equivalence classes of continuous
functions from the unit interval to $\overline{\mathbb D}_R$,
modulo monotonic reparametrizations.
$\gamma$ will represent a particular curve and $\gamma(t)$ a
parametrization of $\gamma$; ${\cal F}$ will represent a set
of curves (more precisely, a closed subset of ${\cal S}_R$).
$\text{d}(\cdot,\cdot)$ will denote the uniform metric
on curves, defined by
\begin{equation} \label{distance}
\text{d} (\gamma_1,\gamma_2) \equiv \inf
\sup_{t \in [0,1]} |\gamma_1(t) - \gamma_2(t)|,
\end{equation}
where the infimum is over all choices of parametrizations
of $\gamma_1$ and $\gamma_2$ from the interval $[0,1]$.
The distance between two closed sets of curves is defined
by the induced Hausdorff metric as follows:
\begin{equation} \label{hausdorff}
\text{dist}({\cal F},{\cal F}') \leq \varepsilon
\Leftrightarrow \forall \, \gamma \in {\cal F}, \, \exists \,
\gamma' \in {\cal F}' \text{ with }
\text{d} (\gamma,\gamma') \leq \varepsilon,
\text{ and vice versa.}
\end{equation}
The space $\Omega_R$ of closed subsets of ${\cal S}_R$
(i.e., collections of curves in $\overline{\mathbb D}_R$)
with the metric~(\ref{hausdorff}) is also a complete
separable metric space.
We denote by ${\cal B}_R$ its Borel $\sigma$-algebra.

For each fixed $\delta>0$, the random curves that we consider
are polygonal paths
on the edges of the hexagonal lattice $\delta {\cal H}$,
dual to the triangular lattice $\delta {\cal T}$.
A superscript $\delta$ is added to indicate that the
curves correspond to a model with a ``short distance cutoff''
of magnitude $\delta$.

We will also consider the complete separable metric space ${\cal S}$
of continuous curves in $\dot{\mathbb R}^2$ with the distance
\begin{equation} \label{Distance}
\text{D} (\gamma_1,\gamma_2) \equiv \inf
\sup_{t \in [0,1]} \Delta(\gamma_1(t),\gamma_2(t)),
\end{equation}
where the infimum is again over all choices of parametrizations
of $\gamma_1$ and $\gamma_2$ from the interval $[0,1]$.
The distance between two closed sets of curves is again
defined by the induced Hausdorff metric as follows:
\begin{equation} \label{hausdorff-D}
\text{Dist}({\cal F},{\cal F}') \leq \varepsilon
\Leftrightarrow \forall \, \gamma \in {\cal F}, \, \exists \,
\gamma' \in {\cal F}' \text{ with }
\text{D} (\gamma,\gamma') \leq \varepsilon,
\text{ and vice versa.}
\end{equation}
The space $\Omega$ of closed sets of $\cal S$
(i.e., collections of curves in $\dot{\mathbb R}^2$)
with the metric~(\ref{hausdorff-D}) is also a complete
separable metric space.
We denote by ${\cal B}$ its Borel $\sigma$-algebra.

When we talk about convergence in distribution of random curves,
we always mean with respect to the uniform metric~(\ref{distance}),
while when we deal with closed collections of curves, we always
refer to the metric~(\ref{hausdorff}) or~(\ref{hausdorff-D}).

\begin{remark} \label{ab}
In this paper, the space $\Omega$ of closed sets of $\cal S$ is
generally used for collections of exploration paths and cluster
boundary loops and their scaling limits, $SLE_6$ paths and
continuum nonsimple loops.
There is one place however, in the statements and proofs of
Lemmas~\ref{double-crossing}, \ref{equal} and~\ref{mushroom},
where we also apply $\Omega$ in essentially the original setting
of Aizenman and Burchard~\cite{aizenman,ab}, i.e., for collections
of blue and yellow simple $\cal T$-paths (see Section~\ref{lapa}
for precise definitions) and their scaling limits.
The slight modification needed to keep track of both the paths
and their colors is easily managed.
\end{remark}

\subsection{Chordal $SLE$ in the Upper Half-Plane} \label{sle1}

The Stochastic Loewner Evolution ($SLE$) was introduced
by Schramm~\cite{schramm} as a tool for studying the
scaling limit of two-dimensional discrete (defined on a
lattice) probabilistic models whose scaling limits are
expected to be conformally invariant.
In this section we define the chordal version of $SLE$;
for more on the subject, the interested reader can consult
the original paper~\cite{schramm} as well as the fine
reviews by Lawler~\cite{lawler1}, Kager and Nienhuis~\cite{kn},
and Werner~\cite{werner4}, and Lawler's book~\cite{lawler2}.

Let $\mathbb H$ denote the upper half-plane.
For a given continuous real function $U_t$ with $U_0 = 0$,
define, for each $z \in \overline{\mathbb H}$, the function
$g_t(z)$ as the solution to the ODE
\begin{equation}
\partial_t g_t(z) = \frac{2}{g_t(z) - U_t},
\end{equation}
with $g_0(z) = z$.
This is well defined as long as $g_t(z) - U_t \neq 0$,
i.e., for all $t < T(z)$, where
\begin{equation}
T(z) \equiv \sup \{ t \geq 0 : \min_{s \in [0,t]} | g_s(z) - U_s| > 0 \}.
\end{equation}
Let $K_t \equiv \{ z \in \overline{\mathbb H} : T(z) \leq t \}$
and let ${\mathbb H}_t$ be the unbounded component of
${\mathbb H} \setminus K_t$; it can be shown that $K_t$ is bounded
and that $g_t$ is a conformal map from ${\mathbb H}_t$ onto $\mathbb H$.
For each $t$, it is possible to write $g_t(z)$ as
\begin{equation}
g_t(z) = z + \frac{2t}{z} + o(\frac{1}{z}),
\end{equation}
when $z \to \infty$.
The family $(K_t, t \geq 0)$ is called the {\bf Loewner chain}
associated to the driving function $(U_t, t \geq 0)$.

\begin{definition} \label{def-sle}
{\bf Chordal $SLE_{\kappa}$} is the Loewner chain $(K_t, t \geq 0)$
that is obtained when the driving function
$U_t = \sqrt{\kappa} B_t$ is $\sqrt{\kappa}$ times a standard
real-valued Brownian motion $(B_t, t \geq 0)$ with $B_0 = 0$.
\end{definition}

For all $\kappa \geq 0$, chordal $SLE_{\kappa}$ is almost surely generated
by a continuous random curve $\gamma$ in the sense that, for all $t \geq 0$,
${\mathbb H}_t \equiv {\mathbb H} \setminus K_t$ is the unbounded connected
component of ${\mathbb H} \setminus \gamma[0,t]$; $\gamma$ is called the
{\bf trace} of chordal $SLE_{\kappa}$.

\subsection{Chordal $SLE$ in an Arbitrary Simply Connected Domain} \label{sle2}

Let $D \subset {\mathbb C}$ ($D \neq {\mathbb C}$) be a simply
connected domain whose boundary is a continuous curve.
By Riemann's mapping theorem, there are (many) conformal maps
from the upper half-plane $\mathbb H$ onto $D$.
In particular, given two distinct points $a,b \in \partial D$
(or more accurately, two distinct prime ends), there exists a
conformal map $f$ from $\mathbb H$ onto $D$ such that $f(0)=a$
and $f(\infty) \equiv \lim_{|z| \to \infty} f(z) = b$.
In fact, the choice of the points $a$ and $b$ on the boundary
of $D$ only characterizes $f(\cdot)$ up to a multiplicative
factor, since $f(\lambda \, \cdot)$ would also do.

Suppose that $(K_t, t \geq 0)$ is a chordal $SLE_{\kappa}$ in
$\mathbb H$ as defined above; we define chordal $SLE_{\kappa}$
$(\tilde K_t, t \geq 0)$ in $D$ from $a$ to $b$ as the
image of the Loewner chain $(K_t, t \geq 0)$ under $f$.
It is possible to show, using scaling properties of
$SLE_{\kappa}$, that the law of $(\tilde K_t, t \geq 0)$
is unchanged, up to a linear time-change, if we replace
$f(\cdot)$ by $f(\lambda \, \cdot)$.
This makes it natural to consider $(\tilde K_t, t \geq 0)$ as
a process from $a$ to $b$ in $D$, ignoring the role of $f$.

We are interested in the case $\kappa = 6$, for which
$(K_t, t \geq 0)$ is generated by a continuous, nonsimple,
non-self-crossing curve $\gamma$ with Hausdorff dimension $7/4$.
We will denote by $\gamma_{D,a,b}$ the image of $\gamma$ under $f$
and call it the trace of chordal $SLE_6$ in $D$ from $a$ to $b$;
$\gamma_{D,a,b}$ is a continuous nonsimple curve inside $D$ from $a$
to $b$, and it can be given a parametrization $\gamma_{D,a,b}(t)$
such that $\gamma_{D,a,b}(0)=a$ and $\gamma_{D,a,b}(1)=b$, so that
we are in the metric framework described in Section~\ref{space}.
It will be convenient to think of $\gamma_{D,a,b}$ as an
oriented path, with orientation from $a$ to $b$.

\section{Construction of the Continuum Nonsimple Loops}
\label{construct}

\subsection{Construction of a Single Loop} \label{single-loop}

\begin{figure}[!ht]
\begin{center}
\includegraphics[width=8cm]{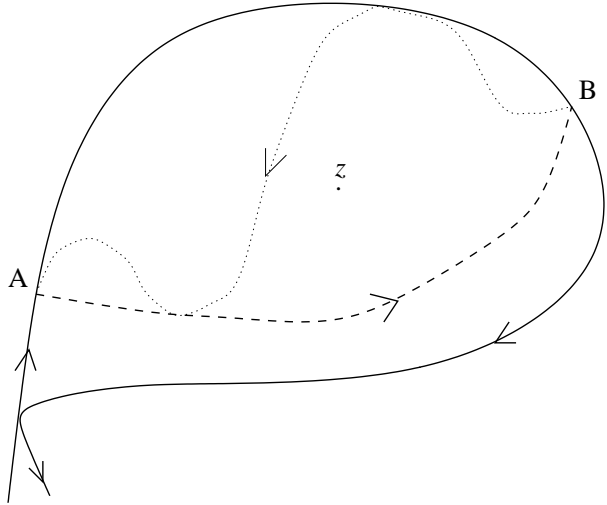}
\caption{Construction of a continuum loop around $z$ in three steps.
A domain $D$ is formed by the solid curved. The dashed curve is an
excursion $\cal E$ (from A to B) of an $SLE_6$ in $D$ that creates
a subdomain $D'$ containing $z$. The dotted curve $\gamma'$ is an
$SLE_6$ in $D'$ from B to A. A loop is formed by $\cal E$ followed
by $\gamma'$.}
\label{fig-sec3}
\end{center}
\end{figure}

As a preview to the full construction, we explain how to
construct a single loop using two $SLE_6$ paths inside a
domain $D$ whose boundary is assumed to have a given
orientation (clockwise or counterclockwise).
This is done in three steps (see Figure~\ref{fig-sec3}), of
which the first consists in choosing two points $a$ and $b$
on the boundary $\partial D$ of $D$ and ``running'' a chordal
$SLE_6$, $\gamma = \gamma_{D,a,b}$, from $a$ to $b$ inside $D$.
As explained in Section~\ref{sle2}, we consider $\gamma$
as an oriented path, with orientation from $a$ to $b$.
The set $D \setminus \gamma_{D,a,b}[0,1]$ is a countable
union of its connected components, which are open and simply
connected.
If $z$ is a deterministic point in $D$, then with probability
one, $z$ is not touched by $\gamma$~\cite{rs} and so it belongs
to a unique domain in $D \setminus \gamma_{D,a,b}[0,1]$
that we denote $D_{a,b}(z)$.

The elements of $D \setminus \gamma_{D,a,b}[0,1]$ can
be conveniently thought of in terms of how a point $z$ in
the interior of the component was first ``trapped'' at some
time $t_1$ by $\gamma[0,t_1]$, perhaps together with either
$\partial_{a,b} D$ or $\partial_{b,a} D$
(the portions of the boundary $\partial D$ from $a$ to $b$
counterclockwise or clockwise respectively):
(1) those components whose boundary contains a segment of
$\partial_{b,a} D$ between two successive visits at
$\gamma_0(z)=\gamma(t_0)$ and $\gamma_1(z)=\gamma(t_1)$ to
$\partial_{b,a} D$ (where here and below $t_0<t_1$), (2) the
analogous components with $\partial_{b,a} D$  replaced by the
other part of the boundary $\partial_{a,b} D$, (3) those components
formed when $\gamma_0(z)=\gamma(t_0)=\gamma(t_1)=\gamma_1(z) \in D$
with $\gamma$ winding about $z$ in a counterclockwise direction
between $t_0$ and $t_1$, and finally (4) the analogous clockwise
components.

We give to the boundary of a domain of type~3 or 4 the orientation
induced by how the curve $\gamma$ winds around the points inside
that domain.
For a domain $D' \ni z$ of type~1 or 2 which is produced by an
``excursion" $\cal E$ from $\gamma_0(z) \in \partial D$ to
$\gamma_1(z) \in \partial D$, the part of the boundary that
corresponds to the inner perimeter of the excursion $\cal E$
(i.e., the perimeter of $\gamma$ seen from $z$) is oriented
according to the direction of $\gamma$, i.e., from $\gamma_0(z)$
to $\gamma_1(z)$.

If we assume that $\partial D$ is oriented from $a$ to $b$
clockwise, then the boundaries of domains of type 2 have
a well defined orientation, while the boundaries of domains of
type 1 do not, since they are composed of two parts which are
both oriented from the beginning to the end of the excursion
that produced the domain.

Now, let $D'$ be a domain of type 1 and let $A$ and $B$ be
respectively the starting and ending point of the excursion
that generated $D'$.
The second step to construct a loop is to run a chordal
$SLE_6$, $\gamma' = \gamma_{D',B,A}$, inside $D'$ from
$B$ to $A$; the third and final step consists in pasting
together $\cal E$ and $\gamma'$.

Running $\gamma'$ inside $D'$ from $B$ to $A$
partitions $D' \setminus \gamma'$ into new domains.
Notice that if we assign an orientation to the boundaries of
these domains according to the same rules used above, all of
those boundaries have a well defined orientation, so that
the construction of loops just presented can be iterated
inside each one of these domains (as well as inside each
of the domains of type~2, 3 and 4 generated by $\gamma_{D,a,b}$
in the first step).
This will be done in the next section.

%
%

\subsection{The Full Construction Inside The Unit Disc}
\label{unit-disc}

In this section we define the Continuum Nonsimple Loop
process inside the unit disc  ${\mathbb D} = {\mathbb D}_1$
via an inductive procedure.
Later, in order to define the continuum nonsimple loops in
the whole plane, the unit disc will be replaced by a growing
sequence of large discs, ${\mathbb D}_R$, with $R\to \infty$
(see Theorem~\ref{thm-therm-lim}).
The basic ingredient in the algorithmic construction, given
in the previous section, consists of a chordal $SLE_6$ path
$\gamma_{D,a,b}$ between two points $a$ and $b$ of the boundary
$\partial D$ of a given simply connected domain $D \subset {\mathbb C}$.

We will organize the inductive procedure in steps,
each one corresponding to one $SLE_6$ inside a certain domain
generated by the previous steps.
To do that, we need to order the domains present at the end
of each step, so as to choose the one to use in the next step.
For this purpose, we introduce a deterministic countable
set of points $\cal P$ that are dense in $\mathbb C$ and
are endowed with a deterministic order (here and below
by deterministic we mean that they are assigned before
the beginning of the construction and are independent
of the $SLE_6$'s).

The first step consists of an $SLE_6$ path,
$\gamma_1 = \gamma_{{\mathbb D},-i,i}$, inside $\mathbb D$
from $-i$ to $i$, which produces many domains that are
the  connected components of the set
$\mathbb D \setminus \gamma_1 [0,1]$.
These domains can be priority-ordered according to the
maximal $x$- or $y$- coordinate distances between points on
their boundaries and using the rank of the points in $\cal P$
(contained in the domains) to break ties, as follows.
For a domain $D$, let $\text{d}_m(D)$ be the maximal $x$-
or $y$-distance between points on its boundary, whichever
is greater.
Domains with larger $\text{d}_m$ have higher priority, and
if two domains have the same $\text{d}_m$, the one containing
the highest ranking point of $\cal P$ from those two domains
has higher priority.
The priority order of domains of course changes as the
construction proceeds and new domains are formed.

The second step of the construction consists of an $SLE_6$
path, $\gamma_2$, that is produced in the domain with highest
priority (after the first step).
Since all the domains that are produced in the construction
are Jordan domains, as explained in the discussion following
Corollary~\ref{jordan}, for all steps we can use the definition
of chordal $SLE$ given in Section~\ref{sle2}.

As a result of the construction, the $SLE_6$ paths are
naturally ordered: $\{ \gamma_j \}_{j \in {\mathbb N}}$.
It will be shown (see especially the proof of Theorem~\ref{thm-convergence}
below) that every domain that is formed during the construction
is eventually used (this is in fact one important requirement in
deciding how to order the domains and therefore how to organize
the construction).

So far we have not explained how to choose the starting
and ending points of the $SLE_6$ paths on the boundaries
of the domains.
In order to do this, we give an orientation to the
boundaries of the domains produced by the construction
according to the rules explained in Section~\ref{single-loop}.
We call {\bf monochromatic} a boundary which gets, as
a consequence of those rules, a well defined (clockwise
or counterclockwise) orientation; the choice of this term
will be clarified when we discuss the lattice version of
the loop construction below.
We will generally take our initial domain ${\mathbb D}_1$
(or ${\mathbb D}_R$) to have a monochromatic boundary
(either clockwise or counterclockwise orientation).

It is easy to see by induction that the boundaries that
are not monochromatic are composed of two ``pieces'' joined
at two special points (call them A and B, as in the example
of Section~\ref{single-loop}), such that one piece is a
portion of the boundary of a previous domain, and the
other is the inner perimeter of an excursion (see again
Section~\ref{single-loop}).
Both pieces are oriented in the same direction, say from A
to B (see Figure~\ref{fig-sec3}).

For a domain whose boundary is not monochromatic, we make
the ``natural'' choice of starting and ending points,
corresponding to the end and beginning of the excursion
that produced the domain (the points B and A respectively,
in the example above).
As explained in Section~\ref{single-loop}, when such
a domain is used with this choice of points on the
boundary, a loop is produced, together with other domains,
whose boundaries are all monochromatic.

For a domain whose boundary is monochromatic, and therefore
has a well defined orientation, there are various procedures
which would yield the ``correct'' distribution for the resulting
Continuum Nonsimple Loop process; one possibility is as follows.

Given a domain $D$, $a$ and $b$ are chosen so that, of all pairs
$(u,v)$ of points in $\partial D$, they maximize $|\text{Re}(u-v)|$
if $|\text{Re}(u-v)| \geq |\text{Im}(u-v)|$, or else they maximize
$|\text{Im}(u-v)|$.
If the choice is not unique, to restrict the number of pairs
one looks at those pairs, among the ones already obtained, that
maximize the other of $\{ |\text{Re}(u-v)|, |\text{Im}(u-v)| \}$.
Notice that this leaves at most two pairs of points; if that's
the case, the pair that contains the point with minimal real
(and, if necessary, imaginary) part is chosen.
The iterative procedure produces a loop every time a domain
whose boundary is not monochromatic is used.
Our basic loop process consists of the collection of all loops
generated by this inductive procedure (i.e., the limiting object
obtained from the construction by letting the number of steps
$k \to \infty$), to which we add a ``trivial" loop for each $z$
in $\mathbb D$, so  that the collection of loops is closed in
the appropriate sense~\cite{ab}.
The Continuum Nonsimple Loop process in the whole plane is
introduced in Theorem~\ref{thm-therm-lim}, Section~\ref{main-results}.
There, a ``trivial" loop for each $z \in {\mathbb C} \cup \infty$
has to be added to make the space of loops closed.

\section{Lattices and Paths} \label{lapa}

\begin{figure}[!ht]
\begin{center}
\includegraphics[width=6cm]{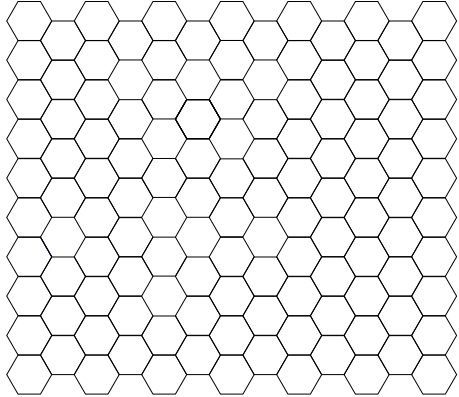}
\caption{Portion of the hexagonal lattice.}
\label{fig1-sec4}
\end{center}
\end{figure}

We will denote by $\cal T$ the two-dimensional triangular lattice,
whose sites we think of as the elementary cells of a regular hexagonal
lattice $\cal H$ embedded in the plane as in Figure~\ref{fig1-sec4}.
A sequence $(\xi_0, \ldots, \xi_n)$ of sites of $\cal T$ such that
$\xi_{i-1}$ and $\xi_i$ are neighbors in $\cal T$ for all
$i= 1, \ldots, n$ and $\xi_i \neq \xi_j$ whenever $i \neq j$
will be called a {\bf $\cal T$-path} and denoted by $\pi$.
If the first and last sites of the path are neighbors in $\cal T$,
the path will be called a {\bf $\cal T$-loop}.

We say that a finite subset $D$ of $\cal T$ is {\bf simply connected}
if both $D$ and ${\cal T} \setminus D$ are connected (by the edges
of $\cal T$).
For a simply connected set $D$ of hexagons, we denote by
$\Delta D$ its {\bf external site boundary}, or {\bf s-boundary}
(i.e., the set of hexagons that do not belong to $D$
but are adjacent to hexagons in $D$), and by $\partial D$ the
topological boundary of $D$ when $D$ is considered as a domain of
$\mathbb C$.
We will call a bounded, simply connected subset $D$ of $\cal T$
a {\bf Jordan set} if its s-boundary $\Delta D$ is a $\cal T$-loop.

For a Jordan set $D \subset {\cal T}$, a vertex $x \in {\cal H}$
that belongs to $\partial D$ can be either of two types, according to whether
the edge incident on $x$ that is not in $\partial D$ belongs to a hexagon
in $D$ or not.
We call a vertex of the second type an {\bf e-vertex} (e for ``external''
or ``exposed'').

Given a Jordan set $D$ and two e-vertices $x,y$ in $\partial D$,
we denote by $\partial_{x,y} D$ the portion of $\partial D$
traversed counterclockwise from $x$ to $y$, and call it the
{\bf right boundary}; the remaining part of the boundary is
denote by $\partial_{y,x} D$ and is called the {\bf left boundary}.
Analogously, the portion of $\Delta_{x,y} D$ of $\Delta D$ whose
hexagons are adjacent to $\partial_{x,y} D$ is called the
{\bf right s-boundary} and the remaining part the {\bf left s-boundary}.

A {\bf percolation configuration}
$\sigma = \{ \sigma(\xi) \}_{\xi \in \cal T} \in \{ -1, +1 \}^{\cal T}$
on $\cal T$ is an assignment of $-1$ (equivalently, yellow) or $+1$
(blue) to each site of $\cal T$.
For a domain $D$ of the plane, the restriction to the subset
$D \cap \cal T$ of $\cal T$ of the percolation configuration
$\sigma$ is denoted by $\sigma_D$.
On the space of configurations $\Sigma = \{ -1,+1 \}^{\cal T}$,
we consider the usual product topology and denote by $\mathbb P$
the uniform measure, corresponding to Bernoulli percolation with
equal density of yellow (minus) and blue (plus) hexagons, which is
critical percolation in the case of the triangular lattice.

A (percolation) {\bf cluster} is a maximal, connected, monochromatic
subset of $\cal T$; we will distinguish between blue (plus) and
yellow (minus) clusters.
The {\bf boundary} of a cluster $D$ is the set of edges of
$\cal H$ that surround the cluster (i.e., its Peierls contour);
it coincides with the topological boundary of $D$ considered as a
domain of $\mathbb C$.
The set of all boundaries is a collection of ``nested'' simple loops
along the edges of $\cal H$.

Given a percolation configuration $\sigma$, we associate an
arrow to each edge of $\cal H$ belonging to the boundary of
a cluster in such a way that the hexagon to the right of the edge
with respect to the direction of the arrow is blue (plus).
The set of all boundaries then becomes a collection of nested,
oriented, simple loops.
A {\bf boundary path} (or {\bf b-path}) $\gamma$ is a sequence
$(e_0, \ldots, e_n)$ of distinct edges of $\cal H$ belonging
to the boundary of a cluster and such that $e_{i-1}$ and $e_i$
meet at a vertex of $\cal H$ for all $i= 1, \ldots, n$.
To each b-path, we can associate a direction according to the
direction of the edges in the path.

Given a b-path $\gamma$, we denote by $\Gamma_B(\gamma)$
(respectively, $\Gamma_Y(\gamma)$) the set of blue (resp.,
yellow) hexagons (i.e., sites of $\cal T$) adjacent to
$\gamma$; we also let
$\Gamma(\gamma) \equiv \Gamma_B(\gamma) \cup \Gamma_Y(\gamma)$.

\subsection{The Percolation Exploration Process and Path} \label{explo}

\begin{figure}[!ht]
\begin{center}
\includegraphics[width=6cm]{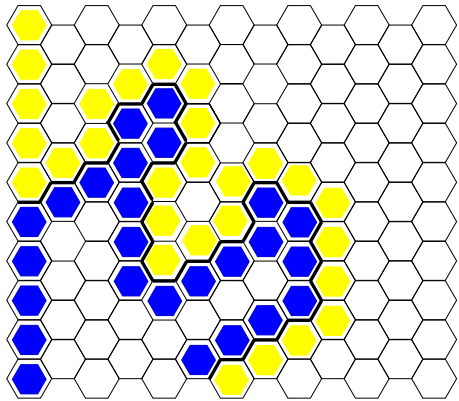}
\caption{Percolation exploration process in a portion of the
hexagonal lattice with $\pm$ boundary conditions on the first
column, corresponding to the boundary of the region where the
exploration is carried out.
The colored hexagons that do not belong to the first column
have been ``explored'' during the exploration process.
The heavy line between yellow (light) and blue (dark) hexagons
is the exploration path produced by the exploration process.}
\label{fig2-sec4}
\end{center}
\end{figure}

For a Jordan set $D \subset {\cal T}$ and two
e-vertices $x,y$ in $\partial D$, imagine coloring blue
all the hexagons in $\Delta_{x,y} D$ and yellow all those
in $\Delta_{y,x} D$.
Then, for any percolation configuration $\sigma_D$ inside
$D$, there is a unique b-path $\gamma$ from $x$ to $y$ which
separates the blue cluster adjacent to $\Delta_{x,y} D$ from
the yellow cluster adjacent to $\Delta_{y,x} D$.
We call $\gamma = \gamma_{D,x,y}(\sigma_D)$ a
{\bf percolation exploration path} (see Figure~\ref{fig2-sec4}).

An exploration path $\gamma$ can be decomposed into
{\bf left excursions} $\cal E$, i.e., maximal b-subpaths
of $\gamma$ that do not use edges of the left boundary $\partial_{y,x} D$.
Successive left excursions are separated by portions of
$\gamma$ that contain only edges of the left boundary $\partial_{y,x} D$.
Analogously, $\gamma$ can be decomposed into {\bf right excursions},
i.e., maximal b-subpaths of $\gamma$ that do not use edges of the right
boundary $\partial_{x,y} D$.
Successive right excursions are separated by portions of
$\gamma$ that contain only edges of the right boundary $\partial_{x,y} D$.

Notice that the exploration path
$\gamma=\gamma_{D,x,y}(\sigma_D)$ only depends on
the percolation configuration $\sigma_D$ inside $D$
and the positions of the e-vertices $x$ and $y$;
in particular, it does not depend on the color of
the hexagons in $\Delta D$, since it is defined by
imposing fictitious $\pm$ boundary conditions on $D$.
To see this more clearly, we next show how to construct
the percolation exploration path dynamically, via
the {\bf percolation exploration process} defined
below.

Given a Jordan set $D \subset {\cal T}$ and two
e-vertices $x,y$ in $\partial D$, assign to $\partial_{x,y} D$
a counterclockwise orientation (i.e., from $x$ to $y$) and
to $\partial_{y,x} D$ a clockwise orientation.
Call $e_x$ the edge incident on $x$ that does not belong to
$\partial D$ and orient it in the direction of $x$;
this is the ``starting edge'' of an exploration procedure
that will produce an oriented path inside $D$ along
the edges of $\cal H$, together with two \emph{nonsimple}
monochromatic paths on $\cal T$.
From $e_x$, the process moves along the edges of
hexagons in $D$ according to the rules below.
At each step there are two possible edges (left or right edge
with respect to the current direction of exploration) to choose
from, both belonging to the same hexagon $\xi$ contained in $D$
or $\Delta D$.

\begin{itemize}
\item If $\xi$ belongs to $D$ and has not been previously ``explored,"
its color is determined by flipping a fair coin and then
the edge to the left (with respect to the direction in which
the exploration is moving) is chosen if $\xi$ is blue (plus),
or the edge to the right is chosen if $\xi$ is yellow (minus).
\item If $\xi$ belongs to $D$ and has been previously explored, the
color already assigned to it is used to choose an edge according to
the rule above.
\item If $\xi$ belongs to the right external boundary $\Delta_{x,y} D$,
the left edge is chosen.
\item If $\xi$ belongs to the left external boundary $\Delta_{y,x} D$,
the right edge is chosen.
\item The exploration process stops when it reaches $b$.
\end{itemize}

We can assign an arrow to each edge in the path in such
a way that the hexagon to the right of the edge with respect
to the arrow is blue; for edges in $\partial D$, we assign the
arrows according to the direction assigned to the boundary.
In this way, we get an oriented path, whose shape and orientation
depend solely on the color of the hexagons explored during
the construction of the path.

%


When we present the discrete construction, we will encounter
Jordan sets $D$ with two e-vertices $x,y \in \partial D$
assigned in some way to be discussed later.
Such domains will have either monochromatic (plus or
minus) boundaries or $\pm$ boundary conditions,
corresponding to having both $\Delta_{x,y} D$ and
$\Delta_{y,x} D$ monochromatic, but of different
colors.

As explained, the exploration path $\gamma_{D,x,y}$
does not depend on the color of $\Delta D$, but the
interpretation of $\gamma_{D,x,y}$ does.
For domains with $\pm$ boundary conditions, the
exploration path represents the interface between the
yellow cluster containing the yellow portion of the
s-boundary of $D$ and the blue cluster containing its
blue portion.

For domains with monochromatic blue (resp., yellow)
boundary conditions, the exploration path represents portions
of the boundaries of yellow (resp., blue)
clusters touching $\partial_{y,x} D$ and adjacent to
blue (resp., yellow) hexagons that are the starting point
of a blue (resp., yellow) path (possibly an empty path)
that reaches $\partial_{x,y} D$, pasted together using
portions of $\partial_{y,x} D$.

In order to study the continuum scaling limit of an exploration
path, we introduce the following definitions.
\begin{definition} \label{approx}
Given a bounded, simply connected domain $D$ of the plane,
we denote by $D^{\delta}$ the largest Jordan set of hexagons
of the scaled hexagonal lattice $\delta {\cal H}$ that is
contained in $D$, and call it the {\bf $\delta$-approximation}
of $D$.
\end{definition}

It is clear that if $D$ is a Jordan domain, then as $\delta \to 0$,
$\partial D^{\delta}$ converges to $\partial D$ in the
metric~(\ref{distance}).

\begin{definition} \label{exp-path}
Let $D$ be a bounded domain of the plane and $D^{\delta}$
its $\delta$-approximation.
For $a,b \in \partial D$, choose the pair $(x_a,x_b)$
of e-vertices in $\partial D^{\delta}$ closest to, respectively,
$a$ and $b$ (if there are two such vertices closest to $a$,
we choose, say, the first one encountered going clockwise
along $\partial D^{\delta}$, and analogously for $b$).
Given a percolation configuration $\sigma$, we define
the {\bf exploration path} $\gamma^{\delta}_{D,a,b}(\sigma)
\equiv \gamma_{D^{\delta},x_a,x_b}(\sigma)$.
\end{definition}

For a fixed $\delta>0$, the measure $\mathbb P$ on percolation
configurations $\sigma$ induces a measure $\mu^{\delta}_{D,a,b}$
on exploration paths $\gamma^{\delta}_{D,a,b}(\sigma)$.
In the continuum scaling limit, $\delta \to 0$, one is
interested in the weak convergence of $\mu^{\delta}_{D,a,b}$ to a
measure $\mu_{D,a,b}$ supported on continuous curves, with respect
to the uniform metric~(\ref{distance}) on continuous curves.

One of the main tools in this paper is the result on convergence
to $SLE_6$ announced by Smirnov~\cite{smirnov} (see also~\cite{smirnov-long}),
whose detailed proof is to appear~\cite{smirnov1}: \emph{The distribution
of $\gamma^{\delta}_{D,a,b}$ converges, as $\delta \to 0$, to that
of the trace of chordal $SLE_6$ inside $D$ from $a$ to $b$, with
respect to the uniform metric~(\ref{distance}) on continuous curves}.

Actually, we will rather use a slightly stronger conclusion,
given as statement (S) at the beginning of Section~\ref{main-results}
below, a version of which, according to~\cite{sw} (see p.~734 there),
and~\cite{smirnov-private-comm}, will be contained in~\cite{smirnov1}.
This stronger statement is that the convergence of the percolation process
to $SLE_6$ takes place \emph{locally uniformly} with respect to the shape
of the domain $D$ and the positions of the starting and ending points $a$
and $b$ on its boundary $\partial D$.
We will use this version of convergence to $SLE_6$ to identify the
Continuum Nonsimple Loop process with the scaling limit of \emph{all}
critical percolation clusters.
Statement (S) is a direct consequence of Corollary~\ref{conv-to-sle-1},
which is proved in Appendix~\ref{convergence-to-sle}.
Although the convergence statements in Corollary~\ref{conv-to-sle-1}
and in (S) are stronger than those in~\cite{smirnov,smirnov-long}, we
note that they are restricted to Jordan domains, a restriction not
present in~\cite{smirnov,smirnov-long}.



Before concluding this section, we give one more definition.
Consider the exploration path $\gamma = \gamma^{\delta}_{D,x,y}$
and the set $\Gamma(\gamma) = \Gamma_Y(\gamma) \cup \Gamma_B(\gamma)$.
The set $D^{\delta} \setminus \Gamma(\gamma)$ is the
union of its connected components (in the lattice sense),
which are simply connected.
If the domain $D$ is large and the e-vertices
$x_a, y_a \in \partial D^{\delta}$ are not too close
to each other, then with high probability the exploration
process inside $D^{\delta}$ will make large excursions
into $D^{\delta}$, so that
$D^{\delta} \setminus \Gamma(\gamma)$ will have more
than one component.
Given a point $z \in {\mathbb C}$ contained in
$D^{\delta} \setminus \Gamma(\gamma)$,
we will denote by $D^{\delta}_{a,b}(z)$ the domain
corresponding to the unique element of
$D^{\delta} \setminus \Gamma(\gamma)$ that contains
$z$ (notice that for a deterministic
$z \in D$, $D^{\delta}_{a,b}(z)$ is well defined
with high probability for $\delta$ small, i.e., when
$z \in D^{\delta}$ and $z \notin \Gamma(\gamma)$).

\subsection{Discrete Loop Construction}

\begin{figure}[!ht]
\begin{center}
\includegraphics[width=8cm]{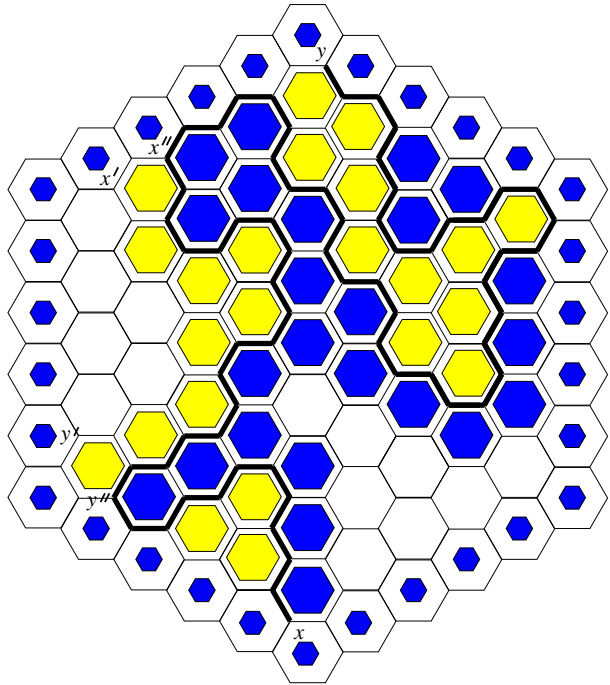}
\caption{First step of the construction of the outer contour
of a cluster of yellow/minus (light in the figure) hexagons
consists of an exploration from the vertex $x$ to $y$ (heavy line).
The outer layer of hexagons does not belong to the domain
where the explorations are carried out, but represents its
monochromatic blue/plus external site boundary.
$x''$ and $y''$ are the ending and starting points of a left
excursion that determines a new domain $D'$, and $x'$ and $y'$
are the vertices where the edges that separate the yellow and
blue portions of the s-boundary of $D'$ intersect $\partial D'$.}
\label{fig3-sec4}
\end{center}
\end{figure}

\begin{figure}[!ht]
\begin{center}
\includegraphics[width=8cm]{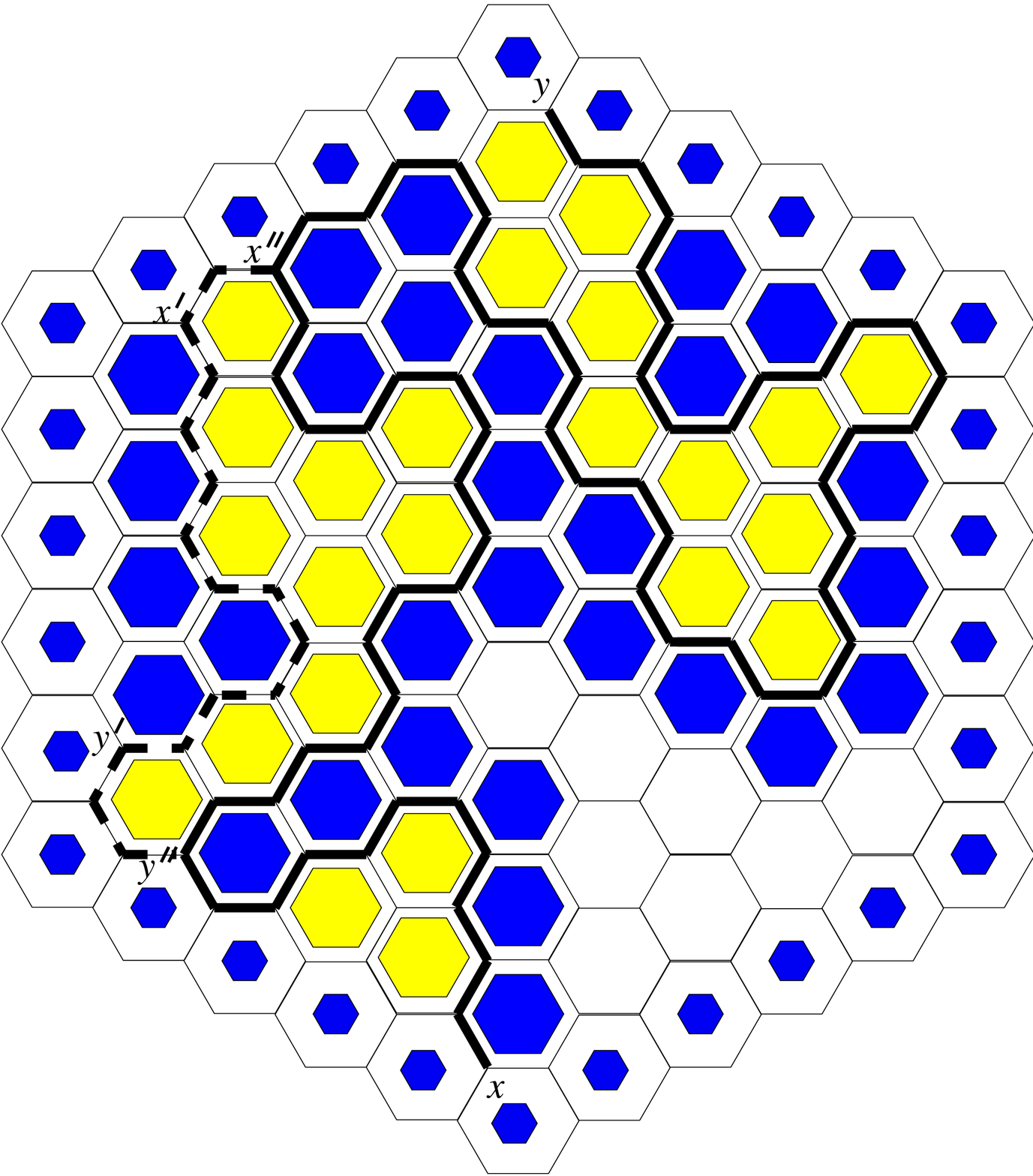}
\caption{Second step of the construction of the outer contour
of a cluster of yellow/minus (light in the figure) hexagons
consisting of an exploration from $x'$ to $y'$ whose resulting
path (heavy broken line) is pasted to the left excursion generated
by the previous exploration with the help of edges (indicated again
by a heavy broken line) of $\partial D$ contained between $x'$ and
$x''$ and between $y'$ and $y''$.}
\label{fig4-sec4}
\end{center}
\end{figure}

Next, we show how to construct, by twice using the
exploration process described in Section~\ref{explo},
a loop $\Lambda$ along the edges of ${\cal H}$
corresponding to the external boundary of a monochromatic
cluster contained in a large, simply connected, Jordan
set $D$ with monochromatic blue (say) boundary conditions
(see Figures~\ref{fig3-sec4} and~\ref{fig4-sec4}).

Consider the exploration path $\gamma = \gamma_{D,x,y}$
and the sets $\Gamma_Y(\gamma)$ and $\Gamma_B(\gamma)$
(see Figure~\ref{fig3-sec4}).
The set $D \setminus \{ \Gamma_Y(\gamma) \cup  \Gamma_B(\gamma) \}$
is the union of its connected components (in the lattice sense),
which are simply connected.
If the domain $D$ is large and the e-vertices
$x, y \in \partial D$ are chosen not too close
to each other, with large probability the exploration
process inside $D$ will make large excursions into
$D$, so that
$D \setminus \{ \Gamma_Y(\gamma) \cup  \Gamma_B(\gamma) \}$
will have many components.

There are four types of components which may be usefully
thought of in terms of their external site boundaries:
(1) those components  whose  site boundary contains both
sites in $\Gamma_Y(\gamma)$ and  $\Delta_{y,x} D$, (2)
the analogous components  with $\Delta_{y,x} D$ replaced
by  $\Delta_{x,y} D$ and $\Gamma_Y(\gamma)$ by
$\Gamma_B(\gamma)$, (3) those components whose site
boundary only contains sites  in $\Gamma_Y(\delta)$,
and finally (4) the analogous components  with
$\Gamma_Y(\gamma)$  replaced by $\Gamma_B(\gamma)$.

Notice that the components of type~1 are the only
ones with $\pm$ boundary conditions, while all other
components have monochromatic s-boundaries.
For a given component $D'$ of type~1, we can identify
the two edges that separate the yellow and blue portions
of its s-boundary.
The vertices $x'$ and $y'$ of $\cal H$ where those
two edges intersect $\partial D'$ are e-vertices and
are chosen to be the starting and ending points of
the exploration path $\gamma_{D',x',y'}$ inside $D'$.

If $x'', y'' \in \partial D$ are respectively the ending
and starting points of the left excursion $\cal E$ of $\gamma_{D,x,y}$
that ``created" $D'$, by pasting together $\cal E$ and
$\gamma_{D',x',y'}$ with the help of the edges of $\partial D$
contained between $x'$ and $x''$ and between $y'$ and $y''$,
we get a loop $\Lambda$ which corresponds to the
boundary of a yellow cluster adjacent to $\partial_{y,x} D$
(see Figure~\ref{fig4-sec4}).
Notice that the path $\gamma_{D',x',y'}$ in general splits
$D'$ into various other domains, all of which have monochromatic
boundary conditions.

%
%

\subsection{Full Discrete Construction} \label{full}

We now give the algorithmic construction for discrete
percolation  which is the analogue of the continuum one.
Each step of the construction is a single percolation
exploration  process; the order of successive steps is
organized as in the continuum construction detailed in
Section~\ref{unit-disc}.
We start with the smallest Jordan set
$D^{\delta}_0 = {\mathbb D}^{\delta}$
of hexagons that covers the unit disc $\mathbb D$.
We will also make use of the countable set $\cal P$ of
points dense in $\mathbb C$ that was introduced earlier.

The first step consists of an exploration process inside
$D^{\delta}_0$.
For this, we need to select two points $x$ and $y$ in
$\partial D^{\delta}_0$ (which identify the starting and
ending edges).
We choose for $x$ the e-vertex closest to $-i$, and for
$y$ the e-vertex closest to $i$ (if there are two such
vertices closest to $-i$, we can choose, say, the one with
smallest real part, and analogously for $i$).
The first exploration produces a path $\gamma^{\delta}_1$
and, for $\delta$ small, many new domains of all four types.
These domains are ordered according to the maximal $x$- or
$y$- distance $\text{d}_m$ between points on their boundaries
and, if necessary, with the help of points in $\cal P$, as in
the continuum case, and that order is used, at each step of the
construction, to determine the next exploration process.
With this choice, the exploration processes and paths are
naturally ordered: $\gamma^{\delta}_1, \gamma^{\delta}_2, \ldots$ .

Each exploration process of course requires choosing a starting
and ending vertex and edge.
For domains of type~1, with a $\pm$ or $\mp$ boundary condition,
the choice is the natural one, explained before.

For a domain $D^{\delta}_k$ (used at the $k$th step) of type
other than~1, and therefore with a monochromatic boundary, the
starting and ending edges are chosen with a procedure that mimics
what is done in the continuum case. 
Once again, the exact procedure used to choose the pair of
points is not important, as long as they are not chosen too
close to each other.
This is clear in the discrete case because the procedure
that we are presenting is only ``discovering'' the cluster
boundaries.
In more precise terms, it is clear that one could couple
the processes obtained with different rules by means of the
same percolation configuration, thus obtaining exactly the
same cluster boundaries.

As in the continuum case, we can choose the following procedure.
(In Theorem~\ref{thm-convergence} we will slightly reorganize the
procedure by using a coupling to the continuum construction to
guarantee that the order of exploration of domains of the discrete
and continuum procedures match despite the rules for breaking ties.)
Given a domain $D$, $x$ and $y$ are chosen so that, of all pairs
$(u,v)$ of points in $\partial D$, they maximize $|\text{Re}(u-v)|$
if $|\text{Re}(u-v)| \geq |\text{Im}(u-v)|$, or else they maximize
$|\text{Im}(u-v)|$.
If the choice is not unique, to restrict the number of pairs
one looks at those pairs, among the ones already obtained, that
maximize the other of $\{ |\text{Re}(u-v)|, |\text{Im}(u-v)| \}$.
Notice that this leaves at most two pairs of points; if that's
the case, the pair that contains the point with minimal real
(and, if necessary, imaginary) part is chosen.

The procedure continues iteratively, with regions that
have monochromatic boundaries playing  the role played
in the first step by the unit disc.
Every time a region with $\pm$ boundary conditions is used,
a new loop, corresponding to the outer boundary contour of
a cluster, is formed by pasting together, as explained in
Section~\ref{single-loop}, the new exploration path and the
excursion containing the region where the last exploration
was carried out.
All the new regions created at a step when a loop is formed
have monochromatic boundary conditions.

\section{Main Results} \label{main-results}

In this section we collect our main results about the Continuum
Nonsimple Loop process.
Before doing that, we state a precise version, called statement (S),
on convergence of exploration paths to $SLE_6$ that we will use in
the proofs of these results, presented in Section~\ref{proofs}.
(A proof of statement (S) is given in Appendix~\ref{convergence-to-sle};
it is an immediate consequence of Corollary~\ref{conv-to-sle-1} there.
The proof relies, among other things, on the result of Smirnov~\cite{smirnov}
concerning convergence of crossing probabilities to Cardy's
formula~\cite{cardy,cardy2} -- see Theorem~\ref{cardy-smirnov}.)
We note that (S) or (Corollary~\ref{conv-to-sle-1}) is both more general
and more special than the convergence statements
in~\cite{smirnov,smirnov-long} --- more general in that the domain
can vary with $\delta$ as $\delta \to 0$, but more special in
the restriction to Jordan domains.

Given a Jordan domain $D$ with two distinct points $a,b \in \partial D$
on its boundary, let $\mu_{D,a,b}$ denote the law of $\gamma_{D,a,b}$,
the trace of chordal $SLE_6$, and let $\mu^{\delta}_{D,a,b}$ denote
the law of the percolation exploration path $\gamma^{\delta}_{D,a,b}$.
Let $X$ be the space of continuous curves inside $D$ from $a$ to $b$.
We define $\rho(\mu_{D,a,b},\mu^{\delta}_{D,a,b}) \equiv
\inf \{\varepsilon>0 : \mu_{D,a,b}(U) \leq
\mu^{\delta}_{D,a,b}(\bigcup_{x \in U}B_{\text{d}}(x,\varepsilon))
+ \varepsilon \text{ for all Borel } U \subset X \}$
(where $B_{\text{d}}(x,\varepsilon)$ denotes the open ball of
radius $\varepsilon$ centered at $x$ in the metric~(\ref{distance}))
and denote by
$\text{d}_{\text{P}}(\mu_{D,a,b},\mu^{\delta}_{D,a,b}) \equiv
\max \{ \rho(\mu_{D,a,b},\mu^{\delta}_{D,a,b}), \rho(\mu^{\delta}_{D,a,b},\mu_{D,a,b}) \}$
the Prohorov distance; weak convergence is equivalent to
convergence in the Prohorov metric.
Statement (S) is the following; it is used in the proofs of all the results
of this section {\it except\/} for Lemmas~\ref{sub-conv}-\ref{boundaries}.
\begin{itemize}
\item[(S)] For Jordan domains, there is convergence in distribution
of the percolation exploration path to the trace of chordal $SLE_6$
that is \emph{locally uniform} in the shape of the boundary with respect
to the uniform metric on continuous curves~(\ref{distance}), and in the
location of the starting and ending points with respect to the Euclidean
metric; i.e., for $(D,a,b)$ a Jordan domain with $a,b \in \partial D$,
$\forall \varepsilon>0$, $\exists \alpha_0=\alpha_0(\varepsilon)$ and
$\delta_0=\delta_0(\varepsilon)$ such that for all $(D',a',b')$ with
$D'$ Jordan and with
$\max{(\text{d}(\partial D, \partial D'),|a-a'|,|b-b'|) \leq \alpha_0}$
and $\delta \leq \delta_0$,
$\text{d}_{\text{P}}(\mu_{D',a',b'},\mu^{\delta}_{D',a',b'}) \leq \varepsilon$.
\end{itemize}

\subsection{Preliminary Results} \label{pre-res}

We first give some important results which
are needed in the proofs of the main theorems.
We start with two lemmas which are consequences
of~\cite{ab}, of standard bounds on the probability
of events corresponding to having a certain number of
monochromatic crossings of an annulus (see Lemma~5
of~\cite{ksz} and Appendix~A of~\cite{lsw5}),
but which do \emph{not} depend on statement (S).

\begin{lemma} \label{sub-conv}
Let $\gamma^{\delta}_{{\mathbb D},-i,i}$ be the percolation
exploration path on the edges of $\delta {\cal H}$ inside
(the $\delta$-approximation of) $\mathbb D$ between
(the e-vertices closest to) $-i$ and $i$.
For any fixed point $z \in {\mathbb D}$, chosen independently
of $\gamma^{\delta}_{{\mathbb D},-i,i}$, as $\delta \to 0$,
$\gamma^{\delta}_{{\mathbb D},-i,i}$ and the boundary
$\partial {\mathbb D}^{\delta}_{-i,i}(z)$ of the domain
${\mathbb D}^{\delta}_{-i,i}(z)$ that contains $z$ jointly
have limits in distribution along subsequences of $\delta$
with respect to the uniform metric~(\ref{distance}) on
continuous curves.
Moreover, any subsequence limit of
$\partial {\mathbb D}^{\delta}_{-i,i}(z)$ is almost surely
a simple loop~\cite{ada}.
\end{lemma}


\begin{lemma} \label{boundaries}
Using the notation of Lemma~\ref{sub-conv},
let $\gamma_{{\mathbb D},-i,i}$ be the limit in distribution
of $\gamma^{\delta}_{{\mathbb D},-i,i}$ as $\delta \to 0$
along some convergent subsequence $\{ \delta_k \}$ and
$\partial {\mathbb D}_{-i,i}(z)$ the boundary of the domain
${\mathbb D}_{-i,i}(z)$ of ${\mathbb D} \setminus \gamma_{D,-i,i}[0,1]$
that contains $z$.
Then, as $k \to \infty$,
$(\gamma^{\delta_k}_{{\mathbb D},-i,i},\partial {\mathbb D}^{\delta_k}_{-i,i}(z))$
converges in distribution to
$(\gamma_{{\mathbb D},-i,i},\partial {\mathbb D}_{-i,i}(z))$.
\end{lemma}


The two lemmas above are important ingredients in the
proof of Theorem~\ref{thm-convergence} below.
The second one says that, for every subsequence limit, the discrete
boundaries converge to the boundaries of the domains generated by
the limiting continuous curve.
If we use statement (S), then the limit $\gamma_{{\mathbb D},-i,i}$
of $\gamma^{\delta_k}_{{\mathbb D},-i,i}$ is the trace of chordal
$SLE_6$ for every subsequence $\delta_k \downarrow 0$, and we can
use Lemmas~\ref{boundaries} and \ref{sub-conv} to deduce that all
the domains produced in the continuum construction are Jordan domains.
The key step in that direction is represented by the following
result, our proof of which relies on (S).

\begin{corollary} \label{jordan}
For any deterministic $z \in {\mathbb D}$,
the boundary $\partial {\mathbb D}_{-i,i}(z)$ of a domain
${\mathbb D}_{-i,i}(z)$ of the continuum construction is
almost surely a Jordan curve.
\end{corollary}

\noindent The corollary says that the domains that appear after
the first step of the continuum construction are Jordan domains.
The steps in the second stage of the continuum construction consist
of $SLE_6$ paths inside Jordan domains, and therefore Corollary~\ref{jordan},
combined with Riemann's mapping theorem and the conformal invariance
of $SLE_6$, implies that the domains produced during the second stage
are also Jordan.
By induction, we deduce that all the domains produced in the
continuum construction are Jordan domains.

We end this section with one more lemma which is another
key ingredient in the proof of Theorem~\ref{thm-convergence};
we remark that its proof requires (S) in a fundamental way.

\begin{lemma} \label{strong-smirnov}
Let $(D,a,b)$ denote a \emph{random} Jordan domain, with $a,b$ two
points on $\partial D$.
Let $\{ (D_k,a_k,b_k) \}_{k \in {\mathbb N}}, \, a_k,b_k \in \partial D_k$,
be a sequence of \emph{random} Jordan domains with points on their
boundaries such that, as $k \to \infty$, $(\partial D_k,a_k,b_k)$
converges in distribution to $(\partial D,a,b)$ with respect to the
uniform metric~(\ref{distance}) on continuous curves, and the Euclidean
metric on $(a,b)$.
For any sequence $\{ \delta_k \}_{k \in {\mathbb N}}$
with $\delta_k \downarrow 0$ as $k \to \infty$,
$\gamma^{\delta_k}_{D_k,a_k,b_k}$ converges in distribution to
$\gamma_{D,a,b}$ with respect to the uniform metric~(\ref{distance})
on continuous curves.
\end{lemma}

\subsection{The Main Theorems} \label{main-thms}

In this section we state the main theorems of this paper and a
corollary, our most important result, that the Continuum Nonsimple
Loop process is the scaling limit of the set of all cluster
boundaries for critical site percolation on the triangular lattice.
The corollary is obtained by combining the first two theorems.
The proofs of all these results rely on statement (S).
As noted before, statement (S) is proved in
Appendix~\ref{convergence-to-sle}.

\begin{theorem} \label{thm-convergence} For
any $k \in {\mathbb N}$, the
first $k$ steps of (a suitably reorganized version of)
the full discrete construction inside the unit disc (of
Section~\ref{full}) converge, jointly in distribution,
to the first $k$ steps of the full continuum construction
inside the unit disc (of Section~\ref{unit-disc}).
Furthermore, the scaling limit of the full (original or
reorganized) discrete construction is the full continuum
construction.

Moreover, if for any fixed $\varepsilon>0$ we let
$K_{\delta}(\varepsilon)$ denote the number of steps
needed to find all the cluster boundaries of Euclidean
diameter larger than $\varepsilon$ in the discrete
construction, then $K_{\delta}(\varepsilon)$ is
bounded in probability as $\delta \to 0$; i.e.,
$\lim_{C \to \infty} \limsup_{\delta \to 0}
{\mathbb P}(K_{\delta}(\varepsilon) > C) = 0$.
This is so in both the original and reorganized
versions of the discrete construction.
\end{theorem}


The second part of Theorem~\ref{thm-convergence}
means that both versions of the discrete
construction used in the theorem find all large
contours in a number of steps which does not diverge
as $\delta \to 0$.
This, together with the first part of the same theorem,
implies that the continuum construction does indeed
describe all macroscopic contours contained inside
the unit disc (with blue boundary conditions) as
$\delta \to 0$.

The construction presented in Section~\ref{unit-disc}
can of course be repeated
for the disc ${\mathbb D}_R$ of radius $R$, for any $R$,
so we should take a ``thermodynamic limit'' by letting
$R \to \infty$.
In this way, we would eliminate the boundary (and
the boundary conditions) and obtain a process on
the whole plane.
Such an extension from the unit disc to the plane
is contained in the next theorem.

Let $P_R$ be the (limiting) distribution of the set of
curves (all continuum nonsimple loops) generated by the
continuum construction inside ${\mathbb D}_R$ (i.e., the
limiting measure, defined by the inductive construction,
on the complete separable metric space ${\Omega}_R$ of
collections of continuous curves in ${\mathbb D}_R$).

For a domain $D$, we denote by $I_D$ the mapping (on $\Omega$
or $\Omega_R$) in which all portions of curves that exit $D$
are removed.
When applied to a configuration of loops in the plane,
$I_D$ gives a set of curves which either start and end
at points on $\partial D$ or form closed loops completely
contained in $D$.
Let $\hat I_D$ be the same mapping lifted to the space
of probability measures on $\Omega$ or $\Omega_R$.

\begin{theorem} \label{thm-therm-lim}
Theorem~\ref{thm-convergence} implies that there
exists a unique probability measure $P$ on the space $\Omega$
of collections of continuous curves in $\dot{\mathbb R}^2$ such
that $P_R \to P$ as $R \to \infty$ in the sense that for every
bounded domain $D$, as $R \to \infty$, $\hat I_D P_R \to \hat I_D P$.
\end{theorem}

\begin{remark}
We remark that we will generally take monochromatic blue
boundary conditions on the disc ${\mathbb D}_R$ of radius $R$.
But one could also take monochromatic boundary conditions
with color depending on $R$ or even non-monochromatic boundary
conditions without any essential change in the results or
the proofs.
\end{remark}

\begin{corollary} \label{scal-lim}
The Continuum Nonsimple Loop process
$P$ in the plane defined in Theorem~\ref{thm-therm-lim} is
the scaling limit of the collection of all boundary contours
for critical site percolation on the triangular lattice.
\end{corollary}

The next theorem states some properties of
the Continuum Nonsimple Loop process in the plane.

\begin{theorem} \label{features}
The Continuum Nonsimple Loop
process in the plane has the following properties,
the first three of which are valid with probability one:
\begin{enumerate}
\item The Continuum Nonsimple Loop process is a random
collection of noncrossing continuous loops in the plane.
The loops can and do touch themselves and each other
many times, but there are no triple points; i.e. no three
or more loops can come together at the same point, and a
single loop cannot touch the same point more than twice,
nor can a loop touch a point where another loop touches
itself.
\item Any deterministic point (i.e., chosen independently
of the loop process) of the plane is surrounded by an
infinite family of nested loops with diameters going to
both zero and infinity; any annulus about that point with
inner radius $r_1 > 0$ and outer radius $r_2 < \infty$
contains only a finite number of those loops.
Consequently, any two distinct deterministic points of the
plane are separated by loops winding around each of them.
\item Any two loops are connected by a finite ``path''
of touching loops.
\item[4.] The Continuum Nonsimple Loop process is conformally
invariant in the sense that, given a Jordan domain $D$
and a conformal homeomorphism $f : D \to D'$ onto $D'$,
the scaling limits, $P_D$ and $P_{D'}$, of the loops
inside $D$ and $D'$ taken with, say, blue boundary
conditions are related by $f * P_D = P_{D'}$.
(Here $f * P_D$ denotes the probability distribution of
the loop process $f(X)$ when $X$ is distributed by $P_D$.)
\end{enumerate}
\end{theorem}

To conclude this section, we show how to recover chordal $SLE_6$
from the Continuum Nonsimple Loop process, i.e., given a (deterministic)
Jordan domain $D$ with two boundary points $a$ and $b$, we give a
construction that uses the continuum nonsimple loops of $P$ to generate
a process distributed like chordal $SLE_6$ inside $D$ from $a$ to $b$.

Remember, first of all, that each continuum nonsimple loop has
either a clockwise or counterclockwise direction, with the set
of all loops surrounding any deterministic point alternating
in direction.
For convenience, let us suppose that $a$ is at the ``bottom"
and $b$ is at the ``top" of $D$ so that the boundary is divided
into a left and right part by these two points.
Fix $\varepsilon>0$ and call $LR(\varepsilon)$ the set of all
the directed segments of loops that connect from the left to
the right part of the boundary touching $\partial D$ at a distance
larger than $\varepsilon$ from both $a$ and $b$, and $RL(\varepsilon)$
the analogous set of directed segments from the right to the left
portion of $\partial D$.
For a fixed $\varepsilon>0$, there is only a finite number
of such segments, and, if they are ordered moving along the
left boundary of $D$ from $a$ to $b$, they alternate in direction
(i.e., a segment in $LR(\varepsilon)$ is followed by one in
$RL(\varepsilon)$ and so on).

Between a segment in $RL(\varepsilon)$ and the next segment
in $LR(\varepsilon)$, there are countably many portions of
loops intersecting $D$ which start and end on $\partial D$
and are maximal in the sense that they are not contained
inside any other portion of loop of the same type; they all
have counterclockwise direction and can be used to make a
``bridge'' between the right-to-left segment and the next
one (in $LR(\varepsilon)$).
This is done by pasting the portions of loops together with
the help of points in $\partial D$ and a limit procedure
to produce a connected (nonsimple) path.

If we do this for each pair of successive segments on both
sides of the boundary of $D$, we get a path that connects
two points on $\partial D$.
By letting $\varepsilon \to 0$ and taking the limit of this
procedure, since almost surely $a$ and $b$ are surrounded
by an infinite family of nested loops with diameters going
to zero, we obtain a path that connects $a$ with $b$;
this path is distributed as chordal $SLE_6$ inside $D$
from $a$ to $b$.
The last claim follows from considering the analogous
procedure for percolation on the discrete lattice
$\delta {\cal H}$, using segments of boundaries.
It is easy to see that in the discrete case this
procedure produces exactly the same path as the
percolation exploration process.
By Corollary~\ref{scal-lim}, the scaling limit of this
discrete procedure is the continuum one described above,
therefore the claim follows from (S).

%
%

\section{Proofs} \label{proofs}

In this section we present the proofs of the results
stated in Section~\ref{main-results}.

\bigskip

\noindent {\bf Proof of Lemma~\ref{sub-conv}.}
The first part of the lemma is a direct consequence
of~\cite{ab}; it is enough to notice that the (random)
polygonal curves $\gamma^{\delta}_{{\mathbb D},-i,i}$
and $\partial {\mathbb D}^{\delta}_{-i,i}(z)$
satisfy the conditions in~\cite{ab} and thus have a
scaling limit in terms of continuous curves, at least
along subsequences of $\delta$.

To prove the second part, we use standard percolation
bounds (see Lemma~5 of~\cite{ksz} and Appendix A
of~\cite{lsw5}) to show that, in the limit $\delta \to 0$,
the loop $\partial {\mathbb D}^{\delta}_{-i,i}(z)$ does
not collapse on itself but remains a  simple loop.

Let us assume that this is not the case and that
the limit $\tilde\gamma$ of
$\partial {\mathbb D}^{\delta_k}_{-i,i}(z)$ along some
subsequence $\{ \delta_k \}_{k \in {\mathbb N}}$ touches
itself, i.e., $\tilde\gamma(t_0)=\tilde\gamma(t_1)$ for
$t_0 \neq t_1$ with positive probability.
If that happens, we can take $\varepsilon>\varepsilon'>0$
small enough so that the annulus
$B(\tilde\gamma(t_1),\varepsilon) \setminus B(\tilde\gamma(t_1),\varepsilon')$
is crossed at least four times by $\tilde\gamma$
(here $B(u,r)$ is the ball of radius $r$ centered at $u$).

Because of the choice of topology, the convergence in
distribution of $\partial {\mathbb D}^{\delta_k}_{-i,i}(z)$
to $\tilde\gamma$ implies that we can find coupled versions of
$\partial {\mathbb D}^{\delta_k}_{-i,i}(z)$ and $\tilde\gamma$
on some probability space $(\Omega',{\cal B}',{\mathbb P}')$ such that
$\text{d}(\partial {\mathbb D}^{\delta}_{-i,i}(z),\tilde\gamma) \to 0$,
for all $\omega' \in \Omega'$ as $k \to \infty$ (see, for example,
Corollary~1 of~\cite{billingsley1}).

Using this coupling, we can choose $k$ large enough (depending
on $\omega'$) so that $\partial {\mathbb D}^{\delta_k}_{-i,i}(z)$
stays in an $\varepsilon'/2$-neighborhood
${\cal N}(\tilde\gamma,\varepsilon'/2) \equiv \bigcup_{u \in \tilde\gamma} B(u,\varepsilon'/2)$
of $\tilde\gamma$.
This event however would correspond to (at least) four
paths of one color (corresponding to the four crossings
by $\partial {\mathbb D}^{\delta_k}_{-i,i}(z)$) and two
of the other color of the annulus
$B(\tilde\gamma(t_1),\varepsilon-\varepsilon'/2) \setminus
B(\tilde\gamma(t_1),3 \, \varepsilon'/2)$.
As $\delta_k \to 0$, we can let $\varepsilon' \to 0$,
in which case the probability of seeing the event
just described somewhere inside $\mathbb D$ goes to
zero~\cite{ksz,lsw5}, leading to a contradiction. \fbox{} \\

\noindent {\bf Proof of Lemma~\ref{boundaries}.}
Let $\{ \delta_k \}_{k \in {\mathbb N}}$ be a convergent
subsequence for $\gamma^{\delta}_{{\mathbb D},-i,i}$
and $\gamma \equiv \gamma_{{\mathbb D},-i,i}$
the limit in distribution of $\gamma^{\delta_k}_{{\mathbb D},-i,i}$
as $k \to \infty$.
For simplicity of notation, in the rest of the proof we will drop
the $k$ and write $\delta$ instead of $\delta_k$.
Because of the choice of topology, the convergence in distribution
of $\gamma^{\delta} \equiv \gamma^{\delta}_{{\mathbb D},-i,i}$ to
$\gamma$ implies that we can find coupled versions of $\gamma^{\delta}$
and $\gamma$ on some probability space $(\Omega',{\cal B}',{\mathbb P}')$
such that $\text{d}(\gamma^{\delta}(\omega'),\gamma(\omega')) \to 0$,
for all $\omega'$ as $k \to \infty$ (see, for example, Corollary~1
of~\cite{billingsley1}).
Using this coupling, our first task will be to prove the following claim:
\begin{itemize}
\item[(C)] For two (deterministic) points $u,v \in {\mathbb D}$,
the probability that ${\mathbb D}_{-i,i}(u) = {\mathbb D}_{-i,i}(v)$ but
${\mathbb D}^{\delta}_{-i,i}(u) \neq {\mathbb D}^{\delta}_{-i,i}(v)$
or vice versa goes to zero as $\delta \to 0$.
\end{itemize}

Let us consider first the case of $u,v$ such that
${\mathbb D}_{-i,i}(u) = {\mathbb D}_{-i,i}(v)$ but
${\mathbb D}^{\delta}_{-i,i}(u) \neq {\mathbb D}^{\delta}_{-i,i}(v)$.
Since ${\mathbb D}_{-i,i}(u)$ is an open subset of $\mathbb C$,
there exists a continuous curve $\gamma_{u,v}$ joining
$u$ and $v$ and a constant $\varepsilon>0$ such that
the $\varepsilon$-neighborhood
${\cal N}(\gamma_{u,v},\varepsilon)$ of the curve is contained
in ${\mathbb D}_{-i,i}(u)$, which implies that $\gamma$ does
not intersect ${\cal N}(\gamma_{u,v},\varepsilon)$.
Now, if $\gamma^{\delta}$ does not intersect
${\cal N}(\gamma_{u,v},\varepsilon/2)$, for
$\delta$ small enough, then there is a $\cal T$-path
$\pi$ of unexplored hexagons connecting the hexagon
that contains $u$ with the hexagon that contains $v$,
and we conclude that
${\mathbb D}^{\delta}_{-i,i}(u) = {\mathbb D}^{\delta}_{-i,i}(v)$.

This shows that the event that ${\mathbb D}_{-i,i}(u) = {\mathbb D}_{-i,i}(v)$
but ${\mathbb D}^{\delta}_{-i,i}(u) \neq {\mathbb D}^{\delta}_{-i,i}(v)$
implies the existence of a curve $\gamma_{u,v}$ whose
$\varepsilon$-neighborhood ${\cal N}(\gamma_{u,v},\varepsilon)$
is not intersected by $\gamma$ but whose $\varepsilon/2$-neighborhood
${\cal N}(\gamma_{u,v},\varepsilon/2)$ is intersected by $\gamma^{\delta}$.
This implies that $\forall u,v \in {\mathbb D}$, $\exists \varepsilon>0$ such
that ${\mathbb P}'({\mathbb D}_{-i,i}(u) = {\mathbb D}_{-i,i}(v) \text{ but }
{\mathbb D}^{\delta}_{-i,i}(u) \neq {\mathbb D}^{\delta}_{-i,i}(v)) \leq
{\mathbb P}'(\text{d}(\gamma^{\delta},\gamma) \geq \varepsilon/2)$.
But the right hand side goes to zero for every $\varepsilon>0$ as
$\delta \to 0$, which concludes the proof of one direction of the claim.

To prove the other direction, we consider two points
$u,v \in D$ such that $D_{-i,i}(u) \neq D_{-i,i}(v)$
but $D^{\delta}_{-i,i}(u) = D^{\delta}_{-i,i}(v)$.
Assume that $u$ is trapped before $v$ by $\gamma$
and suppose for the moment that ${\mathbb D}_{-i,i}(u)$
is a domain of type 3 or 4; the case of a domain of
type 1 or 2 is analogous and will be treated later.
Let $t_1$ be the first time $u$ is trapped by
$\gamma$ with $\gamma(t_0)=\gamma(t_1)$ the double
point of $\gamma$ where the domain ${\mathbb D}_{-i,i}(u)$
containing $u$ is ``sealed off.''
At time $t_1$, a new domain containing $u$ is
created and $v$ is disconnected from $u$.

Choose $\varepsilon>0$ small enough so that neither $u$ nor
$v$ is contained in the ball $B(\gamma(t_1),\varepsilon)$
of radius $\varepsilon$ centered at $\gamma(t_1)$, nor in
the $\varepsilon$-neighborhood
${\cal N}(\gamma[t_0,t_1],\varepsilon)$ of the portion
of $\gamma$ which surrounds $u$. 
Then it follows from the coupling that, for $\delta$
small enough, there are appropriate parameterizations
of $\gamma$ and $\gamma^{\delta}$ such that the portion
$\gamma^{\delta}[t_0,t_1]$ of $\gamma^{\delta}(t)$ is
inside ${\cal N}(\gamma[t_0,t_1],\varepsilon)$, and
$\gamma^{\delta}(t_0)$ and $\gamma^{\delta}(t_1)$ are
contained in $B(\gamma(t_1),\varepsilon)$.

For $u$ and $v$ to be contained in the same domain
in the discrete construction, there must be a
$\cal T$-path $\pi$ of unexplored hexagons
connecting the hexagon that contains $u$ to the
hexagon that contains $v$.
From what we said in the previous paragraph, any
such $\cal T$-path connecting $u$ and $v$
would have to go though a ``bottleneck'' in
$B(\gamma(t_1),\varepsilon)$.

Assume now, for concreteness but without loss of
generality, that ${\mathbb D}_{-i,i}(u)$ is a domain
of type~3, which means that $\gamma$ winds around
$u$ counterclockwise, and consider the hexagons
to the ``left" of $\gamma^{\delta}[t_0,t_1]$.
Those hexagons form a ``quasi-loop'' around $u$
since they wind around it (counterclockwise) and
the first and last hexagons are both contained in
$B(\gamma(t_1),\varepsilon)$.
The hexagons to the left of $\gamma^{\delta}[t_0,t_1]$
belong to the set $\Gamma_Y(\gamma^{\delta})$, which
can be seen as a (nonsimple) path by connecting the
centers of the hexagons in $\Gamma_Y(\gamma^{\delta})$
by straight segments.
Such a path shadows $\gamma^{\delta}$, with the difference
that it can have double (or even triple) points, since
the same hexagon can be visited more than once.
Consider $\Gamma_Y(\gamma^{\delta})$ as a path
$\hat\gamma^{\delta}$ with a given parametrization
$\hat\gamma^{\delta}(t)$, chosen so that
$\hat\gamma^{\delta}(t)$ is inside
$B(\gamma(t_1),\varepsilon)$ when $\gamma^{\delta}(t)$ is,
and it winds around $u$ together with $\gamma^{\delta}(t)$.

Now suppose that there were two times, $\hat t_0$ and
$\hat t_1$, such that $\hat\gamma^{\delta}(\hat t_1)
= \hat\gamma^{\delta}(\hat t_0) \in B(\gamma(t_1),\varepsilon)$
and $\hat\gamma^{\delta}[\hat t_0,\hat t_1]$ winds
around $u$.
This would imply that the ``quasi-loop'' of explored
yellow hexagons around $u$ is actually completed, and
that $D^{\delta}_{a,b}(v) \neq D^{\delta}_{a,b}(u)$.
Thus, for $u$ and $v$ to belong to the same discrete
domain, this cannot happen.

For any $0<\varepsilon'<\varepsilon$, if we take $\delta$
small enough, $\hat\gamma^{\delta}$ will be contained
inside ${\cal N}(\gamma,\varepsilon')$, due to the coupling.
Following the considerations above, the fact that $u$
and $v$ belong to the same domain in the discrete
construction but to different domains in the continuum
construction implies, for $\delta$ small enough, that
there are four disjoint yellow $\cal T$-paths
crossing the annulus $B(\gamma(t_1),\varepsilon)
\setminus B(\gamma(t_1),\varepsilon')$ (the paths
have to be disjoint because, as we said,
$\hat\gamma^{\delta}$ cannot, when coming back to
$B(\gamma(t_1),\varepsilon)$ after winding around
$u$, touch itself inside $B(\gamma(t_1),\varepsilon)$).
Since $B(\gamma(t_1),\varepsilon) \setminus
B(\gamma(t_1),\varepsilon')$ is also crossed
by at least two blue $\cal T$-paths from
$\Gamma_B(\gamma^{\delta})$, there is a total
of at least six $\cal T$-paths, not all of
the same color, crossing the annulus
$B(\gamma(t_1),\varepsilon) \setminus B(\gamma(t_1),\varepsilon')$.

Let us call ${\cal A}_w(\varepsilon,\varepsilon')$ the
event described above, where $\gamma(t_1)=w$; a standard
bound~\cite{ksz} on the probability of six disjoint crossings
(not all of the same color) of an annulus gives that the
probability of ${\cal A}_w(\varepsilon,\varepsilon')$ scales as
$(\frac{\varepsilon'}{\varepsilon})^{2+\alpha}$ with $\alpha>0$.
As $\delta \to 0$, we can let $\varepsilon'$ go to
zero (keeping $\varepsilon$ fixed); when we do this,
the probability of ${\cal A}_w(\varepsilon,\varepsilon')$
goes to zero sufficiently rapidly with $\varepsilon'$
to conclude, like in the proof of Lemma~\ref{sub-conv},
that the probability to see such an event anywhere in
$\mathbb D$ goes to zero.


In the case in which $u$ belongs to a domain of type~1
or 2, let $\cal E$ be the excursion that traps $u$
and $\gamma(t_0) \in \partial {\mathbb D}$ be the point
on the boundary of $\mathbb D$ where $\cal E$ starts and
$\gamma(t_1) \in \partial {\mathbb D}$ the point where
it ends.
Choose $\varepsilon>0$ small enough so that neither $u$
nor $v$ is contained in the balls
$B(\gamma(t_0),\varepsilon)$ and $B(\gamma(t_1),\varepsilon)$
of radius $\varepsilon$ centered at $\gamma(t_0)$ and
$\gamma(t_1)$, nor in the $\varepsilon$-neighborhood
${\cal N}({\cal E},\varepsilon)$ of the excursion $\cal E$.
Because of the coupling, for $\delta$ small enough
(depending on $\varepsilon$), $\gamma^{\delta}$ shadows
$\gamma$ along $\cal E$, staying within
${\cal N}({\cal E},\varepsilon)$.
If this is the case, any $\cal T$-path of unexplored
hexagons connecting the hexagon that contains $u$ with
the hexagon that contains $v$ would have to go through
one of two ``bottlenecks,'' one contained in
$B(\gamma(t_0),\varepsilon)$ and the other in
$B(\gamma(t_1),\varepsilon)$.

Assume for concreteness (but without loss
of generality) that $u$ is in a domain of type 1,
which means that $\gamma$ winds around $u$
counterclockwise.
If we parameterize $\gamma$ and $\gamma^{\delta}$ so that
$\gamma^{\delta}(t_0) \in B(\gamma(t_0),\varepsilon)$
and $\gamma^{\delta}(t_1) \in B(\gamma(t_1),\varepsilon)$,
$\gamma^{\delta}[t_0,t_1]$ forms a ``quasi-excursion''
around $u$ since it winds around it (counterclockwise)
and it starts inside $B_{\varepsilon}(\gamma(t_0))$ and
ends inside $B_{\varepsilon}(\gamma(t_1))$.
Notice that if $\gamma^{\delta}$ touched
$\partial {\mathbb D}^{\delta}$, inside both
$B_{\varepsilon}(\gamma(t_0))$ and
$B_{\varepsilon}(\gamma(t_1))$, this would imply
that the ``quasi-excursion'' is a real excursion and
that $D^{\delta}_{a,b}(v) \neq D^{\delta}_{a,b}(u)$.

For any $0<\varepsilon'<\varepsilon$, if we take $\delta$
small enough, $\gamma^{\delta}$ will be contained inside
${\cal N}(\gamma,\varepsilon')$, due to the coupling.
Therefore, the fact that
${\mathbb D}^{\delta}_{a,b}(v) = {\mathbb D}^{\delta}_{a,b}(u)$
implies, with probability going to one as
$\delta \to 0$, that for $\varepsilon>0$ fixed
and any $0<\varepsilon'<\varepsilon$, $\gamma^{\delta}$
enters the ball $B(\gamma(t_i),\varepsilon')$
and does not touch $\partial {\mathbb D}^{\delta}$
inside the larger ball $B(\gamma(t_i),\varepsilon)$,
for $i=0$ or $1$.
This is equivalent to having at least two yellow and one
blue $\cal T$-paths (contained in ${\mathbb D}^{\delta}$)
crossing the annulus
$B(\gamma(t_i),\varepsilon) \setminus B(\gamma(t_i),\varepsilon')$.
Let us call ${\cal B}_w(\varepsilon,\varepsilon')$ the
event described above, where $\gamma(t_i)=w$; a standard
bound~\cite{lsw5} (this bound can also be derived from the
one obtained in~\cite{ksz}) on the probability of disjoint
crossings (not all of the same color) of a semi-annulus in
the upper half-plane gives that the probability of
${\cal B}_w(\varepsilon,\varepsilon')$ scales as
$(\frac{\varepsilon'}{\varepsilon})^{1+\beta}$ with $\beta>0$.
(We can apply the bound to our case because the unit disc
is a convex subset of the half-plane $\{ x+iy:y>-1 \}$
and therefore the intersection of an annulus centered at
say $-i$ with the unit disc is a subset of the intersection
of the same annulus with the half-plane $\{ x+iy:y>-1 \}$.)
As $\delta \to 0$, we can let $\varepsilon'$ go to zero
(keeping $\varepsilon$ fixed), concluding that the probability
that such an event occurs anywhere on the boundary of the disc
goes to zero.

We have shown that, for two fixed points $u,v \in {\mathbb D}$,
having ${\mathbb D}_{-i,i}(u) \neq {\mathbb D}_{-i,i}(v)$ but
${\mathbb D}^{\delta}_{-i,i}(u) = {\mathbb D}^{\delta}_{-i,i}(v)$
or vice versa implies the occurrence of an event whose
probability goes to zero as $\delta \to 0$, and the
proof of the claim is concluded.

We now introduce the Hausdorff distance $\text{d}_{\text{H}}(A,B)$
between two closed nonempty subsets of $\overline{\mathbb D}$:
\begin{equation} \label{hausdorff-dist}
\text{d}_{\text{H}}(A,B) \equiv \inf \{ \ell \geq 0 : B \subset
\cup_{a \in A} B(a,\ell), \, A \subset
\cup_{b \in B} B(b,\ell) \}.
\end{equation}
With this metric, the collection of closed subsets of
$\overline{\mathbb D}$ is a compact space.
We will next prove that
$\partial {\mathbb D}^{\delta}_{-i,i}(z)$ converges
in distribution to $\partial {\mathbb D}_{-i,i}(z)$
as $\delta \to 0$, in the topology induced
by~(\ref{hausdorff-dist}).
(Notice that the coupling between $\gamma^{\delta}$
and $\gamma$ provides a coupling between
$\partial {\mathbb D}^{\delta}_{-i,i}(z)$ and
$\partial {\mathbb D}_{-i,i}(z)$, seen as
boundaries of domains produced by the two paths.)


We will now use Lemma~\ref{sub-conv} and take a further subsequence
$k_n$ of the $\delta$'s that for simplicity of notation we denote
by $\{ \delta_n \}_{n \in {\mathbb N}}$ such that, as $n \to \infty$,
$\{ \gamma^{\delta_n},\partial {\mathbb D}^{\delta_n}_{-i,i}(z) \}$
converge jointly in distribution to $\{ \gamma,\tilde\gamma \}$,
where $\tilde\gamma$ is a simple loop.
For any $\varepsilon>0$, since $\tilde\gamma$ is a compact set,
we can find a covering of $\tilde \gamma$ by a finite number
of balls of radius $\varepsilon/2$ centered at points on $\tilde\gamma$.
Each ball contains both points in the interior
$\text{int}(\tilde\gamma)$ of $\tilde\gamma$ and in the exterior
$\text{ext}(\tilde\gamma)$ of $\tilde \gamma$, and we can choose
(independently of $n$) one point from $\text{int}(\tilde\gamma)$
and one from $\text{ext}(\tilde\gamma)$ inside each ball.

Once again, the convergence in distribution of
$\partial {\mathbb D}^{\delta_n}_{-i,i}(z)$
to $\tilde\gamma$ implies the existence of a coupling
such that, for $n$ large enough, the selected points
that are in $\text{int}(\tilde\gamma)$ are contained
in ${\mathbb D}^{\delta_n}_{-i,i}(z)$, and those that
are in $\text{ext}(\tilde\gamma)$ are contained in the
complement of $\overline{{\mathbb D}^{\delta_n}_{-i,i}(z)}$.
But by claim (C), each one of the selected points that is contained
in ${\mathbb D}^{\delta_n}_{-i,i}(z)$ is also contained in
${\mathbb D}_{-i,i}(z)$ with probability going to $1$ as $n \to \infty$;
analogously, each one of the selected points contained in the complement
of $\overline{{\mathbb D}^{\delta_n}_{-i,i}(z)}$ is also contained in the
complement of $\overline{{\mathbb D}_{-i,i}(z)}$ with probability
going to $1$ as $n \to \infty$.
This implies that $\partial {\mathbb D}_{-i,i}(z)$ crosses each
one of the balls in the covering of $\tilde\gamma$, and therefore
$\tilde\gamma \subset \cup_ {u \in \partial {\mathbb D}_{-i,i}(z)} B(u,\varepsilon)$.
From this and the coupling between
$\partial {\mathbb D}^{\delta_n}_{-i,i}(z)$ and $\tilde\gamma$,
it follows immediately that, for $n$ large enough,
$\partial {\mathbb D}^{\delta_n}_{-i,i}(z) \subset
\cup_ {u \in \partial {\mathbb D}_{-i,i}(z)} B(u,\varepsilon)$
with probability close to one.

A similar argument (analogous to the previous one but simpler,
since it does not require the use of $\tilde\gamma$), with the roles
of ${\mathbb D}^{\delta_n}_{-i,i}(z)$ and ${\mathbb D}_{-i,i}(z)$
inverted, shows that $\partial {\mathbb D}_{-i,i}(z) \subset
\cup_ {u \in \partial {\mathbb D}^{\delta_n}_{-i,i}(z)} B(u,\varepsilon)$
with probability going to $1$ as $n \to \infty$.
Therefore, for all $\varepsilon>0$,
${\mathbb P}(\text{d}_{\text{H}}(\partial {\mathbb D}^{\delta_n}_{-i,i}(z),
\partial {\mathbb D}_{-i,i}(z)) > \varepsilon) \to 0$ as $n \to \infty$,
which implies convergence in distribution of
$\partial {\mathbb D}^{\delta_n}_{-i,i}(z)$ to
$\partial {\mathbb D}_{-i,i}(z)$, as $\delta_n \to 0$,
in the topology induced by~(\ref{hausdorff-dist}).
But Lemma~\ref{sub-conv} implies that
$\partial {\mathbb D}^{\delta_n}_{-i,i}(z)$
converges in distribution (using~(\ref{distance})) to a
simple loop, therefore $\partial {\mathbb D}_{-i,i}(z)$
must also be a simple loop; and we have convergence in
the topology induced by~(\ref{distance}).

It is also clear that the argument above is independent
of the subsequence $\{ \delta_n \}$, so the limit of
$\partial {\mathbb D}^{\delta}_{-i,i}(z)$ is unique and
coincides with $\partial {\mathbb D}_{-i,i}(z)$.
Hence, we have convergence in distribution of
$\partial {\mathbb D}^{\delta}_{-i,i}(z)$ to
$\partial {\mathbb D}_{-i,i}(z)$, as $\delta \to 0$,
in the topology induced by~(\ref{distance}), and
indeed joint convergence of
$(\gamma^{\delta},\partial {\mathbb D}^{\delta}_{-i,i}(z))$
to $(\gamma,\partial {\mathbb D}_{-i,i}(z))$. \fbox{} \\

\noindent {\bf Proof of Corollary~\ref{jordan}.}
The corollary follows immediately from Lemma~\ref{sub-conv}
and Lemma~\ref{boundaries}, as already seen in the proof
of Lemma~\ref{boundaries}. \fbox {} \\

\noindent {\bf Proof of Lemma~\ref{strong-smirnov}.}
First of all recall that the convergence of $(\partial D_k,a_k,b_k)$
to $(\partial D,a,b)$ in distribution implies the existence of
coupled versions of $(\partial D_k,a_k,b_k)$ and $(\partial D,a,b)$
on some probability space $(\Omega',{\cal B}',{\mathbb P}')$ such that
$\text{d}(\partial D(\omega'),\partial D_k(\omega')) \to 0$,
$a_k(\omega') \to a(\omega')$, $b_k(\omega') \to b(\omega')$
for all $\omega'$ as $k \to \infty$ (see, for example,
Corollary~1 of~\cite{billingsley1}).
This immediately implies that the conditions to apply Rad\'o's theorem
(see Theorem~\ref{rado-thm} of Appendix~\ref{rado}) are satisfied.
%
%
Let $f_k$ be the conformal map that takes the unit disc $\mathbb D$
onto $D_k$ with $f_k(0)=0$ and $f'_k(0)>0$, and let $f$ be the
conformal map from $\mathbb D$ onto $D$ with $f(0)=0$ and
$f'(0)>0$.
Then, by Theorem~\ref{rado-thm}, $f_k$ converges to $f$
uniformly in $\overline{\mathbb D}$, as $k \to \infty$.

Let $\gamma$ (resp., $\gamma_k$) be the chordal $SLE_6$
inside $D$ (resp., $D_k$) from $a$ to $b$ (resp., from $a_k$
to $b_k$), $\tilde \gamma = f^{-1}(\gamma)$, $\tilde a = f^{-1}(a)$,
$\tilde b = f^{-1}(b)$, and $\tilde \gamma_k = f_k^{-1}(\gamma_k)$,
$\tilde a_k = f_k^{-1}(a_k)$, $\tilde b_k = f_k^{-1}(b_k)$.
We note that, because of the conformal invariance of chordal
$SLE_6$, $\tilde\gamma$ (resp., $\tilde\gamma_k$) is distributed
as chordal $SLE_6$ in $\mathbb D$ from $\tilde a$ to $\tilde b$
(resp., from $\tilde a_k$ to $\tilde b_k$).
Since $|a-a_k| \to 0$ and $|b-b_k| \to 0$ for all $\omega'$,
and $f_k \to f$ uniformly in $\overline{\mathbb D}$,
we conclude that $|\tilde a - \tilde a_k| \to 0$
and $|\tilde b - \tilde b_k| \to 0$ for all $\omega'$.

Later we will prove a ``continuity'' property of $SLE_6$
(Lemma~\ref{continuity}) that allows us to conclude that,
under these conditions, $\tilde\gamma_k$ converges in
distribution to $\tilde\gamma$ in the uniform metric~(\ref{distance})
on continuous curves.
Once again, this implies the existence of coupled versions of
$\tilde\gamma_k$ and $\tilde\gamma$ on some probability space
$(\Omega', {\cal B}', {\mathbb P}')$ such that
$\text{d}(\tilde \gamma(\omega'),\tilde \gamma_k(\omega')) \to 0$,
for all $\omega'$ as $k \to \infty$.
Therefore, thanks to the convergence of $f_k$ to $f$
uniformly in $\overline{\mathbb D}$,
$\text{d}(f(\tilde\gamma(\omega')),f_k(\tilde\gamma_k(\omega'))) \to 0$,
for all $\omega'$ as $k \to \infty$.
But since $f(\tilde\gamma_k)$ is distributed as $\gamma_{D_k,a_k,b_k}$
and $f(\tilde\gamma)$ is distributed as $\gamma_{D,a,b}$,
we conclude that, as $k \to \infty$, $\gamma_{D_k,a_k,b_k}$ converges
in distribution to $\gamma_{D,a,b}$ in the uniform metric~(\ref{distance})
on continuous curves.

We now note that (S) implies that, as $\delta \to 0$,
$\gamma^{\delta}_{D_k,a_k,b_k}$ converges in distribution to
$\gamma_{D_k,a_k,b_k}$ \emph{uniformly} in $k$, for $k$ large enough.
Therefore, as $k \to \infty$, $\gamma^{\delta_k}_{D_k,a_k,b_k}$
converges in distribution to $\gamma_{D,a,b}$, and the proof is
concluded. \fbox{}

\begin{lemma} \label{continuity}
Let ${\mathbb D} \subset {\mathbb C}$ be the unit disc, $a$
and $b$ two distinct points on its boundary, and $\gamma$ the
trace of chordal $SLE_6$ inside $\mathbb D$ from $a$ to $b$.
Let $\{ a_k \}$ and $\{ b_k \}$ be two sequences of points in
$\partial {\mathbb D}$ such that $a_k \to a$ and $b_k \to b$.
Then, as $k \to \infty$, the trace $\gamma_k$ of chordal $SLE_6$
inside $\mathbb D$ from $a_k$ to $b_k$ converges in distribution
to $\gamma$ in the uniform topology~(\ref{distance}) on continuous
curves.
\end{lemma}

\noindent {\bf Proof.} Let
$f_k(z)=e^{i \alpha_k} \frac{z-z_k}{1-\bar z_k z}$
be the (unique) linear fractional transformation that takes
the unit disc $\mathbb D$ onto itself, mapping $a$ to $a_k$,
$b$ to $b_k$, and a third point $c \in \partial {\mathbb D}$
distinct from $a$ and $b$ to itself.
$\alpha_k$ and $z_k$ depend continuously on $a_k$ and $b_k$.
As $k \to \infty$, since $a_k \to a$ and $b_k \to b$,
$f_k$ converges uniformly to the identity in $\overline{\mathbb D}$.

Using the conformal invariance of chordal $SLE_6$, we couple
$\gamma_k$ and $\gamma$ by writing $\gamma_k=f_k(\gamma)$.
The uniform convergence of $f_k$ to the identity implies
that $\text{d}(\gamma,\gamma_k) \to 0$ as $k \to \infty$,
which is enough to conclude that $\gamma_k$ converges to
$\gamma$ in distribution. \fbox{} \\

\noindent {\bf Proof of Theorem~\ref{thm-convergence}.}
Let us prove the second part of the theorem first.
We will do this for the original version of the discrete
construction, but essentially the same proof works for the
reorganized version we will describe below, as we will explain
later.
Suppose that at step $k$ of this discrete construction
an exploration process $\gamma^{\delta}_k$ is run inside
a domain $D^{\delta}_{k-1}$, and write
$D^{\delta}_{k-1} \setminus \Gamma(\gamma^{\delta}_k) =
\bigcup_j D^{\delta}_{k,j}$, where $\{ D^{\delta}_{k,j} \}$
are the maximal connected domains of unexplored hexagons
into which $D^{\delta}_{k-1}$ is split by removing the
set $\Gamma(\gamma^{\delta}_k)$ of hexagons explored by
$\gamma^{\delta}_k$.

\begin{figure}[!ht]
\begin{center}
\includegraphics[width=8cm]{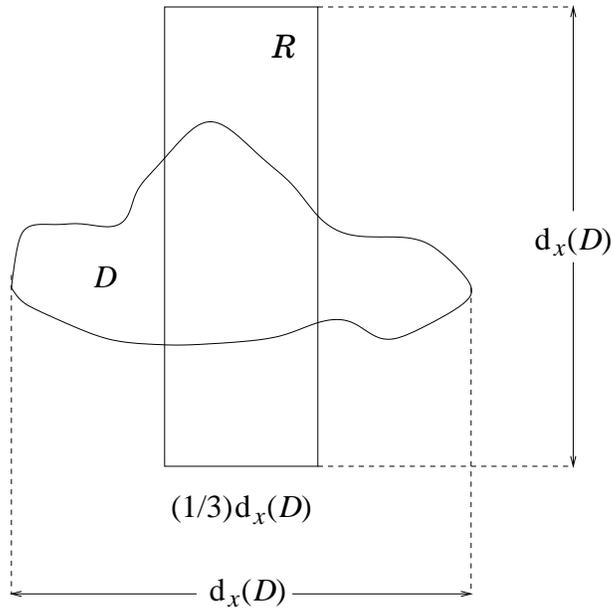}
\caption{Schematic drawing of a domain $D$ with
$\text{d}_x(D) \geq \text{d}_y(D)$ and the associated
rectangle $\cal R$.}
\label{fig3-thm1}
\end{center}
\end{figure}

Let $\text{d}_x(D^{\delta}_{k-1})$ and
$\text{d}_y(D^{\delta}_{k-1})$ be respectively the
maximal $x$- and $y$-distances between pairs of points
in $\partial D^{\delta}_{k-1}$.
Suppose, without loss of generality, that
$\text{d}_x(D^{\delta}_{k-1}) \geq \text{d}_y(D^{\delta}_{k-1})$,
and consider the rectangle $\cal R$ (see Figure~\ref{fig3-thm1})
whose vertical sides are aligned to the $y$-axis, have
length $\text{d}_x(D^{\delta}_{k-1})$, and are each
placed at $x$-distance
$\frac{1}{3} \, \text{d}_x(D^{\delta}_{k-1})$ from
points of $\partial D^{\delta}_{k-1}$ with minimal
or maximal $x$-coordinate in such a way that the
horizontal sides of $\cal R$ have length
$\frac{1}{3} \, \text{d}_x(D^{\delta}_{k-1})$;
the bottom and top sides of $\cal R$ are placed in
such a way that they are at equal $y$-distance from
the points of $\partial D^{\delta}_{k-1}$ with minimal
or maximal $y$-coordinate, respectively.

It follows from the Russo-Seymour-Welsh lemma~\cite{russo,sewe}
(see also~\cite{kesten,grimmett}) that the probability to
have two vertical $\cal T$-crossings of $\cal R$ of different
colors is bounded away from zero by a positive constant $p_0$
that does not depend on $\delta$ (for $\delta$ small enough).
If that happens, then
$\max_j \text{d}_x(D^{\delta}_{k,j}) \leq \frac{2}{3} \text{d}_x(D^{\delta}_{k-1})$.
The same argument of course applies to the maximal
$y$-distance when
$\text{d}_y(D^{\delta}_{k-1}) \geq \text{d}_x(D^{\delta}_{k-1})$.
We can summarize the above observation in the following lemma.
\begin{lemma} \label{smaller}
Suppose that at step $k$ of the full discrete construction
an exploration process $\gamma^{\delta}_k$ is run inside
a domain $D^{\delta}_{k-1}$.
If $\text{\emph{d}}_x(D^{\delta}_{k-1}) \geq \text{\emph{d}}_y(D^{\delta}_{k-1})$,
then for $\delta$ small enough (i.e., $\delta \leq C \, \text{\emph{d}}_x(D^{\delta}_{k-1})$
for some constant $C$),
$\max_j \text{\emph{d}}_x(D^{\delta}_{k,j}) \leq \frac{2}{3} \text{\emph{d}}_x(D^{\delta}_{k-1})$
with probability at least $p_0$ independent of $\delta$.
The same holds for the maximal $y$-distances when
$\text{\emph{d}}_y(D^{\delta}_{k-1}) \geq \text{\emph{d}}_x(D^{\delta}_{k-1})$.
\end{lemma}

%

Here is another lemma that will be useful later on.
\begin{lemma} \label{daughters}
Two ``daughter" subdomains, $D^{\delta}_{k,j}$ and $D^{\delta}_{k,j'}$,
either have disjoint s-boun\-da\-ries, or else their common s-boundary
consists of exactly two adjacent hexagons (of the same color) where
the exploration path $\gamma^{\delta}_k$ came within $2$ hexagons of
touching itself just when completing the s-boundary of one of the
two subdomains.
\end{lemma}

\noindent {\bf Proof.} Suppose that the two daughter subdomains
have s-boundaries $\Delta D^{\delta}_{k,j}$ and $\Delta D^{\delta}_{k,j'}$
that are not disjoint and let $S = \{ \xi_1,\ldots,\xi_i \}$ be the set
of (sites of $\cal T$ that are the centers of the) hexagons that belong
to both s-boundaries.
$S$ can be partitioned into subsets consisting of single hexagons
that are not adjacent to any another hexagon in $S$ and groups of
hexagons that form simple $\cal T$-paths (because the s-boundaries
of the two subdomains are simple $\cal T$-loops).
Let $\{ \xi_l,\ldots,\xi_m \}$ be such a subset of hexagons of $S$
that form a simple $\cal T$-path $\pi_0=(\xi_l,\ldots,\xi_m)$.
Then there is a $\cal T$-path $\pi_1$ of hexagons in $\Delta D^{\delta}_{k,j}$
that goes from $\xi_l$ to $\xi_m$ without using any other hexagon of $\pi_0$
and a different $\cal T$-path $\pi_2$ in $\Delta D^{\delta}_{k,j'}$ that
goes from $\xi_m$ to $\xi_l$ without using any other hexagon of $\pi_0$.
But then, all the hexagons in $\pi_0$ other than $\xi_l$ and $\xi_m$ are
``surrounded" by $\pi_1 \cup \pi_2$ and therefore cannot have been explored
by the exploration process that produced $D^{\delta}_{k,j}$ and $D^{\delta}_{k,j'}$,
and cannot belong to $\Delta D^{\delta}_{k,j}$ or $\Delta D^{\delta}_{k,j'}$,
leading to a contradiction, unless $\pi_0 = (\xi_l,\xi_m)$.
Similar arguments lead to a contradiction if $S$ is partitioned into more
than one subset.

If $\xi_i \in S$ is not adjacent to any other hexagon in $S$, then it
is adjacent to two other hexagons of $\Delta D^{\delta}_{k,j}$ and two
hexagons of $\Delta D^{\delta}_{k,j'}$.
Since $\xi_i$ has only six neighbors and neither the two hexagons of
$\Delta D^{\delta}_{k,j}$ adjacent to $\xi_i$ nor those of $\Delta D^{\delta}_{k,j'}$
can be adjacent to each other, each hexagon of $\Delta D^{\delta}_{k,j}$
is adjacent to one of $\Delta D^{\delta}_{k,j'}$.
But then, as before, $\xi_i$ is ``surrounded" by
$\{ \Delta D^{\delta}_{k,j} \cup \Delta D^{\delta}_{k,j'} \} \setminus \xi_i$
and therefore cannot have been explored by the exploration process that produced
$D^{\delta}_{k,j}$ and $D^{\delta}_{k,j'}$, and cannot belong to $\Delta D^{\delta}_{k,j}$
or $\Delta D^{\delta}_{k,j'}$, leading once again to a contradiction.
The proof is now complete, since the only case remaining is the one where
$S$ consists of a single pair of adjacent hexagons as stated in the lemma. \fbox{} \\

With these lemmas, we can now proceed with the proof of the second
part of the theorem.
Lemma~\ref{smaller} tells us that large domains are ``chopped" with
bounded away from zero probability ($\geq p_0 > 0$), but we need to
keep track of domains of diameter larger than $\varepsilon$ in such
a way as to avoid ``double counting" as the lattice construction proceeds.
More accurately, we will keep track of domains $\tilde D^{\delta}$ having
$\text{d}_m(\tilde D^{\delta}) \geq \frac{1}{\sqrt{2}} \, \varepsilon$,
since only these can have diameter larger than $\varepsilon$.
To do so, we will associate with each domain $\tilde D^{\delta}$ having
$\text{d}_m(\tilde D^{\delta}) \geq \frac{1}{\sqrt{2}} \, \varepsilon$
that we encounter as we do the lattice construction a non-negative
integer label.
The first domain is $D_0^{\delta} = {\mathbb D}^{\delta}$ (see the beginning
of Section~\ref{full}) and this gets label $1$.
After each exploration process in a domain $\tilde D^{\delta}$ with
$\text{d}_m(\tilde D^{\delta}) \geq \frac{1}{\sqrt{2}} \, \varepsilon$,
if the number $\tilde m$ of ``daughter" subdomains $\tilde D^{\delta}_j$ with
$\text{d}_m(\tilde D_j^{\delta}) \geq \frac{1}{\sqrt{2}} \, \varepsilon$
is $0$, then the label of $\tilde D^{\delta}$ is no longer used, if instead
$\tilde m \geq 1$, then one of these $\tilde m$ subdomains (chosen by any
procedure -- e.g., the one with the highest priority for further exploration)
is assigned the \emph{same} label as $\tilde D^{\delta}$ and the rest are
assigned the next $\tilde m - 1$ integers that have never before been used
as labels.
Note that once all domains have $\text{d}_m < \frac{1}{\sqrt{2}} \, \varepsilon$,
there are no more labelled domains.

%

\begin{lemma} \label{labels}
Let $M^{\delta}_{\varepsilon}$ denote the total number of labels
used in the above procedure; then for any fixed $\varepsilon>0$,
$M^{\delta}_{\varepsilon}$ is bounded in probability as $\delta \to 0$;
i.e.,
$\lim_{M \to \infty} \limsup_{\delta \to 0}
{\mathbb P}(M^{\delta}_{\varepsilon} > M) = 0$.
\end{lemma}

\noindent {\bf Proof.} Except for $D^{\delta}_0$, every domain comes with
(at least) a ``physically correct" monochromatic ``half-boundary" (notice
that we are considering s-boundaries and that a half-boundary coming from
the ``artificially colored" boundary of $D^{\delta}_0$ is not considered
a physically correct monochromatic half-boundary).
Let us assume, without loss of generality, that $M^{\delta}_{\varepsilon}>1$.
If we associate with each label the ``last" (in terms of steps of the discrete
construction) domain which used that label (its daughter subdomains all had
$\text{d}_m < \frac{1}{\sqrt{2}} \, \varepsilon$), then we claim that it
follows from Lemma~\ref{daughters} that (with high probability) any two such
last domains that are labelled have disjoint s-boundaries.
This is a consequence of the fact that the two domains are subdomains of two
``ancestors" that are distinct daughter subdomains of the \emph{same} domain
(possibly $D^{\delta}_0$) and whose s-boundaries are therefore (by Lemma~\ref{daughters})
either disjoint or else overlap at a pair of hexagons where an exploration path
had a close encounter of distance two hexagons with itself.
But since we are dealing only with macroscopic domains (of diameter at least
order $\varepsilon$), such a close encounter would imply, like in Lemmas~\ref{sub-conv}
and~\ref{boundaries}, the existence of six crossings, not all of the same color,
of an annulus whose outer radius can be kept fixed while the inner radius is
sent to zero together with $\delta$.
The probability of such an event goes to zero as $\delta \to 0$ and hence
the unit disc $\mathbb D$ contains, with high probability, at least
$M^{\delta}_{\varepsilon}$ disjoint monochromatic $\cal T$-paths of diameter
at least $\frac{1}{\sqrt{2}} \, \varepsilon$, corresponding to the physically
correct half-boundaries of the $M^{\delta}_{\varepsilon}$ labelled domains.

Now take the collection of squares $s_j$ of side length
$\varepsilon'>0$ centered at the sites $c_j$ of a scaled
square lattice $\varepsilon' {\mathbb Z}^2$ of mesh size
$\varepsilon'$, and let $N(\varepsilon')$ be the number
of squares of side $\varepsilon'$ needed to cover the
unit disc.
Let $\varepsilon'<\varepsilon/2$ and consider the
event $\{ M^{\delta}_{\varepsilon} \geq 6 \, N(\varepsilon') \}$, which
implies that, with high probability, the unit disc contains at least
$6 \, N(\varepsilon')$ disjoint monochromatic $\cal T$-paths of diameter at
least $\frac{1}{\sqrt{2}} \, \varepsilon$ and that, for at least one $j=j_0$,
the square $s_{j_0}$ intersects at least six disjoint monochromatic
$\cal T$-paths of diameter larger that $\frac{1}{\sqrt{2}} \, \varepsilon$,
so that the ``annulus" $B(c_{j_0},\frac{1}{2 \sqrt{2}} \, \varepsilon) \setminus s_{j_0}$
is crossed by at least six disjoint monochromatic $\cal T$-paths contained
inside the unit disc.

If all these $\cal T$-paths crossing
$B(c_{j_0},\frac{1}{2 \sqrt{2}} \, \varepsilon) \setminus s_{j_0}$
have the same color, say blue, then since they are portions of boundaries
of domains discovered by exploration processes, they are ``shadowed" by
exploration paths and therefore between at least one pair of blue $\cal T$-paths,
there is at least one yellow $\cal T$-path crossing
$B(c_{j_0},\frac{1}{2 \sqrt{2}} \, \varepsilon) \setminus s_{j_0}$.
Therefore, whether the original monochromatic $\cal T$-paths are
all of the same color or not,
$B(c_{j_0},\frac{1}{2 \sqrt{2}} \, \varepsilon) \setminus s_{j_0}$
is crossed by at least six disjoint monochromatic $\cal T$-paths
not all of the same color contained in the unit disc.
%
%
Let $g(\varepsilon,\varepsilon')$ denote the $\limsup$ as $\delta \to 0$
of the probability that such an event happens anywhere inside the unit disc.
We have shown that the event $\{ M^{\delta}_{\varepsilon} \geq 6 \, N(\varepsilon') \}$
implies a ``six-arms" event unless not all labelled domains have disjoint s-boundaries.
But the latter also implies a ``six-arms" event, as discussed before; therefore
\begin{equation}
\limsup_{\delta \to 0} {\mathbb P}(M^{\delta}_{\varepsilon} \geq 6 \, N(\varepsilon'))
\leq 2 \, g(\varepsilon,\varepsilon').
\end{equation}
Since $B(c_{j_0},\frac{1}{2 \sqrt{2}} \, \varepsilon) \setminus B(c_{j_0},\frac{1}{\sqrt{2}} \, \varepsilon')
\subset B(c_{j_0},\frac{1}{2 \sqrt{2}} \, \varepsilon) \setminus s_{j_0}$,
bounds in~\cite{ksz} imply that, for $\varepsilon$ fixed,
$g(\varepsilon,\varepsilon') \to 0$ as $\varepsilon' \to 0$,
which shows that
\begin{equation}
\lim_{M \to \infty} \limsup_{\delta \to 0} {\mathbb P}(M^{\delta}_{\varepsilon} > M) = 0
\end{equation}
and concludes the proof of the lemma. \fbox{} \\

Now, let $N^{\delta}_i$ denote the number of distinct domains that had
label $i$ (this is equal to the number of steps that label $i$ survived).
Let us also define $H(\varepsilon)$ to be the smallest integer $h \geq 1$
such that $(\frac{2}{3})^h<\frac{1}{\sqrt{2}} \, \varepsilon$ and $G_h$ to
be the random variable corresponding to how many Bernoulli trials (with
probability $p_0$ of success) it takes to have $h$ successes.
Then, we may apply (sequentially) Lemma~\ref{smaller} to conclude that
for any $i$
\begin{equation}
{\mathbb P}(N^{\delta}_i \geq k+1) \leq {\mathbb P}(G_{H(\varepsilon)} + G'_{H(\varepsilon)} \geq k),
\end{equation}
where $G'_h$ is an independent copy of $G_h$.

Now let $\tilde N_1(\varepsilon),\tilde N_2(\varepsilon),\ldots$ be i.i.d.
random variables equidistributed with $G_{H(\varepsilon)} + G'_{H(\varepsilon)}$.
Let $\tilde K_{\delta}(\varepsilon)$ be the number of steps needed so that
all domains left to explore have $\text{d}_m<\frac{1}{\sqrt{2}} \, \varepsilon$.
Then, for any positive integer $M$,
\begin{equation}
{\mathbb P}(\tilde K_{\delta}(\varepsilon)>C)
\leq {\mathbb P}(M^{\delta}_{\varepsilon} \geq M+1)
+{\mathbb P}(\tilde N_1(\varepsilon)+\ldots+\tilde N_M(\varepsilon) \geq C).
\end{equation}
Notice that, for fixed $M$,
${\mathbb P}(\tilde N_1(\varepsilon)+\ldots+\tilde N_M(\varepsilon) \geq C) \to 0$
as $C \to \infty$.
Moreover, for any $\hat\varepsilon>0$, by Lemma~\ref{labels},
we can choose $M_0=M_0(\hat\varepsilon)$ large enough so that
$\limsup_{\delta \to 0} {\mathbb P}(M^{\delta}_{\varepsilon} > M_0) < \hat\varepsilon$.
So, for any $\hat\varepsilon>0$, it follows that
\begin{equation}
\limsup_{C \to \infty} \limsup_{\delta \to 0} {\mathbb P}(\tilde K_{\delta}(\varepsilon)>C)
< \hat\varepsilon,
\end{equation}
which implies that
\begin{equation}
\lim_{C \to \infty} \limsup_{\delta \to 0} {\mathbb P}(\tilde K_{\delta}(\varepsilon)>C)
= 0.
\end{equation}

To conclude this part of the proof, notice that the discrete construction
cannot ``skip'' a contour and move on to explore its interior, so that all the
contours with diameter larger than $\varepsilon$ must have been found by step
$k$ if all the domains present at that step have diameter smaller than $\varepsilon$.
Therefore, $K_{\delta}(\varepsilon) \leq \tilde K_{\delta}(\varepsilon)$, which
shows that $K_{\delta}(\varepsilon)$ is bounded in probability as $\delta \to 0$. \\

For the first part of the theorem, we need to prove, for
any fixed $k \in {\mathbb N}$, joint convergence in distribution
of the first $k$ steps of a suitably reorganized discrete
construction to the first $k$ steps of the continuum one.
Later we will explain why this reorganized construction has
the same scaling limit as the one defined in Section~\ref{full}.
For each $k$, the first $k$ steps of the reorganized discrete
construction will be coupled to the first $k$ steps of the continuum
one with suitable couplings in order to obtain the convergence in
distribution of those steps of the discrete construction to the analogous
steps of the continuum one; the proof will proceed by induction in $k$.
We will explain how to reorganize the discrete construction as we go along;
in order to explain the idea of the proof, we will consider first the
cases $k=1$, $2$ and $3$, and then extend to all $k>3$. \\


\noindent $k=1$. The first step of the continuum construction
consists of an $SLE_6$ $\gamma_1$ from $-i$ to $i$ inside
$\mathbb D$.
Correspondingly, the first step of the discrete construction
consists of an exploration path $\gamma^{\delta}_1$ inside
${\mathbb D}^{\delta}$ from the e-vertex closest to $-i$ to
the e-vertex closest to $i$.
The convergence in distribution of $\gamma^{\delta}_1$ to
$\gamma_1$ is covered by statement (S). \\

\noindent $k=2$.
The convergence in distribution of the percolation exploration
path to chordal $SLE_6$ implies that we can couple $\gamma^{\delta}_1$
and $\gamma_1$ generating them as random variables on some probability
space $(\Omega',{\cal B}',{\mathbb P}')$
such that
$\text{d}(\gamma_1(\omega'),\gamma^{\delta}_1(\omega')) \to 0$
for all $\omega'$ as $k \to \infty$ (see, for example,
Corollary~1 of~\cite{billingsley1}).

Now, let $D_1$ be the domain generated by $\gamma_1$ that
is chosen for the second step of the continuum construction,
and let $c_1 \in {\cal P}$ be the highest ranking
point of $\cal P$ contained in $D_1$.
For $\delta$ small enough, $c_1$ is also contained in
${\mathbb D}^{\delta}$; let $D_1^{\delta} = D_1^{\delta}(c_1)$
be the unique connected component of the set
${\mathbb D}^{\delta} \setminus \Gamma(\gamma_1^{\delta})$
containing $c_1$ (this is well-defined with probability
close to $1$ for small $\delta$); $D_1^{\delta}$ is the domain
where the second exploration process is to be carried out.
From the proof of Lemma~\ref{boundaries}, we know that
the boundaries $\partial D_1^{\delta}$ and $\partial D_1$
of the domains $D_1^{\delta}$ and $D_1$ produced respectively
by the path $\gamma_1^{\delta}$ and $\gamma_1$ are close
with probability close to one for $\delta$ small enough.


For the next step of the discrete construction, we choose
the two e-vertices $x_1$ and $y_1$ in $\partial D_1^{\delta}$
that are closest to the points $a_1$ and $b_1$ of $\partial D_1$
selected for the coupled continuum construction (if the choice
is not unique, we can select the e-vertices with any rule to
break the tie) and call $\gamma^{\delta}_2$ the percolation
exploration path inside $D^{\delta}_1$ from $x_1$ to $y_1$.
It follows from ~\cite{ab} that
$\{ \gamma_1^{\delta},\partial D^{\delta}_1,\gamma^{\delta}_2 \}$
converge jointly in distribution along some subsequence
to some limit
$\{ \tilde\gamma_1,\partial \tilde D_1,\tilde\gamma_2 \}$.
We already know that $\tilde\gamma_1$ is distributed like $\gamma_1$
and we can deduce from the joint convergence in distribution of
$(\gamma^{\delta}_1,\partial D_1^{\delta})$ to $(\gamma_1,\partial D_1)$
(Lemma~\ref{boundaries}), that $\partial \tilde D_1$ is distributed
like $\partial D_1$.
Therefore, if we call $\gamma_2$ the $SLE_6$ path inside $D_1$ from
$a_1$ to $b_1$, Lemma~\ref{strong-smirnov}
implies that $\tilde\gamma_2$ is distributed like $\gamma_2$ and indeed
that, as $\delta \to 0$,
$\{ \gamma^{\delta}_1,\partial D^{\delta}_1,\gamma^{\delta}_2 \}$ converge
jointly in distribution to $\{ \gamma_1,\partial D_1,\gamma_2 \}$. \\

\noindent $k = 3$. So far, we have proved the convergence
in distribution of the (paths and boundaries produced in the)
first two steps of the discrete construction to the (paths and
boundaries produced in the) first two steps of the discrete
construction.
The third step of the continuum construction consists of an $SLE_6$
path $\gamma_3$ from $a_2 \in \partial D_2$ to $b_2 \in \partial D_2$,
inside the domain $D_2$ with highest priority after the second step
has been completed.
Let $c_2 \in {\cal P}$ be the highest ranking point of $\cal P$ contained
in $D_2$, $D_2^{\delta}$ the domain of the discrete construction containing
$c_2$ after the second step of the discrete construction has been completed
(this is well defined with probability close to $1$ for small $\delta$),
and choose the two e-vertices $x_2$ and $y_2$ in $\partial D_2^{\delta}$
that are closest to the points $a_2$ and $b_2$ of $\partial D_2$ selected
for the coupled continuum construction (if the choice is not unique, we
can select the e-vertices with any rule to break the tie).
The third step of the discrete construction consists of an exploration
path $\gamma_3^{\delta}$ from $x_2$ to $y_2$ inside $D_2^{\delta}$.

It follows from ~\cite{ab} that
$\{ \gamma_1^{\delta},\partial D^{\delta}_1,\gamma^{\delta}_2,\partial D_2^{\delta},\gamma_3^{\delta} \}$
converge jointly in distribution along some subsequence
to some limit
$\{ \tilde\gamma_1,\partial\tilde D_1,\tilde\gamma_2,\partial\tilde D_2,\tilde\gamma_3 \}$.
We already know that $\tilde\gamma_1$ is distributed like $\gamma_1$,
$\partial\tilde D_1$ like $\partial D_1$ and $\tilde\gamma_2$ like $\gamma_2$,
and we would like to apply Lemma~\ref{strong-smirnov} to conclude
that $\tilde\gamma_3$ is distributed like $\gamma_3$
and indeed that, as $\delta \to 0$,
$(\gamma_1^{\delta},\partial D^{\delta}_1,\gamma^{\delta}_2,\partial D_2^{\delta},\gamma_3^{\delta})$
converges in distribution to $(\gamma_1,\partial D_1,\gamma_2,\partial D_2,\gamma_3)$.
In order to do so, we have to first show that $\partial\tilde D_2$
is distributed like $\partial D_2$.
If $D^{\delta}_2$ is a subset of ${\mathbb D}^{\delta} \setminus \Gamma(\gamma^{\delta}_1)$,
this follows from Lemma~\ref{boundaries}, as in the previous case,
but if the s-boundary of $D^{\delta}_2$ contains hexagons of
$\Gamma(\gamma^{\delta}_2)$, then we cannot use Lemma~\ref{boundaries}
directly, although the proof of the lemma can be easily adapted to
the present case, as we now explain.

\begin{figure}[!ht]
\begin{center}
\includegraphics[width=8cm]{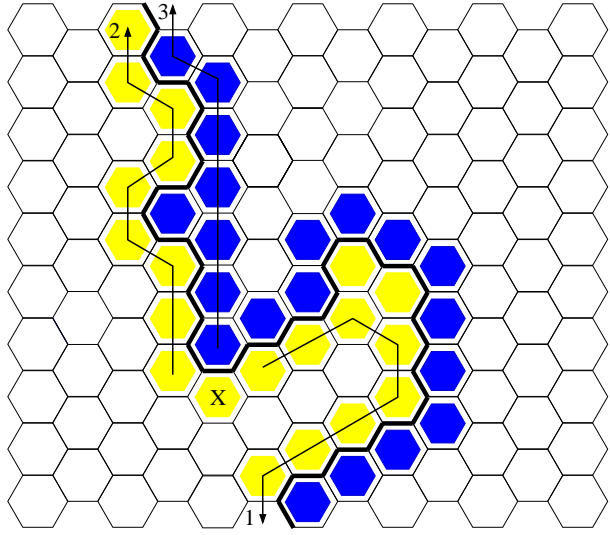}
\caption{Hexagon X, in the s-boundary of the domain $D_j^{\delta}$
to the left of the exploration path indicated by a heavy line,
has three neighbors that are the starting points of two
disjoint yellow $\cal T$-paths (denoted 1 and 2) belonging to
the s-boundary of $D_j^{\delta}$ and one blue $\cal T$-path
(denoted 3) belonging to the adjacent percolation cluster.}
\label{fig2-thm1}
\end{center}
\end{figure}

Indeed, the only difference is in the proof of claim (C) and is
due to the fact that, when dealing with a domain of type~1 or 2,
we cannot use the bound on the probability of three disjoint
crossings of a semi-annulus because the domains we are dealing
with may not be convex (like the unit disc).
On the other hand, the discrete domains like $D^{\delta}_1$ and
$D^{\delta}_2$ where we have to run exploration processes at
various steps of the discrete construction are themselves generated
by previous exploration processes, so that any hexagon of the
s-boundary of such a domain has three adjacent hexagons which are
the starting points of three disjoint $\cal T$-paths (two of
one color and one of the other).
Two of these $\cal T$-paths belong to the s-boundary
of the domain, while the third belongs to the adjacent
percolation cluster (see Figure~\ref{fig2-thm1}).
This allows us to use the bound on the probability of
\emph{six} disjoint crossings of an annulus.

To see this, let $\pi_1, \pi_2$ be the $\cal T$-paths
contained in the s-boundary of the discrete domain (i.e.,
$D^{\delta}_1$ in the present context) and $\pi_3$ the
$\cal T$-path belonging to the adjacent cluster, all starting
from hexagons adjacent to some hexagon $\xi$ (centered at $u$)
in the s-boundary of $D^{\delta}_1$.
For $0< \varepsilon' < \varepsilon$ and $\delta$ small
enough, let ${\cal A}_u(\varepsilon,\varepsilon')$ be the
event that the exploration path $\gamma_2^{\delta}$
enters the ball $B(u,\varepsilon')$ without touching
$\partial D^{\delta}_1$ inside the larger ball $B(u,\varepsilon)$.
${\cal A}_u(\varepsilon,\varepsilon')$ implies having (at least)
three disjoint $\cal T$-paths (two of one color and one of the other),
$\pi_4, \pi_5$ and $\pi_6$, contained in $D^{\delta}_1$ and crossing
the annulus $B(u,\varepsilon) \setminus B(u,\varepsilon')$, with
$\pi_4, \pi_5$ and $\pi_6$ disjoint from $\pi_1, \pi_2$ and $\pi_3$.
Hence, ${\cal A}_u(\varepsilon,\varepsilon')$ implies the event that
there are (at least) six disjoint crossings (not all of the same color)
of the annulus $B(u,\varepsilon) \setminus B(u,\varepsilon')$.

Once claim (C) is proved, the rest of the proof of
Lemma~\ref{boundaries} applies to the present case.
Therefore, we have convergence in distribution of
$\partial D_2^{\delta}$ to $\partial D_2$, which allows
us to use Lemma~\ref{strong-smirnov} and conclude that
$(\gamma_1^{\delta},\gamma_2^{\delta},\gamma_3^{\delta})$
converges in distribution to $(\gamma_1,\gamma_2,\gamma_3)$. \\

\noindent $k > 3$.
We proceed by induction in $k$, iterating the steps explained
above; there are no new difficulties; all steps for $k \geq 4$
are analogous to the case $k = 3$. \\

To conclude the proof of the theorem, we need to show
that the scaling limit of the original full discrete
construction defined in Section~\ref{full} is the same
as that of the reorganized one just used in the proof
of the first part of the theorem.
In order to do so, we can couple the two constructions
by using the same percolation configuration for both,
so that the two constructions have at their disposal
the same set of loops to discover.
We proved above that the original discrete construction
finds \emph{all} the ``macroscopic" loops, so we have to
show that this is true also for the reorganized version
of the discrete construction.
This is what we will do next, using essentially the same
arguments as those employed for the original discrete
construction; we present these arguments for the sake
of completeness since there are some changes.

Consider the reorganized discrete construction described
above, where the starting and ending points of the exploration
processes at each step are chosen to be close to those of
the corresponding (coupled) continuum construction.
Suppose that at step $k$ of this discrete construction
an exploration process $\gamma^{\delta}_k$ is run inside
a domain $D^{\delta}_{k-1}$, and write
$D^{\delta}_{k-1} \setminus \Gamma(\gamma^{\delta}_k) =
\bigcup_j D^{\delta}_{k,j}$, where $\{ D^{\delta}_{k,j} \}$
are the connected domains into which $D^{\delta}_{k-1}$
is split by the set $\Gamma(\gamma^{\delta}_k)$ of hexagons
explored by $\gamma^{\delta}_k$.

Let $\text{d}_x(D_{k-1})$ (resp., $\text{d}_x(D^{\delta}_{k-1})$)
and $\text{d}_y(D_{k-1})$ (resp., $\text{d}_y(D^{\delta}_{k-1})$)
be respectively the maximal $x$- and $y$-distance between pairs
of points in $\partial D_{k-1}$ (resp., $\partial D^{\delta}_{k-1}$).
%
%
%
If
$\text{d}_x(D^{\delta}_{k-1}) \geq \text{d}_y(D^{\delta}_{k-1})$
and the e-vertices on $\partial D^{\delta}_{k-1}$
are chosen to be closest to two points of $\partial D_{k-1}$
with maximal $x$-distance, then the same construction and
argument spelled out earlier in the first part of the proof
(corresponding to the second part of the theorem) show that
$\max_j \text{d}_x(D^{\delta}_{k,j}) \leq \frac{2}{3} \text{d}_x(D^{\delta}_{k-1})$
with bounded away from zero probability.

%

If the e-vertices on $\partial D^{\delta}_{k-1}$ are
chosen to be closest to two points of $\partial D_{k-1}$
with maximal $x$-distance but
$\text{d}_x(D^{\delta}_{k-1}) \leq \text{d}_y(D^{\delta}_{k-1})$,
then consider the rectangle ${\cal R}'$ whose vertical sides
are aligned to the $y$-axis, have length $\text{d}_y(D^{\delta}_{k-1})$,
and are each placed at the same $x$-distance from the
points of $\partial D^{\delta}_{k-1}$ with minimal
or maximal $x$-coordinate in such a way that the
horizontal sides of ${\cal R}'$ have length
$\frac{1}{3} \, \text{d}_y(D^{\delta}_{k-1})$;
the bottom and top sides of ${\cal R}'$ are placed in
such a way that they touch the points of
$\partial D^{\delta}_{k-1}$ with minimal or maximal
$y$-coordinate, respectively.
Notice that, because of the coupling between the continuum
and discrete constructions, for any $\tilde\varepsilon >0$, for
$k$ large enough, $|\text{d}_x(D_{k-1}^{\delta}) - \text{d}_x(D_{k-1})| \leq \tilde\varepsilon$
and $|\text{d}_y(D_{k-1}^{\delta}) - \text{d}_y(D_{k-1})| \leq \tilde\varepsilon$.
Since in the case under consideration we have
$\text{d}_y(D_{k-1}^{\delta}) \geq \text{d}_x(D_{k-1}^{\delta})$ and
$\text{d}_x(D_{k-1}) \geq \text{d}_y(D_{k-1})$, for $\delta$ large enough,
we must also have
$|\text{d}_y(D_{k-1}^{\delta}) - \text{d}_x(D_{k-1}^{\delta})| \leq 2 \, \tilde\varepsilon$.
Once again, it follows from the Russo-Seymour-Welsh lemma
that the probability to have two vertical $\cal T$-crossings
of ${\cal R}'$ of different colors is bounded away from zero by
a positive constant that does not depend on $\delta$ (for $\delta$
small enough).
If that happens, then
$\max_j \text{d}_x(D^{\delta}_{k,j}) \leq
\frac{2}{3} \, \text{d}_x(D^{\delta}_{k-1}) + \frac{1}{3} \, \tilde\varepsilon$.

All other cases are handled in the same way, implying that
the maximal $x$- and $y$-distances of domains that appear
in the discrete construction have a positive probability
(bounded away from zero) to decrease by (approximately) a
factor $2/3$ at each step of the discrete construction in
which an exploration process is run in that domain.


With this result at our disposal, the rest of the proof, that
for any $\varepsilon>0$ the number of steps needed to find
all the loops of diameter larger than $\varepsilon$ is bounded
in probability as $\delta \to 0$ (which implies that \emph{all}
the ``macroscopic" loops are discovered), proceeds exactly like
for the original discrete construction. \fbox{} \\

\noindent {\bf Proof of Theorem~\ref{thm-therm-lim}.}
First of all, we want to show that $P^D \equiv \hat I_D P_R$
does not depend on $R$, provided $D$ is strictly contained
in ${\mathbb D}_R$ and
$\partial D \cap \partial {\mathbb D}_R = \emptyset$.
In order to do this, we assume that the above
conditions are satisfied for the pair $D,R$ and
show that $\hat I_D P_R = \hat I_D P_{R'}$ for all $R' > R$.

Take two copies of the scaled hexagonal lattice,
$\delta{\cal H}$ and $\delta{\cal H}'$, their
dual lattices $\delta{\cal T}$ and $\delta{\cal T}'$,
and two percolation configurations,
$\sigma_{{\mathbb D}_R}$ and ${\sigma'}_{{\mathbb D}_{R'}}$,
both with blue boundary conditions and coupled in
such a way that
$\sigma_{{\mathbb D}_R} = \sigma'_{{\mathbb D}_R}$.
The laws of the boundaries of $\sigma$ and $\sigma'$
are also coupled, in such a way that the boundaries
or portions of boundaries contained inside $D$ are
identical for all small enough $\delta$.
Therefore, letting $\delta \to 0$ and using the
convergence of the percolation boundaries inside
${\mathbb D}_R$ and ${\mathbb D}_{R'}$ to the
continuum nonsimple loop processes $P_R$ and $P_{R'}$
respectively, we conclude that $\hat I_D P_R = \hat I_D P_{R'}$.

From what we have just proved, it follows that the probability
measures $P^{{\mathbb D}_R}$ on $(\Omega_R,{\cal B}_R)$, for
$R \in {\mathbb R}_+$, satisfy the consistency conditions
$P^{{\mathbb D}_{R_1}} = \hat I_{{\mathbb D}_{R_1}} P^{{\mathbb D}_{R_2}}$
for all $R_1 \leq R_2$.
Since $\Omega_R$, $\Omega$ are complete separable metric
spaces, the measurable spaces $(\Omega_R,{\cal B}_R)$,
$(\Omega,{\cal B})$ are standard Borel spaces and so we
can apply Kolmogorov's extension theorem (see, for example,
\cite{durrett}) and conclude that there exists a unique
probability measure on $(\Omega,{\cal B})$ with
$P^{{\mathbb D}_R} = \hat I_{{\mathbb D}_R} P$ for all
$R \in {\mathbb R}_+$.
It follows that, for $R'>R$ and all $D$ strictly contained
in ${\mathbb D}_R$ and such that
$\partial D \cap \partial {\mathbb D}_R = \emptyset$,
$\hat I_D P_R = P^D = \hat I_D P_{R'} =
\hat I_D \hat I_{{\mathbb D}_R} P_{R'} =
\hat I_D P^{{\mathbb D}_R} =
\hat I_D \hat I_{{\mathbb D}_R} P = \hat I_D P$, which
concludes the proof. \fbox{} \\

\noindent {\bf Proof of Corollary~\ref{scal-lim}.}
The corollary is an immediate consequence of
Theorems~\ref{thm-convergence} and~\ref{thm-therm-lim},
where the full scaling limit is intended in the topology
induced by~(\ref{hausdorff-D}). \fbox{} \\

\noindent {\bf Proof of Theorem~\ref{features}.}
{\it 1.} The fact that the Continuum Nonsimple Loop process
is a random collection of noncrossing continuous loops is a
direct consequence of its definition.
The fact that the loops touch themselves is a
consequence of their being constructed out of $SLE_6$,
while the fact that they touch each other follows
from the observation that a chordal $SLE_6$ path
$\gamma_{D,a,b}$ touches $\partial D$ with probability
one.
Therefore, each new loop in the continuum construction
touches one or more previous ones (many times).

The nonexistence of triple points follows directly from
Lemma~5 of~\cite{ksz} on the number of crossings of an
annulus, combined with Corollary~\ref{scal-lim}, which
allows to transport discrete results to the continuum case.
In fact, a triple point would imply, for discrete percolation,
at least six crossings (not all of the same color) of an
annulus whose ratio of inner to outer radius goes to zero
in the scaling limit, leading to a contradiction.

{\it 2.} This follows from straightforward
Russo-Seymour-Welsh type arguments for percolation
(for more details, see, for example, Lemma~3
of~\cite{ksz}), combined with Corollary~\ref{scal-lim}.

{\it 3.} 
Combining Russo-Seymour-Welsh type arguments for
percolation (see, for example, Lemma~3 of~\cite{ksz})
with Corollary~\ref{scal-lim}, we know that $P$-a.s.
there exists a (random) $R^* = R^*(R)$, with $R^* < \infty$,
such that ${\mathbb D}_R$ is surrounded by a continuum
nonsimple loop contained in ${\mathbb D}_{R^*}$.
From (the proof of) Theorem~\ref{thm-therm-lim}, we also
know that
$\hat I_{{\mathbb  D}_{R''}} P = P^{{\mathbb D}_{R''}} =
I_{{\mathbb D}_{R''}} P_{R'}$ for all $R' > R''$.
This implies that by taking $R'$ large enough
and performing the continuum construction inside ${\mathbb D}_{R'}$,
we have a positive probability of generating a loop $\lambda$
contained in the annulus ${\mathbb D}_{R''} \setminus {\mathbb D}_R$,
with $R'>R''>R$.
If that is the case, all the loops contained inside
${\mathbb D}_R$ are connected, by construction, to the
loop $\lambda$ surrounding ${\mathbb D}_R$ by a finite
sequence (a ``path'') of loops (remember that in the
continuum construction each loop is generated by pasting
together portions of $SLE_6$ paths inside domains whose
boundaries are determined by previously formed loops or
excursions).
Therefore, any two loops contained inside ${\mathbb D}_R$
are connected to each other by a ``path" of loops.

Using again the fact that
$\hat I_{{\mathbb  D}_{R''}} P = P^{{\mathbb D}_{R''}} =
I_{{\mathbb D}_{R''}} P_{R'}$ for all $R' > R''$, and letting first
$R'$ and then $R''$ go to $\infty$, we see from the discussion above
(with $R \to \infty$ as well) that any two loops are connected by a
finite ``path" of intermediate loops, $P$-a.s.

{\it 4.} In order to prove the claim, we will define a
discrete construction inside $D'$ coupled to the continuum
construction inside $D$, by means of the conformal map $f$
from $D$ to $D'$.
Roughly speaking, this new discrete construction for $D'$
is one in which the $(x,y)$ pairs at each step are chosen
to be closest to the $(f(a),f(b))$ points in $D'$ mapped
from $D$ via $f$, where the pairs $(a,b)$ are those that
appear at the corresponding steps of the continuum
construction inside $D$.

More precisely, let $\gamma_1$ be the first $SLE_6$ path
in $D$ from $a_1$ to $b_1$.
Because of the conformal invariance of $SLE_6$, the image
$f(\gamma_1)$ of $\gamma_1$ under $f$ is a path distributed
as the trace of chordal $SLE_6$ in $D'$ from $f(a_1)$ to $f(b_1)$.
Therefore, the exploration path $\gamma^{\delta}_1$ inside
$D'$ from $x_1$ to $y_1$, chosen to be closest to $f(a_1)$
and $f(b_1)$ respectively, converges in distribution to $f(\gamma_1)$,
as $\delta \to 0$, which means that there exist a coupling
between $\gamma^{\delta}_1$ and $f(\gamma_1)$ such that the
paths stay close for $\delta$ small.

We see already that one can use the same strategy as in the
proof of the first part of Theorem~\ref{thm-convergence},
and obtain a discrete construction whose exploration paths
are coupled to the $SLE_6$ paths $f(\gamma_k)$ that are the
images of the paths $\gamma_k$ in $D$.
Then, for this discrete construction, the scaling limit of
the exploration paths will be distributed as the images of
the $SLE_6$ paths in $D$.

In order to conclude the proof, we just have to show that
the discrete construction inside $D'$ defined above finds
all the boundaries in a number of steps that is bounded in
probability as $\delta \to 0$ (this is equivalent to the
second part of Theorem~\ref{thm-convergence}).
To do that, we use the second part of Theorem~\ref{thm-convergence},
which implies that, for any fixed $\varepsilon > 0$ and
$C < \infty$, the number of steps of the discrete construction
in $D$ that are necessary to ensure that only domains with
diameter less than $\varepsilon / C$ are present is bounded
in probability as $\delta \to 0$.
Since $f$ (can be extended to a function that) is continuous
in the compact set $\overline D$, $f$ is uniformly continuous
and so we can now choose $C=C(f) < \infty$ such that any subdomain
of $D$ of diameter at most $\varepsilon / C$ is mapped by $f$
to a subdomain of $D'$ of diameter at most $\varepsilon$.
This, combined with the coupling between $SLE_6$ paths and
exploration paths inside $D'$, assures that the number of
steps necessary for the new discrete construction inside
$D'$ to find all the loops of diameter at least $\varepsilon$
is bounded in probability as $\delta \to 0$.

Therefore, the scaling limit, as $\delta \to 0$, of this
new discrete construction for $D'$ gives the measure
$P_{D'}$.
It follows by construction that $f * P_D = P_{D'}$,
which concludes the proof.~\fbox{}

%

\appendix

\refstepcounter{section}
\section*{Appendix \thesection: Convergence of the Percolation Exploration Path} \label{convergence-to-sle}

In this first appendix, we provide a detailed proof of statement (S).
%
%
%
The existence of subsequential limits for the percolation exploration
path, which follows from the work of Aizenman and Burchard~\cite{ab},
means that the proof can be divided into two parts: first we will give
a characterization of chordal $SLE_6$ in terms of two properties that
determine it uniquely; then we will show that any subsequential scaling
limit of the percolation exploration path satisfies these two properties.

The characterization part will follow from known properties of hulls
and of $SLE_6$ (see~\cite{lsw7} and~\cite{werner}).
The second part will follow from an extension of Smirnov's result about the
convergence of crossing probabilities to Cardy's formula~\cite{cardy,cardy2}
(see Theorem~\ref{strong-cardy} below) for sequences of Jordan domains $D_k$,
with the domain $D_k$ changing together with the mesh $\delta_k$ of the lattice,
combined with the proof of a certain spatial Markov property for subsequential
limits of percolation exploration hulls (Theorem~\ref{spatial-markov}).
We note that although Theorem~\ref{strong-cardy} represents only a slight
extension to Smirnov's result on convergence of crossing probabilities, this
extension and its proof play a major role in the technically important
Lemmas~\ref{double-crossing}, \ref{equal} and~\ref{mushroom}, which control
the ``close encounters" of exploration paths to domain boundaries.
The proof of Theorem~\ref{strong-cardy} is modelled after a simpler
geometric argument involving only rectangles used in~\cite{cns1}.

Let $D'$ be a bounded simply connected domain containing the origin
whose boundary $\partial D'$ is a continuous curve.
Let $f:{\mathbb D} \to D'$ be the (unique) conformal map from the
unit disc to $D'$ with $f(0)=0$ and $f'(0)>0$; note that by
Theorem~\ref{cont-thm} of Appendix~\ref{rado}, $f$ has a continuous
extension to $\overline{\mathbb D}$.
Let $z_1,z_2,z_3,z_4$ be four points of $\partial D'$ (or more
accurately, four prime ends) in counterclockwise order -- i.e.,
such that $z_j=f(w_j), \,\,\, j=1,2,3,4$, with $w_1,\ldots,w_4$
in counterclockwise order.
Also, let $\eta = \frac{(w_1-w_2)(w_3-w_4)}{(w_1-w_3)(w_2-w_4)}$.
Cardy's formula~\cite{cardy,cardy2} for the probability $\Phi_{D'}(z_1,z_2;z_3,z_4)$
of a ``crossing" inside $D'$ from the counterclockwise arc $\overline{z_1 z_2}$ to
the counterclockwise arc $\overline{z_3 z_4}$ is
\begin{equation} \label{cardy-formula}
\Phi_{D'}(z_1,z_2;z_3,z_4) = \frac{\Gamma(2/3)}{\Gamma(4/3) \Gamma(1/3)} \eta^{1/3} {}_2F_1(1/3,2/3;4/3;\eta),
\end{equation}
where ${}_2F_1$ is a hypergeometric function.

For a given mesh $\delta>0$, the probability of a crossing inside $D'$
from the counterclockwise arc $\overline{z_1 z_2}$ to the counterclockwise
arc $\overline{z_3 z_4}$ is the probability of the existence of a blue
$\cal T$-path contained in $(D')^{\delta}$, the $\delta$-approximation of
$D'$ (see Definition~\ref{approx} above), that starts at a hexagon adjacent
to one intersecting $\overline{z_1 z_2}$ and ends at a hexagon adjacent
to one intersecting $\overline{z_3 z_4}$.
Smirnov proved the following major theorem, concerning the conjectured
behavior~\cite{cardy,cardy2} of crossing probabilities in the scaling limit.
\begin{theorem} \emph{(Smirnov~\cite{smirnov})} \label{cardy-smirnov}
In the limit $\delta \to 0$, the crossing probability becomes conformally
invariant and converges to Cardy's formula~(\ref{cardy-formula}).
\end{theorem}

\begin{remark}
We actually only need Theorem~\ref{cardy-smirnov} for Jordan domains,
as can be seen by a careful reading of the proof of Theorem~\ref{strong-cardy}.
We note that Smirnov does not restrict attention to that case.
\end{remark}

Let us now specify the objects that we are interested in.
Suppose $D$ is a simply connected domain whose boundary $\partial D$
is a continuous curve, and $a,b$ are two distinct points in $\partial D$
(or more accurately, two distinct prime ends), and let $\tilde\mu_{D,a,b}$
be a probability measure on continuous (non-self-crossing) curves
$\tilde\gamma = \tilde\gamma_{D,a,b}:[0,\infty) \to \overline D$
with $\tilde\gamma(0)=a$ and
$\tilde\gamma(\infty) \equiv \lim_{t \to \infty} \tilde\gamma(t) = b$
(we remark that the use of $[0,\infty)$ instead of $[0,1]$ for the time
parametrization is purely for convenience).
Let $D_t \equiv D \setminus \tilde K_t$ denote the (unique) connected
component of $D \setminus \tilde\gamma[0,t]$ whose closure contains $b$,
where $\tilde K_t$, the {\bf filling} of $\tilde\gamma[0,t]$, is a closed
connected subset of $\overline D$.
$\tilde K_t$ is called a {\bf hull} if it satisfies the condition
\begin{equation} \label{hulls}
\overline{\tilde K_t \cap D} = \tilde K_t.
\end{equation}
We will generally be interested in curves $\tilde\gamma$ such that $\tilde K_t$
is a hull for each $t$, although we normally only consider $\tilde K_T$ at
certain stopping times $T$.
(An example of such a curve that we are particularly interested in is the
trace of chordal $SLE_6$.)

\begin{figure}[!ht]
\begin{center}
\includegraphics[width=8cm]{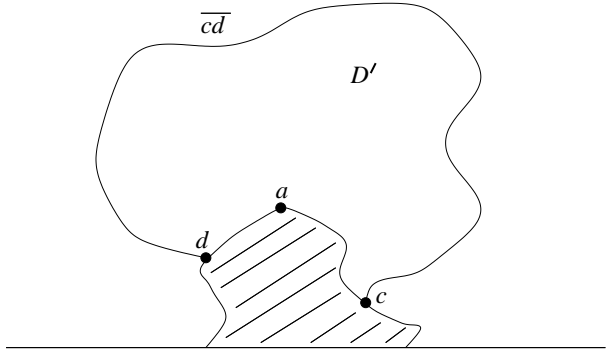}
\caption{$D$ is the upper half-plane $\mathbb H$ with the shaded portion removed,
$C'$ is an unbounded subdomain, and $D' = D \setminus C'$ is indicated in the figure.}
\label{fig1-sec7}
\end{center}
\end{figure}

\begin{figure}[!ht]
\begin{center}
\includegraphics[width=8cm]{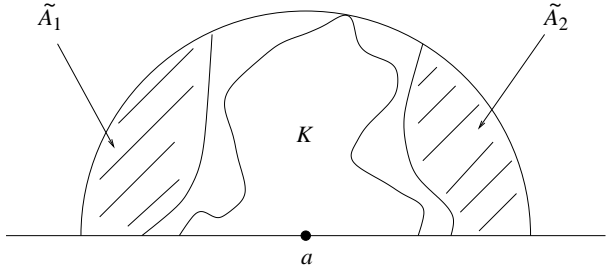}
\caption{Example of a hull $K$ and a set $\tilde A_1 \cup \tilde A_2$ in $\cal A$.
Here, $D = {\mathbb H}$ and $D'$ is the semi-disc centered at $a$.}
\label{fig2-sec7}
\end{center}
\end{figure}

Let $C' \subset D$ be a closed subset of $\overline D$ such that $a \notin C'$,
$b \in C'$, and $D' = D \setminus C'$ is a bounded simply connected domain
whose boundary contains the counterclockwise arc $\overline {cd}$ that does
not belong to $\partial D$ (except for its endpoints $c$ and $d$ -- see
Figure~\ref{fig1-sec7}).
Let $T'=\inf \{ t:\tilde K_t \cap C' \neq \emptyset \}$ be the first time
that $\tilde\gamma(t)$ hits $C'$ and assume that the filling $\tilde K_{T'}$
of $\tilde\gamma[0,T']$ is a hull; we denote by $\tilde\nu_{D',a,c,d}$
the distribution of $\tilde K_{T'}$.
To explain what we mean by the distribution of a hull, consider the set
$\tilde{\cal A}$ of closed subsets $\tilde A$ of $\overline{D'}$ that do not
contain $a$ and such that $\partial \tilde A \setminus \partial D'$ is a simple
(continuous) curve contained in $D'$ starting at a point on $\partial D' \cap D$
and ending at a point on $\partial D$ (see Figure~\ref{fig2-sec7}).
Let $\cal A$ be the set of closed subsets of $\overline{D'}$ of the form
$\tilde A_1 \cup \tilde A_2$, where $\tilde A_1, \tilde A_2 \in \tilde{\cal A}$
and $\tilde A_1 \cap \tilde A_2 = \emptyset$.

For a given $C'$ and corresponding $T'$, let $\cal K$ be the set whose
elements are possible hulls at time $T'$; we claim that the events
$\{ K \in {\cal K}: K \cap A = \emptyset \}$, for $A \in {\cal A}$,
form a $\pi$-system $\Pi$ (i.e., they are closed under finite intersections;
we also include the empty set in $\Pi$), and we consider the $\sigma$-algebra
$\Sigma=\sigma(\Pi)$ generated by these events.
To see that $\Pi$ is closed under pairwise intersections, notice that, if
$A_1, A_2 \in {\cal A}$, then $\{ K \in {\cal K}: K \cap A_1 = \emptyset \}
\cap \{ K \in {\cal K}: K \cap A_2 = \emptyset \} =
\{ K \in {\cal K}: K \cap \{ A_1 \cup A_2 \} = \emptyset \}$
and $A_1 \cup A_2 \in {\cal A}$ (or else
$\{ K \in {\cal K} : K \cap \{ A_1 \cup A_2 \} = \emptyset \}$ is empty).
We are interested in probability spaces of the form $({\cal K},\Sigma,{\mathbb P}^*)$.

We say that the exit distribution of $\tilde\gamma(t)$ is determined by Cardy's formula
if, for any $C'$ and any counterclockwise arc $\overline{xy}$ of $\overline{cd}$, the
probability that $\tilde\gamma$ hits $C'$ at time $T'$ on $\overline{xy}$ is given by
\begin{equation}
{\mathbb P}^*(\tilde\gamma(T') \in \overline{xy}) = \Phi_{D'}(a,c;x,d) - \Phi_{D'}(a,c;y,d).
\end{equation}


It is easy to see that if the hitting distribution of $\tilde\gamma(t)$ is
determined by Cardy's formula, then the probabilities of events in $\Pi$
are also determined by Cardy's formula in the following way.
Let $A \in {\cal A}$ be the union of $\tilde A_1, \tilde A_2 \in \tilde{\cal A}$,
with $\partial \tilde A_1 \setminus \partial D'$ given by a curve from
$u_1 \in \partial D' \cap D$ to $v_1 \in \partial D$ and
$\partial \tilde A_2 \setminus \partial D'$ given by a curve from
$u_2 \in \partial D' \cap D$ to $v_2 \in \partial D$; then, assuming that
$a$, $v_1$, $u_1$, $u_2$, $v_2$ are ordered counterclockwise around $\partial D'$,
\begin{equation}
{\mathbb P}^*(\tilde K_{T'} \cap A = \emptyset) =
\Phi_{D' \setminus A}(a,v_1;u_1,v_2,) - \Phi_{D' \setminus A}(a,v_1;u_2,v_2).
\end{equation}
Since $\Pi$ is a $\pi$-system, the probabilities of the events in $\Pi$ determine
uniquely the distribution of the hull in the sense described above.
Therefore, if we let $\gamma_{D,a,b}$ denote the trace of chordal $SLE_6$
inside $D$ from $a$ to $b$, $K_t$ its hull up to time $t$, and
$\tau = \inf \{ t : K_t \cap C' \neq \emptyset \}$ the first time that
$\gamma_{D,a,b}$ hits $C'$, we have the following simple but useful lemma.

\begin{lemma} \label{hull}
With the notation introduced above, if the hitting distribution 
of $\tilde\gamma_{D,a,b}$ is determined by Cardy's formula and $\tilde K_{T'}$ is a hull,
then $\tilde K_{T'}$ is distributed like the hull $K_{\tau}$ of chordal $SLE_6$.
\end{lemma}

\noindent {\bf Proof.} It is enough to note that the hitting distribution for
chordal $SLE_6$ is determined by Cardy's formula~\cite{lsw1}. \fbox{} \\

Now let $\tilde f_0$ be a conformal map from $D$ to the upper half-plane
$\mathbb H$ such that $\tilde f_0(a)=0$ and $\lim_{z \to b} \tilde f_0(z)=\infty$
(these two conditions determine $\tilde f_0$ only up to a multiplicative constant).
For $\varepsilon>0$ fixed, let
$C(u,\varepsilon) = \{ z:|u-z|<\varepsilon \} \cap {\mathbb H}$
denote the semi-ball of radius $\varepsilon$ centered at $u$ on
the real line and let $\tilde T_1=\tilde T_1(\varepsilon)$ denote
the first time $\tilde\gamma(t)$ hits $D \setminus \tilde G_1$,
where $\tilde G_1 \equiv \tilde f_0^{-1}(C(0,\varepsilon))$.
Define recursively $\tilde T_{j+1}$ as the first time
$\tilde\gamma[\tilde T_j,\infty)$ hits $\tilde D_{\tilde T_j} \setminus \tilde G_{j+1}$,
where $\tilde D_{\tilde T_j} \equiv D \setminus \tilde K_{\tilde T_j}$,
$\tilde G_{j+1} \equiv \tilde f^{-1}_{\tilde T_j}(C(0,\varepsilon))$, and
$\tilde f_{\tilde T_j}$ is a conformal map from $\tilde D_{\tilde T_j}$ to
$\mathbb H$ that maps $\tilde\gamma(\tilde T_j)$ to $0$ and $b$ to $\infty$.
We also define $\tilde\tau_{j+1} \equiv \tilde T_{j+1} - \tilde T_j$,
so that $\tilde T_j = \tilde\tau_1+\ldots+\tilde\tau_j$.
We note that, like $\tilde f_0$, the conformal maps $\tilde f_{\tilde T_j}$
are only defined up to a multiplicative factor.

Notice that $\tilde G_{j+1}$ is a bounded simply connected domain chosen
so that the conformal transformation which maps $\tilde D_{\tilde T_j}$
to $\mathbb H$ maps $\tilde G_{j+1}$ to the semi-ball $C(0,\varepsilon)$
centered at the origin on the real line.
With these definitions, consider the (discrete-time) stochastic process
$\tilde X_j \equiv (\tilde K_{\tilde T_j}, \tilde\gamma(\tilde T_j))$ for $j=1,2,\ldots$;
we say that $\tilde K_t$ satisfies the {\bf spatial Markov property} if each
$\tilde K_{\tilde T_j}$ is a hull and $\tilde X_j$ for $j=1,2,\ldots$ is a
Markov chain (for any choice of the multitplicative factors for
$\tilde f_0, \tilde f_{\tilde T_1}, \tilde f_{\tilde T_2}, \ldots$).
Notice that the trace of chordal $SLE_6$ satisfies the spatial Markov property,
due to the conformal invariance and Markovian properties~\cite{schramm} of $SLE_6$.

Next, we give the main characterization theorem.
\begin{theorem} \label{characterization}
If the filling process $\tilde K_t$ of $\tilde\gamma_{D,a,b}$ satisfies the spatial
Markov property and its hitting distribution is determined by Cardy's formula, then
$\tilde\gamma_{D,a,b}$ is distributed like the trace $\gamma_{D,a,b}$ of chordal
$SLE_6$ inside $D$ started at $a$ and aimed at $b$.
\end{theorem}

\noindent {\bf Proof.} Since the trace $\gamma_{D,a,b}$ of chordal $SLE_6$ in a
bounded Jordan domain $D$ is defined (up to a linear time change) as $f(\gamma)$,
where $\gamma = \gamma_{{\mathbb H},0,\infty}$ is the trace of chordal $SLE_6$
in the upper half-plane started at $0$ and $f$ is any conformal map from the upper
half-plane $\mathbb H$ to $D$ such that $f(0)=a$ and $f(\infty)=b$, it is enough to
show that $\hat\gamma = f^{-1}(\tilde\gamma_{D,a,b})$ is distributed like the trace
of chordal $SLE_6$ in the upper half-plane.
Let $\hat K_t$ denote the filling of $\hat\gamma(t)$ at time $t$ and let $\hat g_t(z)$
be the unique conformal transformation that maps ${\mathbb H} \setminus \hat K_t$
onto $\mathbb H$ with the following expansion at infinity:
\begin{equation} \label{conf-map1}
\hat g_t(z) = z + \frac{\hat a(t)}{z} + o(\frac{1}{z}).
\end{equation}
We choose to parametrize $\hat\gamma(t)$ so that $t = \frac{\hat a(t)}{2}$
(this is often called parametrization by capacity, $\hat a(t)$ being the
capacity of the filling up to time $t$).

We want to compare $\hat\gamma(t)$ with the trace $\gamma(t)$ of chordal $SLE_6$
in the upper half-plane parametrized in the same way (i.e., with $a(t) = 2 \, t$),
so that, if $K_t$ denotes the filling of $\gamma$ at time $t$,
${\mathbb H} \setminus K_t$ is mapped onto $\mathbb H$ by a conformal
transformation with the following expansion at infinity:
\begin{equation} \label{conf-map2}
g_t(z) = z + \frac{2t}{z} + o(\frac{1}{z}).
\end{equation}

Our strategy will be to construct suitable polygonal approximations
$\hat\gamma_{\varepsilon}$ and $\gamma_{\varepsilon}$ of $\hat\gamma$ and
$\gamma$ which converge, as $\varepsilon \to 0$, to the original curves
(in the uniform metric on continuous curves~(\ref{distance})), and show that
$\hat\gamma_{\varepsilon}$ and $\gamma_{\varepsilon}$ have the same distribution.
This implies that the distributions of $\hat\gamma$ and $\gamma$ must coincide.

Let us describe first the construction for $\gamma_{\varepsilon}(t)$; we use
exactly the same construction for $\hat\gamma_{\varepsilon}(t)$.
We remark that the important features in the construction of the polygonal
approximations are the spatial Markov property of the fillings and Cardy's
formula, which are valid for both $\gamma$ and $\hat\gamma$.

For $\varepsilon>0$ fixed, as above let
$C(u,\varepsilon) = \{ z:|u-z|<\varepsilon \} \cap {\mathbb H}$
denote the semi-ball of radius $\varepsilon$ centered at $u$ on the real line.
Let $T_1=T_1(\varepsilon)$ denote the first time $\gamma(t)$ hits
${\mathbb H} \setminus G_1$, where $G_1 \equiv C(0,\varepsilon)$, and
define recursively $T_{j+1}$ as the first time $\gamma[T_j,\infty)$
hits ${\mathbb H}_{T_j} \setminus G_{j+1}$, where
${\mathbb H}_{T_j} = {\mathbb H} \setminus K_{T_j}$ and
$G_{j+1} \equiv g^{-1}_{T_j}(C(g_{T_j}(\gamma(T_j)),\varepsilon))$.
Notice that $G_{j+1}$ is a bounded simply connected domain chosen so that the
conformal transformation which maps ${\mathbb H}_{T_j}$ to $\mathbb H$ maps
$G_{j+1}$ to the semi-ball $C(g_{T_j}(\gamma(T_j)),\varepsilon)$ centered at
the point of the real line where the ``tip" $\gamma(T_j)$ of the hull $K_{T_j}$
is mapped.
The spatial Markov property and the conformal invariance of the hull of $SLE_6$
imply that if we write $T_j = \tau_1+\ldots+\tau_j$, with
$\tau_{j+1} \equiv T_{j+1} - T_j$, the $\tau_j$'s are i.i.d. random variables,
and also that the distribution of $K_{T_{j+1}}$ is the same as the distribution
of $K_{T_j} \cup g^{-1}_{T_j}(K'_{T_1} + g_{T_j}(\gamma(T_j)))$, where $K'_{T_1}$
is equidistributed with $K_{T_1}$, but also is independent of $K_{T_1}$.
The polygonal approximation $\gamma_{\varepsilon}$ is obtained by joining, for
all $j$, $\gamma(T_j)$ to $\gamma(T_{j+1})$ with a straight segment, where the
speed $\gamma_{\varepsilon}'(t)$ is chosen to be constant.

Now let $\hat T_1 = \hat T_1(\varepsilon)$ denote the first time
$\hat\gamma(t)$ hits ${\mathbb H} \setminus \hat G_1$, where
$\hat G_1 \equiv C(0,\varepsilon)$, and define recursively
$\hat T_{j+1}$ as the first time $\hat\gamma[\hat T_j,\infty)$ hits
$\hat{\mathbb H}_{\hat T_j} \setminus \hat G_{j+1}$, where
$\hat{\mathbb H}_{\hat T_j} \equiv {\mathbb H} \setminus \hat K_{\hat T_j}$ and
$\hat G_{j+1} \equiv \hat g^{-1}_{T_j}(C(\hat g_{\hat T_j}(\hat\gamma(\hat T_j)),\varepsilon))$.
We also define $\hat\tau_{j+1} \equiv \hat T_{j+1} - \hat T_j$, so that
$\hat T_j = \hat\tau_1+\ldots+\hat\tau_j$.
Once again, $\hat G_{j+1}$ is a bounded simply connected domain chosen so that the
conformal transformation which maps $\hat{\mathbb H}_{\hat T_j}$ to $\mathbb H$ maps
$\hat G_{j+1}$ to the semi-ball $C(\hat g_{\hat T_j}(\hat\gamma(\hat T_j)),\varepsilon)$
centered at the point on the real line where the ``tip" $\hat\gamma(\hat T_j)$ of the
hull $\hat K_{\hat T_j}$ is mapped.
The polygonal approximation $\hat\gamma_{\varepsilon}$ is obtained by joining,
for all $j$, $\hat\gamma(\hat T_j)$ to $\hat\gamma(\hat T_{j+1})$ with a straight
segment, where the speed $\hat\gamma_{\varepsilon}'(t)$ is chosen to be constant.

Consider the sequence of times
$\tilde T_j$ defined in the natural way so that
$\tilde\gamma(\tilde T_j) = f(\hat\gamma(\hat T_j))$
and the (discrete-time) stochastic processes
$\hat X_j \equiv (\hat K_{\hat T_j}, \hat\gamma(\hat T_j))$ and
$\tilde X_j \equiv (\tilde K_{\tilde T_j}, \tilde\gamma(\tilde T_j))$
related by $\hat X_j = f^{-1}(\tilde X_j)$.
If for $x \in {\mathbb R}$ we let $\theta[x]$ denote the translation
that maps $x$ to $0$ and define the family of conformal maps
$\tilde f_{\tilde T_j} = \theta[g_{\hat T_j}(\hat\gamma(\hat T_j))] \circ g_{\hat T_j} \circ f^{-1}$
from $D \setminus \tilde K_{\tilde T_j}$ to $\mathbb H$, then
$\tilde f_{\tilde T_j}$ sends $\tilde\gamma(\tilde T_j)$ to $0$ and $b$ to $\infty$,
and $\tilde T_{j+1}$ is the first time $\tilde\gamma[\tilde T_j,\infty)$
hits $\tilde{\mathbb H}_{\hat T_j} \setminus \tilde G_{j+1}$, where
$\tilde{\mathbb H}_{\tilde T_j} = {\mathbb H} \setminus \tilde K_{\hat T_j}$
and $\tilde G_{j+1} = \tilde f_{\tilde T_j}^{-1}(C(0,\varepsilon))$.
Therefore, $\{ \tilde T_j \}$ is a sequence of stopping times like
those used in the definition of the spatial Markov property and, thanks
to the relation $\hat X_j = f^{-1}(\tilde X_j)$, the fact that $\tilde K_t$
satisfies the spatial Markov property implies that $\hat X_j$ is a Markov chain.
We also note that the fact that the hitting distribution of $\tilde\gamma(t)$
is determined by Cardy's formula implies the same for the hitting distribution
of $\hat\gamma(t)$, thanks to the conformal invariance of Cardy's formula.
We next use these properties 
to show that $\hat\gamma_{\varepsilon}$ is distributed like $\gamma_{\varepsilon}$.

To do so, we first note that the conformal transformations $g_{T_j}$ and $\hat g_{\hat T_j}$
are random and that their distributions are functionals of the distributions of the hulls
$K_{T_j}$ and $\hat K_{\hat T_j}$, since there is a one-to-one correspondence between hulls
and conformal maps (with the normalization we have chosen in~(\ref{conf-map1})--(\ref{conf-map2})).
Therefore, since $\hat K_{\hat T_1}$ is distributed like $K_{T_1}$ (see Lemma~\ref{hull}),
$g_{T_1}$ and $\hat g_{\hat T_1}$ have the same distribution, which also implies that
$\hat T_1$ is distributed like $T_1$ because, due to the parametrization by capacity of
$\gamma$ and $\hat\gamma$, $2 T_1$ is exactly the coefficient of the term $1/z$ in the
expansion at infinity of $g_{T_1}$, and $2 \hat T_1$ is exactly the coefficient of the
term $1/z$ in the expansion at infinity of $\hat g_{\hat T_1}$.
Moreover, it is also clear that $\hat\gamma(\hat T_1)$ is distributed like $\gamma(T_1)$,
because their distributions are both determined by Cardy's formula, and so
$\hat g_{\hat T_1}(\hat\gamma(\hat T_1))$ is distributed like $g_{T_1}(\gamma(T_1))$.
Notice that the law of the hull $\hat K_{\hat T_1}$ is conformally invariant because,
by Lemma~\ref{hull}, it coincides with the law of the $SLE_6$ hull $K_{T_1}$.

Using now the Markovian character of $\hat X_j$, which implies that,
conditioned on $\hat X_1=(\hat K_{\hat T_1},\hat\gamma(\hat T_1))$,
$\hat K_{\hat T_2} \setminus \hat K_{\hat T_1}$ and $\hat\gamma(\hat T_2)$
are determined by Cardy's formula in $\hat G_2$, from the fact that
$\hat K_{\hat T_1}$ is equidistributed with $K_{T_1}$ and therefore $\hat G_2$
is equidistributed with $G_2$, we obtain that the hull $\hat K_{\hat T_2}$ is
distributed like $K_{T_2}$ and
its ``tip" $\hat\gamma(\hat T_2)$ is distributed like the ``tip" $\gamma(T_2)$
of the hull $K_{T_2}$, and we can conclude that the joint distribution of
$\{ \hat\gamma(\hat T_1),\hat\gamma(\hat T_2) \}$ is the same as the joint
distribution of $\{ \gamma(T_1),\gamma(T_2) \}$.
It also follows immediately that $\hat g_{\hat T_2}$ is equidistributed
with $g_{T_2}$ and $\hat\tau_2$ is equidistributed with $\tau_2$ or indeed
with $\tau_1$.

By repeating this argument recursively, using at each step the Markovian character
of the hulls and tips, we obtain that, for all $j$, the joint distribution of
$\{ \hat\gamma(\hat T_1),\ldots,\hat\gamma(\hat T_j) \}$ is the same as the joint
distribution of $\{ \gamma(T_1),\ldots,\gamma(T_j) \}$.
This immediately implies that $\hat\gamma_{\varepsilon}$ has the same distribution
as $\gamma_{\varepsilon}$.

In order to conclude the proof, we just have to show that, as $\varepsilon \to 0$,
$\hat\gamma_{\varepsilon}$ converges to $\hat\gamma$ and $\gamma_{\varepsilon}$
to $\gamma$ in the uniform metric on continuous curves~(\ref{distance}).
This, however, follows easily from properties of continuous curves, if we can
show that the time intervals $\hat T_{j+1} - \hat T_j = \hat\tau_{j+1}$ and
$T_{j+1}-T_j=\tau_{j+1}$ go to $0$ as $\varepsilon \to 0$.
To see this, we recall that $\hat\tau_{j+1}$ and $\tau_{j+1}$ are distributed like
$\tau_1$ and use Lemma~2.1 of~\cite{lsw6}, which implies the (deterministic) bound
$\tau_1(\varepsilon) \leq \frac{1}{2} \, \varepsilon^2$, which follows from the
well-know bound $a(t) \leq \varepsilon^2$ for the capacity $a(t) = 2 \, t$
of~(\ref{conf-map2}). \fbox{}

\begin{remark} \label{remark-difference}
The procedure for constructing the polygonal approximations of $\hat\gamma$
and $\gamma$ and the recursive strategy for proving that they have the same
distribution include significant modifications to the sketched argument for
convergence of the percolation exploration process to chordal $SLE_6$ given
by Smirnov in~\cite{smirnov} and \cite{smirnov-long}.
One modification is that we use ``conformal semi-balls" instead of balls
(see~\cite{smirnov,smirnov-long}) to define the sequences of stopping times
$\{ \hat T_j \}$ and $\{ T_j \}$.
Since the paths we are dealing with touch themselves (or almost do), if one
were to use balls, some of them would intersect multiple disjoint pieces of
the past hull, making it impossible to use Cardy's formula in the
``triangular setting" proposed by (Carleson and) Smirnov and used here.
The use of conformally mapped semi-balls ensures, thanks to the choice of the
conformal maps, that the domains used to define the stopping times intersect
a single piece of the past hull.
This seems to be a natural choice (exploiting the conformal invariance)
to obtain a good polygonal approximation of the paths while still being
able to use Cardy's formula to determine hitting distributions.
\end{remark}


We will next prove a version of Smirnov's result (Theorem~\ref{cardy-smirnov}
above) extended to cover the convergence of crossing probabilities to Cardy's
formula for the case of sequences of domains.
The statement of Theorem~\ref{strong-cardy} below is certainly
not optimal, but it is sufficient for our purposes.
We remark that a weaker statement restricted, for instance, only
to Jordan domains would not be sufficient -- see
Figure~\ref{cut-point-fig} and the discussion referring to it in
the proof of Theorem~\ref{spatial-markov} below.
First, we introduce some definitions that will simplify the notation
in the rest of the paper.

We will consider bounded simply connected domains $D$ whose boundaries
$\partial D$ are continuous curves.
Let $f:{\mathbb D} \to D$ be the (unique) conformal map from the
unit disc to $D$ with $f(0)=0$ and $f'(0)>0$; note that by
Theorem~\ref{cont-thm} of Appendix~\ref{rado}, $f$ has a continuous
extension to $\overline{\mathbb D}$.
Let $a,c,d$ be three points of $\partial D$ (or more accurately,
three prime ends) in counterclockwise order -- i.e., such that
$a=f(a^*)$, $c=f(c^*)$ and $d=f(d^*)$, with $a^*$, $c^*$ and $d^*$
in counterclockwise order.
We will call $D$ {\bf admissible} with respect to $(a,c,d)$ if
the counterclockwise arcs $\overline{ac}$, $\overline{cd}$ and
$\overline{da}$ are simple curves, $\overline{cd}$ does not
touch the interior of either $\overline{ac}$ or $\overline{da}$,
and from each point in $\overline{cd}$ there is a path to infinity
that does not cross $\partial D$.

Notice that, according to our definition, the interiors of the
arcs $\overline{ac}$ and $\overline{da}$ can touch.
If that happens, the double-points of the boundary (belonging to
both arcs) are counted twice and considered as two distinct points
(and are two different prime ends).
The significance of the notion of admissible is that certain domains
arising naturally in the proof of Theorem~\ref{spatial-markov}
(and thus in the proof of convergence to $SLE_6$) are not Jordan
but are admissible -- see Figure~\ref{cut-point-fig}.

Consider a sequence of admissible domains $D_k$ with $j$ distinct
selected points $u^1_k, \ldots, u^j_k$ on $\partial D_k$ on each
of their boundaries.
If $D$ is an admissible domain with $j$ selected distinct points
$u^1,\ldots,u^j$ on its boundary such that, as $k \to \infty$,
$\text{d}(\partial D_k,\partial D) \to 0$ and $|u^1_k-u^1|,\ldots,|u^j_k-u^j| \to 0$,
we say that $(D_k,u^1_k,\ldots,u^j_k)$ converges to $(D,u^1,\ldots,u^j)$
and write $(D_k,u^1_k,\ldots,u^j_k) \to (D,u^1,\ldots,u^j)$.


\begin{theorem} \label{strong-cardy}
Consider a sequence $\{(D_k,a_k,b_k,c_k,d_k)\}$ of domains containing the
origin, admissible with respect to $(a_k,c_k,d_k)$, and with $b_k$
belonging to the interior of the counterclockwise arc $\overline{c_k d_k}$.
Assume that $(D_k,a_k,b_k,c_k,d_k)$ converges,
as $k \to \infty$, to $(D,a,b,c,d)$, where $D$ is a domain containing
the origin, admissible with respect to $(a,c,d)$, and $b$ belongs
to the interior of the counterclockwise arc $\overline{cd}$.
Then, for any sequence $\delta_k \downarrow 0$, the probability
$\Phi^{\delta_k}_k (\equiv \Phi^{\delta_k}_{D_k})$ of a blue crossing
inside $D_k$ from the counterclockwise segment $\overline{a_k c_k}$
of $\partial D_k$ to the counterclockwise segment $\overline{b_k d_k}$
of $\partial D_k$ converges, as $k \to \infty$, to Cardy's formula
$\Phi_D$ (see~(\ref{cardy-formula})) for a crossing inside $D$ from
the counterclockwise segment $\overline{a c}$ of $\partial D$ to the
counterclockwise segment $\overline{c d}$ of $\partial D$.
%
\end{theorem}

\noindent {\bf Proof.}
We will construct for each small $\varepsilon>0$,
two domains with boundary points, $(\tilde D, \tilde a, \tilde b, \tilde c, \tilde d)$
and $(\hat D, \hat a, \hat b, \hat c, \hat d)$, approximating $(D,a,b,c,d)$
in such a way that not only does
$\tilde\Phi_{\varepsilon} \equiv \Phi_{\tilde D}(\tilde a, \tilde c; \tilde b, \tilde d)
\stackrel{\varepsilon \to 0}{\longrightarrow} \Phi \equiv \Phi_D(a,c;b,d)$ (by
the continuity of Cardy's formula -- see Lemma~\ref{cont-cardy}) and the same
for $\hat\Phi_{\varepsilon}$, but also so that
\begin{equation} \label{bounds}
\tilde\Phi_{\varepsilon} = \liminf_{k \to \infty} \Phi^{\delta_k}_{\tilde D(\varepsilon)}
\leq \liminf_{k \to \infty} \Phi^{\delta_k}_k \leq \limsup_{k \to \infty} \Phi^{\delta_k}_k
\leq \limsup_{k \to \infty} \Phi^{\delta_k}_{\hat D(\varepsilon)} = \hat\Phi_{\varepsilon}.
\end{equation}
This yields the desired result.
The construction of the approximating domains uses fairly straightforward
conformal mapping arguments.
We provide details for $\tilde D$; the construction of $\hat D$ is analogous.
Before providing the construction details for $\tilde D$, we give several
paragraphs of overview.

\begin{figure}[!ht]
\begin{center}
\includegraphics[width=8cm]{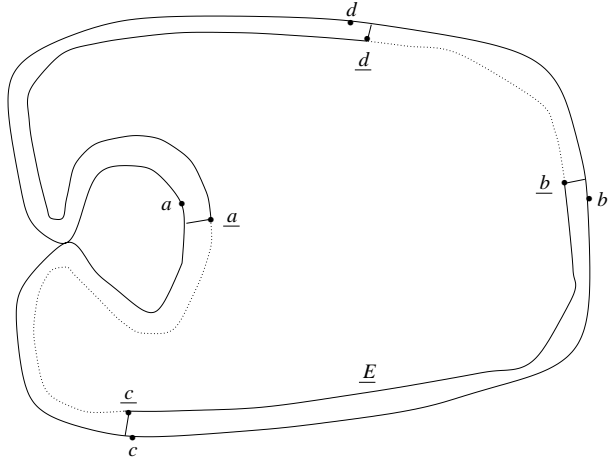}
\caption{Schematic figure for the construction of $\underline E$
and of the new simple loop $\partial G$ of which $\overline E$ will
be an outer approximation.
The outer loop represents $\partial D$ and the inner one $\underline E$.
The dotted portions of $\underline E$ are not used to construct
$\partial G$, which is obtained pasting together the remaining
portions of $\underline E$, portions of $\partial D$ and the
segments connecting $\underline E$ and $\partial D$.
}
\label{fig1-thm6}
\end{center}
\end{figure}

To construct $\tilde D$, we will first need to take an inner approximation
$\underline E = \underline E(D,\varepsilon)$ of $\partial D$, where
$\underline E$ is a simple loop surrounded by $\partial D$ and with
$\text{d}(\partial D, \underline E) \leq \varepsilon$.
We will then construct a Jordan domain $G(\varepsilon)$ whose boundary
$\partial G$ is composed of pieces of $\partial D$ and of $\underline E$,
plus four segments joining $\partial D$ with $\underline E$
(see Figure~\ref{fig1-thm6}), as we explain below.

Next, we will take an outer approximation $\overline E = \overline E(G,\varepsilon)$
of $\partial G$, where $\overline E$ is a simple loop surrounding $\partial G$
and with $\text{d}(\partial G, \overline E) \leq \varepsilon$.
We will also need four simple curves $\{ \partial_a,\partial_b,\partial_c,\partial_d \}$
in the interior of the (topological) annulus between $\underline E$ and $\overline E$
and connecting their endpoints
$\{ (\underline a,\overline a),(\underline b,\overline b),(\underline c,\overline c),
(\underline d,\overline d) \}$ on $\underline E$ and $\overline E$ with each of the
four curves touching $\partial D$ at exactly one point which is either in the interior
of the counterclockwise segment $\overline{ac}$ (for $\partial_a$ and $\partial_c$) or
else the counterclockwise segment $\overline{bd}$ (for $\partial_b$ and $\partial_d$).
Furthermore each of these connecting curves is close to its corresponding point $a,b,c$,
or $d$; i.e., $\text{d}(\partial_a,a) \leq \varepsilon$, etc.
(see Figure~\ref{fig1-thm6} where each connecting segment represents ``half'' of one
of these connecting curves).

We will take $\tilde a = \overline a$, $\tilde b = \overline b$, $\tilde c = \overline c$,
and $\tilde d = \overline d$ with $\partial \tilde D$ the concatenation of:
$\partial_a$ from $\underline a$ to $\overline a$, the portion of $\overline E$
from $\overline a$ to $\overline c$ counterclockwise, $\partial_c$ from $\overline c$
to $\underline c$, the portion of $\underline E$ from $\underline c$ to $\underline b$
counterclockwise, $\partial_b$ from $\underline b$ to $\overline b$, the portion
of $\overline E$ from $\overline b$ to $\overline d$ counterclockwise, $\partial_d$
from $\overline d$ to $\underline d$, and the portion of $\underline E$ from
$\underline d$ to $\underline a$ counterclockwise.
It is important that (for fixed $\varepsilon$ and $\tilde D$) there is a strictly
positive minimal distance between $\underline E$ and $\partial D$, and between
$\partial\tilde D$ and the union of the two counterclockwise segments $\overline{cb}$
and $\overline{da}$ of $\partial D$.
These features will guarantee (as we explain with more detail below) that for fixed
$\varepsilon$, once $k$ is large enough, a continuous curve within $\tilde D$ that
corresponds to the crossing event whose probability is
$\Phi_{\tilde D}(\tilde a,\tilde c;\tilde b,\tilde d)$ must have a subpath
corresponding to the crossing event in $D_k$ whose probability is
$\Phi_{D_k}(a_k,c_k;b_k,d_k)$.
This is the key feature of $\tilde D$, which will yield the first inequality
of~(\ref{bounds}).

\begin{figure}[!ht]
\begin{center}
\includegraphics[width=6cm]{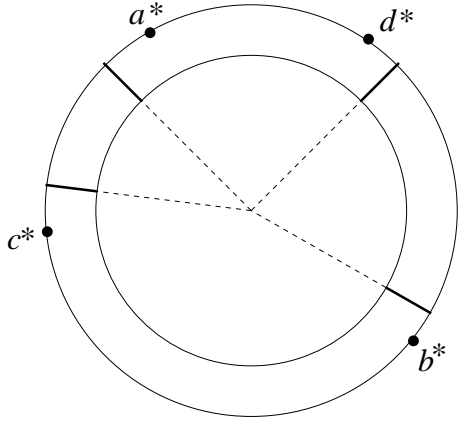}
\caption{The figure shows $\partial^*_a, \partial^*_b, \partial^*_c, \partial^*_d$
represented as heavy segments between the unit circle and the circle of radius
$1-\varepsilon'$ near $a^*,b^*,c^*,d^*$.}
\label{fig2-thm6}
\end{center}
\end{figure}

\begin{figure}[!ht]
\begin{center}
\includegraphics[width=6cm]{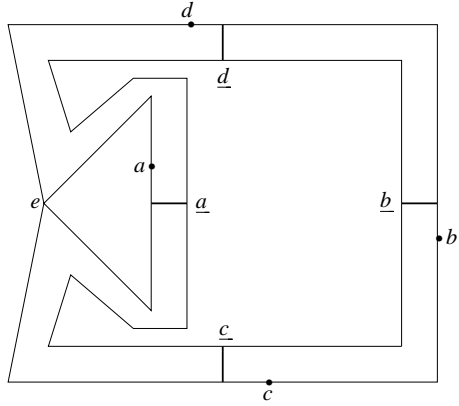}
\caption{The outer loop is the boundary $\partial D$ of a bounded
simply connected domain $D$.
The boundary has a double point at $e$.
The inner loop is $\underline E = f(\partial {\mathbb D}_{1-\varepsilon'})$,
where $f$ is the conformal map from $\mathbb D$ to $D$.
The four segments are the images under $f$ of the segments
$\partial^*_a, \partial^*_b, \partial^*_c, \partial^*_d$
of Figure~\ref{fig2-thm6}.}
\label{inner-bd}
\end{center}
\end{figure}

We now give a more detailed explanation of the construction of $\tilde D$.
We first construct the parts of $\partial\tilde D$ that are inside $\partial D$.
Let $f$ be the conformal map from $\mathbb D$ onto $D$ with $f(0)=0$ and $f'(0)>0$,
and consider the image $f(\partial{\mathbb D}_{1-\varepsilon'})$ of the circle
$\partial{\mathbb D}_{1-\varepsilon'} = \{ z : |z| = 1 - \varepsilon' \}$ under $f$
and the inverse images, $a^*,b^*,c^*,d^*$, under $f^{-1}$ of $a,b,c,d$.
Let $\partial^*_a(\varepsilon',\varepsilon_a)$ be the straight segment between
$e^{-i \varepsilon_a} a^*$ on the unit circle $\partial {\mathbb D}$, and
$(1-\varepsilon') e^{-i \varepsilon_a} a^*$ on the circle
$\partial {\mathbb D}_{1-\varepsilon'}$, and define $\partial^*_b,\partial^*_c$,
and $\partial^*_d$ similarly,  but using \emph{clockwise} rotations by
$e^{+i \varepsilon_c}$ and $e^{+i \varepsilon_d}$ for $c^*$ and $d^*$
(see Figure~\ref{fig2-thm6}).
$f(\partial {\mathbb D}_{1-\varepsilon'})$ is a candidate for $\underline E(D,\varepsilon)$
and $f(\partial^*_{\sharp}(\varepsilon',\varepsilon_{\sharp}))$ is a candidate for half of
$\partial_{\sharp}$ (where ${\sharp}=a$ or $b$ or $c$ or $d$), so we must choose
$\varepsilon'$ and the $\varepsilon_{\sharp}$'s small enough so that
$\text{d}(\partial D, f({\mathbb D}_{1-\varepsilon'})) \leq \varepsilon$,
$\text{d}(f(\partial^*_a(\varepsilon',\varepsilon_a)),a) \leq \varepsilon$, etc.
We then define $\underline a = f((1-\varepsilon')e^{-i \varepsilon_a} a^*)$
and similarly for $\underline b$, $\underline c$ and $\underline d$
(see Figure~\ref{inner-bd}).

\begin{figure}[!ht]
\begin{center}
\includegraphics[width=6cm]{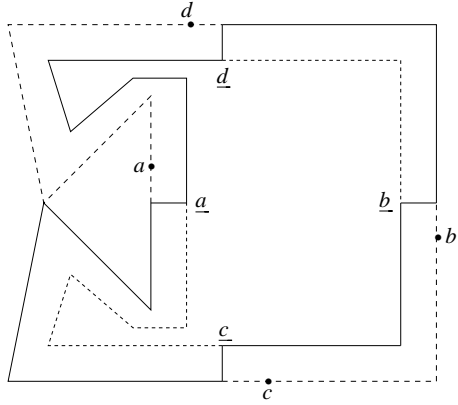}
\caption{The pieces of $\partial D$ (heavier lines) and
of $\underline E$ used in the construction of $\partial G$
are indicated by full lines, the ones that are not used
by dashed lines.
The boundary $\partial G$ of the new Jordan domain $G$
is the resulting full line loop.}
\label{intermediate1}
\end{center}
\end{figure}

Consider now the Jordan domain $G=G(\varepsilon)$ whose boundary $\partial G$
is given by the concatenation of: $f(\partial^*_a)$ from $\underline a$ to
$f(e^{-i \varepsilon_a} a^*)$, the portion of $\partial D$ from
$f(e^{-i \varepsilon_a} a^*)$ to $f(e^{i \varepsilon_c} c^*)$ counterclockwise,
$f(\partial^*_c)$ from $f(e^{i \varepsilon_c} c^*)$ to $\underline{c}$, the
portion of $\underline E$ from $\underline c$ to $\underline{b}$ counterclockwise,
$f(\partial^*_b)$ from $\underline b$ to $f(e^{-i \varepsilon_b} b^*)$, the portion
of $\partial D$ from $f(e^{-i \varepsilon_b} b^*)$ to $f(e^{i \varepsilon_d} d^*)$
counterclockwise, $f(\partial^*_d)$ from $f(e^{i \varepsilon_d} d^*)$ to
$\underline d$, and the portion of $\underline E$ from $\underline d$ to
$\underline a$ counterclockwise (see Figure~\ref{intermediate1}).

\begin{figure}[!ht]
\begin{center}
\includegraphics[width=6cm]{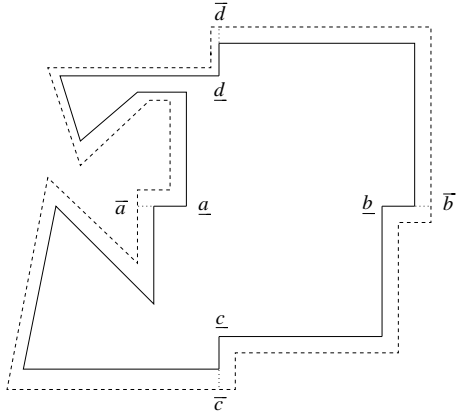}
\caption{The inner (full) loop is $\partial G$, the outer (dashed)
one is its outer approximation $\overline E$. The four dotted segments
between $\partial G$ and $\overline E$, ending at $\overline a$,
$\overline b$, $\overline c$ and $\overline d$, are the other halves
of $\partial_a$, $\partial_b$, $\partial _c$ and $\partial_d$,
respectively.}
\label{intermediate2}
\end{center}
\end{figure}

The exterior of this new Jordan domain $G$ is a connected domain for which
we can do a construction analogous to the one for the original domain $D$
using a conformal map from $\mathbb D$ to obtain candidates for
$\overline E(D,\varepsilon)$ and for the exterior halves of the
$\partial_{\sharp}$'s.
To do this, we use a conformal map from $\mathbb D$ to the exterior of
$\partial G$ and use $\underline a,\underline b,\underline c,\underline d$
as replacements for $a,b,c,d$ (see Figure~\ref{intermediate2}).

\begin{figure}[!ht]
\begin{center}
\includegraphics[width=6cm]{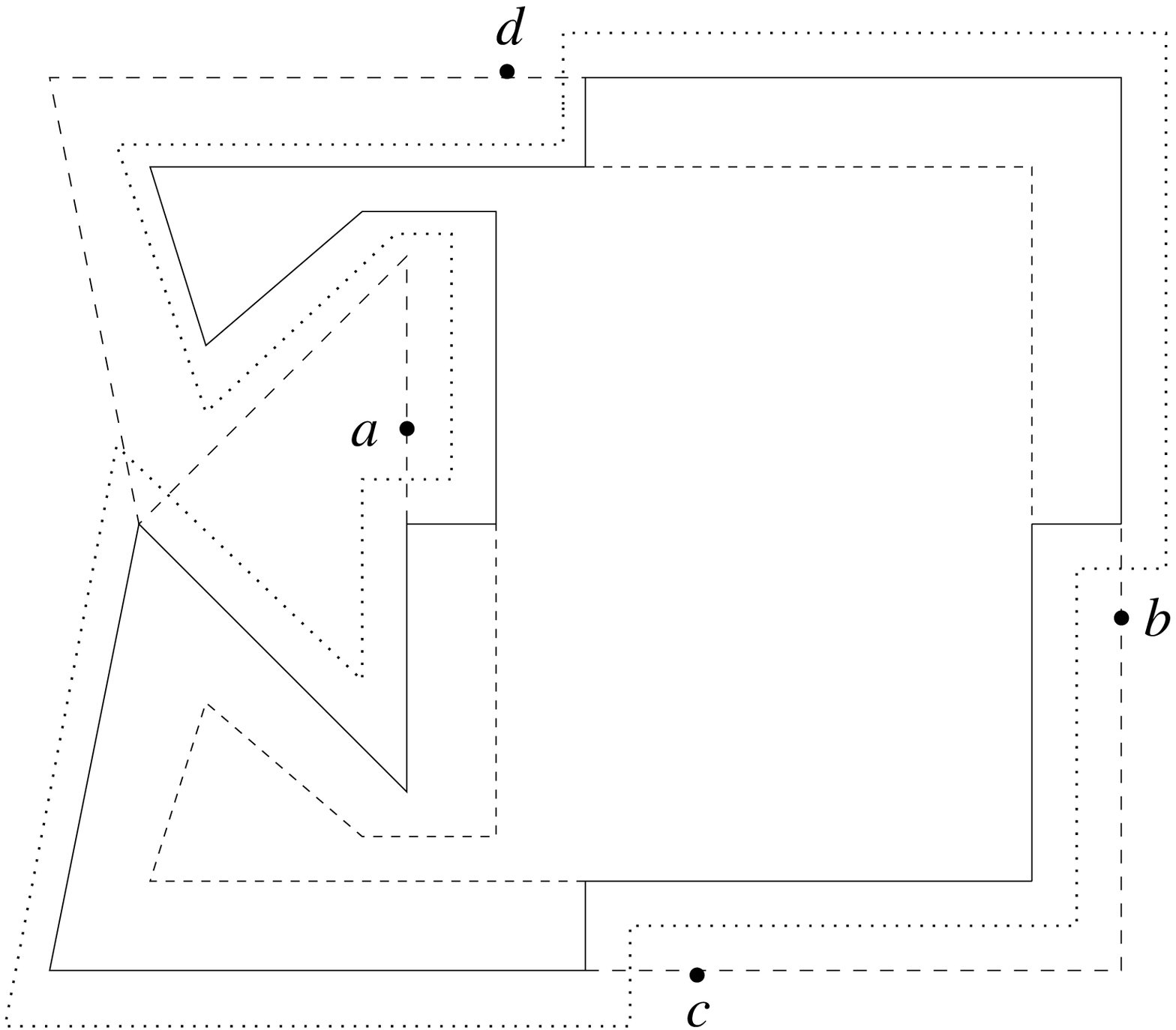}
\caption{The outer approximation $\overline E$ (dotted loop)
of $\partial G$ is chosen close enough to $\partial G$
(full loop) so that it stays between $a$ and $f(e^{-i \varepsilon_a} a^*)$,
between $b$ and $f(e^{-i \varepsilon_b} b^*)$, between $c$ and
$f(e^{i \varepsilon_c} c^*)$, and between $d$ and
$f(e^{i \varepsilon_d} d^*)$.}
\label{outer-bd-fig}
\end{center}
\end{figure}

Finally, we use the freedom to choose the exterior replacements for
$\varepsilon'$ and the $\varepsilon_{\sharp}$'s differently from the
interior values to make sure that the interior and exterior halves
of the $\partial_{\sharp}$'s match up.
We also choose the exterior values for $\varepsilon'$ small enough so
that $\overline E$ stays between $a$ and $f(e^{-i \varepsilon_a} a^*)$,
between $b$ and $f(e^{-i \varepsilon_b} b^*)$, between $c$ and
$f(e^{i \varepsilon_c} c^*)$, between $d$ and $f(e^{i \varepsilon_d} d^*)$
(see Figure~\ref{outer-bd-fig}).

Once all the pieces of $\partial\tilde D(\varepsilon)$ are available,
they are put together as explained above, and as show in
Figures~\ref{intermediate3} and~\ref{fig3-thm6} below.

\begin{figure}[!ht]
\begin{center}
\includegraphics[width=6cm]{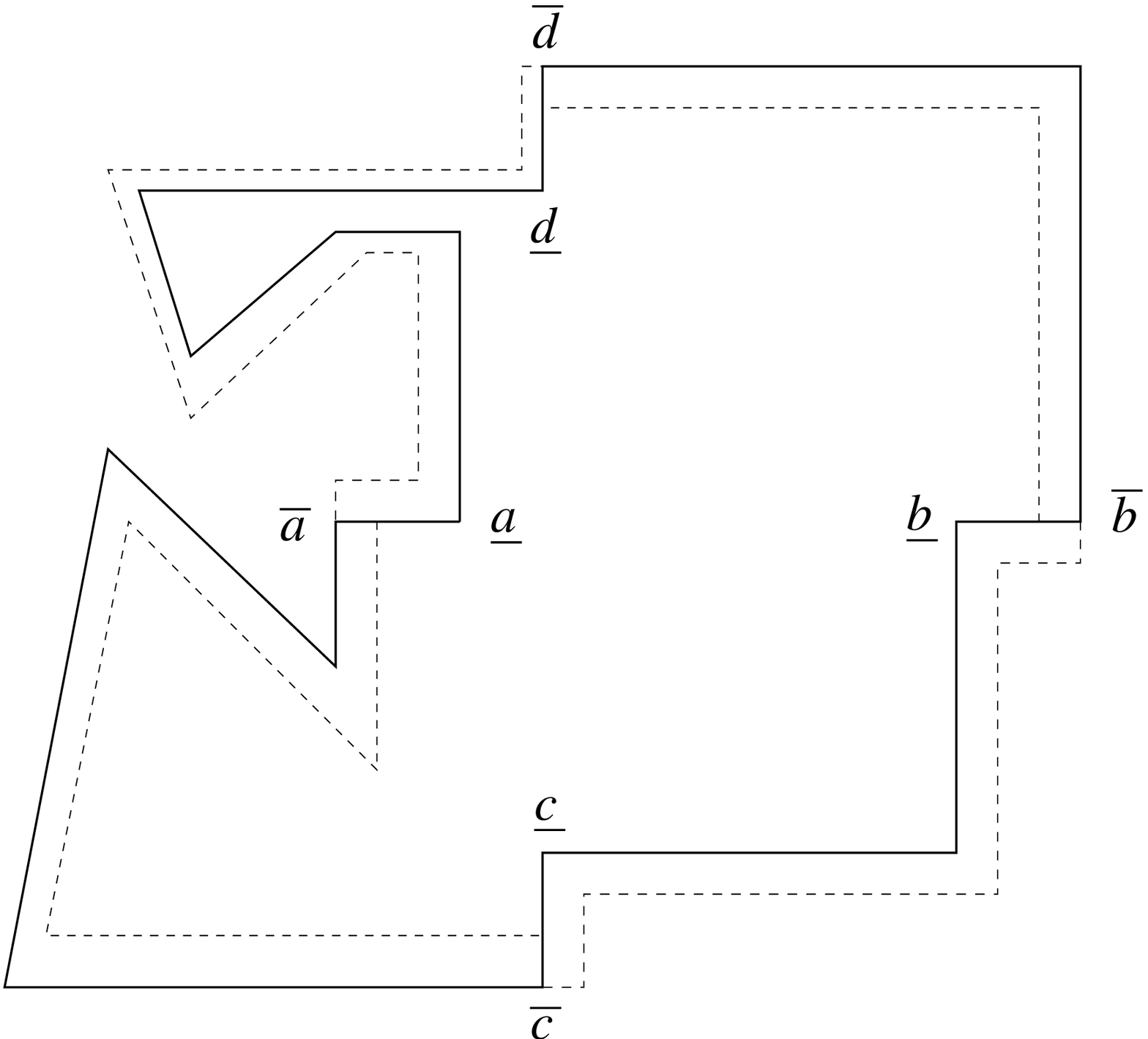}
\caption{Final step of the construction of the simple loop $\partial\tilde D$
(shown as a full line loop) using pieces of $\partial G$ and of $\overline E$
(see Figure~\ref{intermediate2}).
The pieces of $\partial G$ and of $\overline E$ that are not used are dashed.}
\label{intermediate3}
\end{center}
\end{figure}

\begin{figure}[!ht]
\begin{center}
\includegraphics[width=6cm]{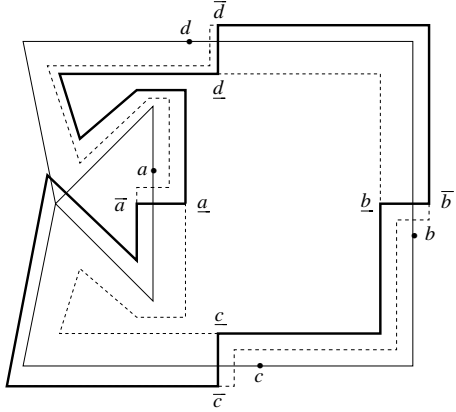}
\caption{The figure shows all the curves and special points involved
in the construction of the simple loop $\partial\tilde D$, whose steps
have been detailed in the text.
$\partial\tilde D$ is indicated by a heavy full line.}
\label{fig3-thm6}
\end{center}
\end{figure}

It should be clear that for a given approximation $\partial\tilde D$ of $\partial D$
constructed as described above there is a strictly positive $\tilde\varepsilon$ such
that the distance between $\partial D$ and the portions of $\partial\tilde D$ that
belong to $\underline E$ and $\overline E$ is not smaller than $\tilde\varepsilon$,
and the distance between the union $\overline{cb} \cup \overline{da}$, two
counterclockwise segments of $\partial D$, and
$\partial_a \cup \partial_b \cup \partial_c \cup \partial_d$ is also not smaller
than $\tilde\varepsilon$.
On the other hand, for any $\tilde\varepsilon>0$, there exists
$k_0=k_0(\tilde\varepsilon)$ such that for all $k \geq k_0$, $\partial D_k$ is
contained inside the $\tilde\varepsilon$-neighborhood of $\partial D$ with
the counterclockwise segment $\overline{c_k b_k}$ (resp., $\overline{d_k a_k}$)
in the $\tilde\varepsilon$-neighborhood of the counterclockwise segment
$\overline{cb}$ (resp., $\overline{da}$).
This implies that for $k$ large enough, any blue path crossing inside $\tilde D$ from
the counterclockwise segment $\overline{\overline a \, \overline c}$ of $\partial \tilde D$
to the counterclockwise segment $\overline{\overline b \, \overline d}$ of $\partial \tilde D$
must have a subpath that stays inside $D_k$ and crosses between the counterclockwise
segment $\overline{a_k c_k}$ of $\partial D_k$ and the counterclockwise segment
$\overline{b_k d_k}$ of $\partial D_k$.
Therefore the crossing probability $\tilde \Phi^{\delta_k}_{\tilde D(\varepsilon)}$
is a lower bound for $\Phi^{\delta_k}_k$ for all $k \geq k_0$, so that
\begin{equation} \label{bound}
\liminf_{k \to \infty} \Phi^{\delta_k}_k \geq
\lim_{k \to \infty} \tilde\Phi^{\delta_k}_{\tilde D(\varepsilon)} = \tilde\Phi_{\varepsilon},
\end{equation}
as desired (the equality uses Smirnov's result, Theorem~\ref{cardy-smirnov},
for fixed $\tilde D(\varepsilon)$).

We now note that as $\varepsilon \to 0$,
$(\tilde D, \tilde a, \tilde b, \tilde c, \tilde d) \to (D,a,b,c,d)$.
This allows us to use the continuity of Cardy's formula
(Lemma~\ref{cont-cardy} in Appendix~\ref{rado}) to obtain
\begin{equation} \label{cardy2}
\lim_{\varepsilon \to 0} \tilde\Phi_{\varepsilon} = \Phi.
\end{equation}
From this and~(\ref{bound}) it follows that
\begin{equation} \label{upper}
\liminf_{k \to \infty} \Phi^{\delta_k}_k \geq \Phi.
\end{equation}

The remaining part of the proof involves defining a domain $\hat D$
analogous to $\tilde D$ but with the property that the probability
$\hat\Phi^{\delta_k}_{\varepsilon}$ of an appropriate crossing,
such that
\begin{equation}
\lim_{\varepsilon \to 0} \lim_{k \to \infty}
\hat \Phi^{\delta_k}_{\hat D(\varepsilon)} = \Phi,
\end{equation}
is an upper bound for $\Phi^{\delta_k}_k$ for all $k$ large enough.
(The details of the construction of $\hat D$ are analogous to
those of $\tilde D$; we leave them to the reader.)
This shows that
\begin{equation} \label{lower}
\limsup_{k \to \infty} \Phi^{\delta_k}_k \leq \Phi,
\end{equation}
which, combined with~(\ref{upper}), implies
\begin{equation}
\lim_{k \to \infty} \Phi^{\delta_k}_k = \Phi
\end{equation}
and concludes the proof. \fbox{} \\

We note that Theorem~\ref{strong-cardy}, combined with the continuity
of Cardy's formula in the shape of the domain (for admissible domains)
and positions of the four
points on the boundary (Lemma~\ref{cont-cardy}), implies the convergence
of crossing probabilities to Cardy's formula \emph{locally uniformly}
in the shape of the domain with respect to the uniform metric on curves,
and in the location of the four points on the boundary with respect to
the Euclidean metric; i.e., for $(D,a,b,c,d)$ an admissible domain with
$a,b,c,d \in \partial D$ (with the notation used in Theorem~\ref{strong-cardy}),
$\forall \varepsilon>0$, $\exists \alpha_0=\alpha_0(\varepsilon)$ and
$\delta_0=\delta_0(\varepsilon)$ such that for all admissible domains
$(D',a',b',c',d')$ with
$\max{(\text{d}(\partial D, \partial D'),|a-a'|,|b-b'|,|c-c'|,|d-d'|) \leq \alpha_0}$
and $\delta \leq \delta_0$,
$|\Phi_{D'}(a',c';b',d') - \Phi^{\delta}_{D'}(a',c';b',d')| \leq \varepsilon$,
where $\Phi_{D'}(a',c';b',d')$ is Cardy's formula and
$\Phi^{\delta}_{D'}(a',c';b',d')$ is the corresponding crossing probability.

Our next task is to show that the filling of any subsequential scaling limit
of the percolation exploration process satisfies the spatial Markov property.
Let us start with some notation.
First of all, suppose that $D$ is a Jordan domain, and $a$ and $b$ are
two points on $\partial D$.
Define the {\bf discrete filling} (or simply filling) at time $t$ of a
percolation exploration path $\gamma^{\delta}_{D,a,b}$ (with a given
parametrization) inside (the $\delta$-approximation $D^{\delta}$ of) $D$
from (the e-vertex closest to) $a$ to (the e-vertex closest to) $b$ to be
the union of the hexagons explored up to time $t$ and those unexplored
hexagons from which it is not possible to reach $b$ without crossing an
explored hexagon or $\partial D$ (in other words, this is the set of
hexagons that at time $t$ have been explored or are disconnected from
$b$ by the exploration path).

Consider a sequence $\{(D_k,a_k,b_k)\}$ of Jordan domains such that
$(D_k,a_k,b_k) \to (D,a,b)$, where $a_k,b_k \in \partial D_k$ are
two distinct points on $\partial D_k$.
Denote by $\gamma^{\delta}_k \equiv \gamma^{\delta}_{D_k,a_k,b_k}$
the percolation exploration path inside (the $\delta$-approximation
$D^{\delta}_k$ of) $D_k$ from (the e-vertex closest to) $a_k$ to
(the e-vertex closest to) $b_k$.

Notice that we can couple the paths $\gamma^{\delta}_k$ simultaneously for
all values of $k$ and $\delta$ by using the same percolation configuration
to generate all of them.
We can then apply the results of~\cite{ab} to conclude that there exists a
subsequence $\delta_k \downarrow 0$ such that the law of $\gamma^{\delta_k}_k$
converges to some limiting law for a process $\tilde\gamma$ supported
on (H\"older) continuous curves inside $D$ from $a$ to $b$.
The filling $\tilde K_t$ of $\tilde\gamma[0,t]$, appearing in the next
theorem, is defined just above Equation~(\ref{hulls}).

%

\begin{theorem} \label{spatial-markov}
For any subsequential limit process $\tilde\gamma$ of the percolation
exploration path $\gamma^{\delta}_k$ defined above, the filling
$\tilde K_t$ of $\tilde\gamma[0,t]$, as a process, satisfies the spatial
Markov property.
\end{theorem}

\noindent {\bf Proof.} Let $\delta_k \downarrow 0$ be a subsequence such
that the law of $\gamma^{\delta_k}_k$ converges to some limiting law
supported on continuous curves $\tilde\gamma$ inside $D$ from $a$ to $b$.
We will prove the spatial Markov property by showing that
$(\tilde K_{\tilde T_j}, \tilde\gamma(\tilde T_j))$ as defined in the proof
of Theorem~\ref{characterization} are \emph{jointly} distributed like the
corresponding $SLE_6$ hull variables, which do have the spatial Markov property.
For each $k$, let $K^k_t$ denote the filling at time $t$ of $\gamma^{\delta_k}_k$
(with some parametrization -- we do not need to worry about the choice of
parametrization here).
It follows from the Markovian character of the percolation exploration
process that, for all $k$, the filling $K^k_t$ of the percolation exploration
path $\gamma^{\delta_k}_k$ satisfies a suitably adapted (to the discrete
setting) spatial Markov property.
(In fact, the percolation exploration path satisfies a stronger property -- roughly
speaking, that the future of the path given the filling of the past is distributed
as a percolation exploration path in the original domain from which the filling
of the past has been removed.)

To be more precise, let $f^k_0$ be a conformal transformation that maps $D_k$
to $\mathbb H$ such that $f^k_0(a_k)=0$ and $f^k_0(b_k) = \infty$ and let
$T^k_1=T^k_1(\varepsilon)$ denote the first exit time of
$\gamma^{\delta_k}_k(t)$ from $G^k_1 \equiv (f_0^k)^{-1}(C(0,\varepsilon))$ defined
as the first time an explored hexagon intersects the image under $(f_0^k)^{-1}$
of the semi-circle $\{ z : |z| = \varepsilon \} \cap {\mathbb H}$.
Define recursively $T^k_{j+1}$ as the first exit time of $\gamma^{\delta_k}_k[T^k_j,\infty)$
from $G^k_{j+1} \equiv (f^k_{T^k_j})^{-1}(C(0,\varepsilon))$, where $f^k_{T^k_j}$
is a conformal map from $D_k \setminus K^k_{T^k_j}$ to $\mathbb H$ that maps
$\gamma^{\delta_k}_k(T^k_j)$ to $0$ and $b_k$ to $\infty$.
We also define $\tau^k_{j+1} \equiv T^k_{j+1} - T^k_j$, so that
$T^k_j=\tau^k_1+\ldots+\tau^k_j$, and the (discrete-time) stochastic process
$X^k_j \equiv (K^k_{T^k_j},\gamma^{\delta_k}_k(T^k_j))$ for $j=1,2,\ldots$ .
The Markovian character of the percolation exploration process implies that,
for every $k$, $X^k_j$ is a Markov process (in $j$).
In order to study the limit as $k \to \infty$ of $X^k_1,X^k_2,\ldots$,
we first need to analyze in more detail the mappings $f^k_{T^k_j}$.

The conformal transformations $f^k_{T^k_j}$ are defined as follows, with
an arbitrary multiplicative factor $\lambda_j$.
We choose
$f^k_{T^k_j}$ to be the composition $\lambda_j \, \psi^k_j \circ \phi^k_j$
of two maps, where $\lambda_j>0$, $\phi^k_j$ is the conformal transformation
that maps $D_k \setminus K^k_{T^k_j}$ onto $\mathbb D$ with $\phi^k_j(0)=0$
and $(\phi^k_j)'(0)>0$ (we are assuming for simplicity that the domain
$D_k \setminus K^k_{T^k_j}$ contains the origin; if that is not the case,
one can think of a translated domain that does contain the origin),
while $\psi^k_j$ is the inverse of the transformation
\begin{equation}
w = e^{i \theta^k_j} \left( \frac{(z+1)-z^k_j}{(z+1)-\overline{z^k_j}} \right)
\end{equation}
that maps $\overline{\mathbb H}$ onto $\overline{\mathbb D}$, where $\theta^k_j$
is chosen so that $e^{i \theta^k_j} = \phi^k_j(b_k)$ and $z^k_j$ can be chosen
so that $|1-z^k_j|=1$, $\text{Im}(z^k_j)>0$ and
$\phi^k_j(b_k) (\frac{1-z^k_j}{1-\overline{z^k_j}}) = \phi^k_j(\gamma^{\delta_k}_k(T^k_j))$,
which means that $\lambda_j \, \psi^k_j$ maps $\phi^k_j(\gamma^{\delta_k}_k(T^k_j))$ to $0$
and $\phi^k_j(b_k)$ to $\infty$, so that $f^k_{T^k_j} = \lambda_j \, \psi^k_j \circ \phi^k_j$
indeed maps $\gamma^{\delta_k}_k(T^k_j)$ to $0$ and $b_k$ to $\infty$.

Since $\gamma^{\delta_k}_k$ converges in distribution to $\tilde\gamma$, we
can find two coupled versions of $\gamma^{\delta_k}_k$ and $\tilde\gamma$
on some probability space $(\Omega',{\cal B}',{\mathbb P}')$ such that
$\gamma^{\delta_k}_k$ converges to $\tilde\gamma$ for all $\omega' \in \Omega'$;
in the rest of the proof we work with these new versions which, with a slight
abuse of notation, we denote with the same names as the original ones.
Let $\tilde f_0$ be a conformal transformation that maps $D$ to
$\mathbb H$ such that $\tilde f_0(a)=0$ and $\tilde f_0(b) = \infty$ and
let $\tilde T_1=\tilde T_1(\varepsilon)$ denote the first time $\tilde\gamma(t)$
hits $D \setminus \tilde G_1$, with $\tilde G_1 \equiv \tilde f_0^{-1}(C(0,\varepsilon))$.
Define recursively $\tilde T_{j+1}$ as the first time $\tilde\gamma(t)$ hits
$D \setminus \tilde G_{j+1}$, with
$\tilde G_{j+1} \equiv \tilde f_{\tilde T_j}^{-1}(C(0,\varepsilon))$, where
$\tilde f_{\tilde T_j}$ is a conformal map from $D \setminus \tilde K_{\tilde T_j}$
to $\mathbb H$ that maps $\tilde\gamma(\tilde T_j)$ to $0$ and $b$ to $\infty$.
We also define $\tilde\tau_j \equiv \tilde T_{j+1} - \tilde T_j$, so that
$\tilde T_j = \tilde\tau_1+\ldots+\tilde\tau_j$, and the (discrete-time)
stochastic process $\tilde X_j \equiv (\tilde K_{\tilde T_j},\tilde\gamma(\tilde T_j))$.
As above, we choose the conformal transformation $\tilde f_{\tilde T_j}$ to be the
composition $\lambda_j \, \tilde\psi_j \circ \tilde\phi_j$ of two maps, where $\lambda_j>0$,
$\tilde\phi_j$ is the conformal transformation that maps $D \setminus \tilde K_{T_j}$
onto $\mathbb D$ with $\tilde\phi_j(0)=0$ and $\tilde\phi'_j(0)>0$ (once again, we
are assuming for simplicity that the domain $D \setminus \tilde K_{T_j}$ contains
the origin; if that is not the case, one can think of a translated domain that does
contain the origin), while $\tilde\psi_j$ is the inverse of the transformation
\begin{equation}
w = e^{i \tilde\theta_j} \left( \frac{(z+1)-\tilde z_j}{(z+1)-\overline{\tilde z_j}} \right)
\end{equation}
that maps $\overline{\mathbb H}$ onto $\overline{\mathbb D}$, where $\tilde\theta_j$
is chosen so that $e^{i \tilde\theta_j} = \tilde\phi_j(b)$ and $\tilde z_j$ can
be chosen so that $|1-\tilde z_j|=1$, $\text{Im}(\tilde z_j)>0$ and
$\tilde\phi_j(b) (\frac{1-\tilde z_j}{1-\overline{\tilde z_j}}) = \tilde\phi_j(\tilde\gamma(\tilde T_j))$,
which means that $\lambda_j \, \tilde\psi_j$ maps $\tilde\phi_j(\tilde\gamma(\tilde T_j))$
to $0$ and $\tilde\phi_j(b)$ to $\infty$, so that
$\tilde f_{\tilde T_j} = \lambda_j \, \tilde\psi_j \circ \tilde\phi_j$
indeed maps $\tilde\gamma(\tilde T_j)$ to $0$ and $b$ to $\infty$.

Analogous quantities can be defined for the trace of chordal $SLE_6$.
For clarity, they will be indicated here by the superscript $SLE_6$;
e.g., $f_j^{SLE_6}$, $K_{T_j}^{SLE_6}$, $G_j^{SLE_6}$ and $X_j^{SLE_6}$.
%
%
We want to show recursively that, for any $j$, as $k \to \infty$,
$\{ X^k_1,\ldots,X^k_j \}$ converge jointly in distribution to
$\{ \tilde X_1,\ldots,\tilde X_j \}$.
By recursively applying Theorem~\ref{strong-cardy} and Lemma~\ref{hull},
we can then conclude that $\{ \tilde X_1,\ldots,\tilde X_j \}$ are
jointly equidistributed with the corresponding $SLE_6$ hull variables
(at the corresponding stopping times) $\{ X_1^{SLE_6},\ldots,X_j^{SLE_6} \}$.
Since the latter do satisfy the spatial Markov property, so will the
former, as desired.

The zeroth step consists in noticing that the convergence of $(D_k,a_k,b_k)$
to $(D,a,b)$ as $k \to \infty$ allows us to apply Rad\'o's theorem
(i.e., Theorem~\ref{rado-thm} of Appendix~\ref{rado}) to show that $(\phi^k_0)^{-1}$
converges to $\tilde\phi_0^{-1}$ uniformly in $\overline{\mathbb D}$.
This, together with the convergence of $a_k$ to $a$ and $b_k$ to $b$, implies
that $\phi^k_0(a_k)$ converges to $\tilde\phi_0(a)$ and $\phi^k_0(b_k)$ to
$\tilde\phi_0(b)$.
Therefore, we also have the convergence of $\lambda_0 \, \psi^k_0$ to
$\lambda_0 \, \tilde\psi_0$ and we can conclude that $(f^k_0)^{-1}$ converges
to $\tilde f_0^{-1}$ uniformly on compact subsets of $\overline{\mathbb H}$,
which implies that the boundary
$\partial G^k_1$ of $G^k_1=(f^k_0)^{-1}(C(0,\varepsilon))$ converges to the
boundary $\partial\tilde G_1$ of $\tilde G_1 = \tilde f_0^{-1}(C(0,\varepsilon))$
in the uniform metric on continuous curves.

Starting from there, the first step of our recursion argument is organized
as follows:
\begin{itemize}
\item[(1)] $K^k_{T^k_1} \to \tilde K_{\tilde T_1}$ by ``number of arms"
percolation bounds~\cite{ksz} and Lemma~\ref{double-crossing} below,
but also $K^k_{T^k_1} \to K_{T_1}^{SLE_6}$ by Lemma~\ref{hull}
(Theorem~\ref{strong-cardy} is used here).
\item[(2)] $D_k \setminus K^k_{T^k_1} \to D \setminus \tilde K_{\tilde T_1}$,
but also $D_k \setminus K^k_{T^k_1} \to D \setminus K_{T_1}^{SLE_6}$, by (1).
\item[(3)] $f^k_{T^k_1} \to \tilde f_{\tilde T_1}$, but also
$f^k_{T^k_1} \to f_{T_1}^{SLE_6}$, by Corollary~\ref{cor-unif-conv}.
\item[(4)] $G^k_2 \to \tilde G_2$, but also $G^k_2 \to G_2^{SLE_6}$, by (3).
\end{itemize}
At this point, we are in the same situation as at the zeroth step,
but with $G^k_1$, $\tilde G_1$ and $G_1^{SLE_6}$ replaced by $G^k_2$,
$\tilde G_2$ and $G_2^{SLE_6}$ respectively, and we can proceed by
recursion.
As explained above, the theorem then follows from the fact that
the $SLE_6$ hull variables do posses the spatial Markov property.

In what follows we will show that the ``number of arms"
bounds~\cite{ksz} and Lemma~\ref{double-crossing} below imply the
convergence of $K^k_{T^k_j}$ to $\tilde K_{\tilde T_j}$, and that
the conditions to apply Theorem~\ref{strong-cardy}, Lemma~\ref{double-crossing}
and Corollary~\ref{cor-unif-conv} are always satisfied.
(This last point boils down to showing that the domains
$D \setminus \tilde K_{\tilde T_j}$ and $\tilde G_j$ are
admissible for all $j$.)
%
%
%
We begin by showing first that, as $k \to \infty$, $K^k_{T^k_1}$
converges in distribution to $\tilde K_{\tilde T_1}$ (which also
implies that $\gamma^{\delta_k}_k(T^k_1)$ converges in distribution
to $\tilde\gamma(\tilde T_1)$), from which it follows that $X^k_1$
converges in distribution to $\tilde X_1$.

\begin{figure}[!ht]
\begin{center}
\includegraphics[width=8cm]{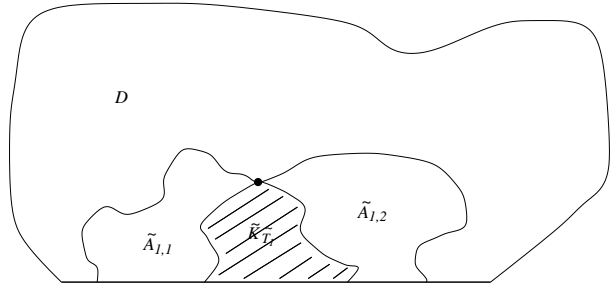}
\caption{Schematic figure representing
$\tilde G_1 \setminus \tilde K_{\tilde T_1} = \tilde A_{1,1} \cup \tilde A_{1,2}$.}
\label{fig1-thm7}
\end{center}
\end{figure}

Consider $G^k_1 \setminus K^k_{T^k_1}$ and $\tilde G_1 \setminus \tilde K_{\tilde T_1}$;
they are both composed of two domains (which ``meet" at $\gamma^{\delta_k}_k(T^k_1)$
and $\tilde\gamma(\tilde T_1)$ respectively), which we denote by $A^k_{1,1}$
and $A^k_{1,2}$ and by $\tilde A_{1,1}$ and $\tilde A_{1,2}$, respectively
(see Figure~\ref{fig1-thm7}).

It follows from~\cite{ab} that for some further subsequence $k_n$ of the $k$'s
(which we denote by simply replacing $k$ by $n$),
$(\gamma^{\delta_n}_n,\partial A^n_{1,1},\partial A^n_{1,2})$ converge jointly
in distribution to some limit; we already know that $\gamma^{\delta_n}_n$ must
converge to $\tilde\gamma$ and want to use this fact and a suitably adapted
triangular array version of Lemma~\ref{boundaries} (whose validity does not
rely on statement (S)) to conclude that the limit is unique and coincides
with $(\tilde\gamma,\partial\tilde A_{1,1},\partial\tilde A_{1,2})$.
(Notice that $A^n_{1,1}$ and $A^n_{1,2}$ are two of the (many) domains in $G^n_1$
produced by the exploration process started at $a_n$ and stopped when it first
hits the image under $(f^n_0)^{-1}$ of the semi-circle
$\{ z : |z|=\varepsilon \} \cap {\mathbb H}$, so we are in a context close to
that of Lemma~\ref{boundaries}).

First of all, we need to show that the scaling limit of $K^{\delta_n}_{T^n_1}$
touches the image of the semi-circle $\{ z : |z|=\varepsilon \} \cap {\mathbb H}$
under $(f_0)^{-1}$ at a single point.
This follows immediately from the definition of the stopping time $T^n_1$
for every fixed $n$ (with the map $(f^n_0)^{-1}$), but it could fail to be
true in the limit $n \to \infty$.
The fact that it holds true in the limit is a direct consequence of
Lemma~\ref{double-crossing} below (in fact, a simpler version concerning
Jordan domains would suffice here, but not when we iterate the argument
-- see below), which also implies that the single point at which the scaling
limit of $K^{\delta_n}_{T^n_1}$ touches the image of the semi-circle
$\{ z : |z|=\varepsilon \} \cap {\mathbb H}$ under $(f_0)^{-1}$ coincides
with the limit of $\gamma^{\delta_n}(T^n_1)$ and with $\tilde\gamma(\tilde T_1)$.

Therefore, if we remove the single point $\tilde\gamma(\tilde T_1)$,
the scaling limit of the boundary of $K^{\delta_n}_{T^n_1}$ splits
into a left and a right part (corresponding to the scaling limit of
the leftmost yellow and the rightmost blue $\cal T$-paths of hexagons
explored by $\gamma^{\delta_n}_n$, respectively)
that do not touch the image of the semi-circle
$\{ z : |z|=\varepsilon \} \cap {\mathbb H}$ under $(f_0)^{-1}$.

Moreover, Lemma~\ref{touching} below implies that if $\gamma^{\delta_n}_n$
has a ``close encounter" with
$\partial D_n$, then it touches $\partial D^{\delta_n}_n$.
Analogously, the standard bound on the probability of six crossings
of an annulus~\cite{ksz}, used repeatedly before, implies that wherever
$\gamma^{\delta_n}_n$ has a ``close encounter" with itself, there is
touching (see the proof of Lemma~\ref{sub-conv}).
These two observations assure that the scaling limit of
$K^{\delta_n}_{T^n_1}$ is almost surely a filling (of
$\tilde\gamma[0,\tilde T_1]$), i.e., a closed connected set whose
complement in $D$ is simply connected.
From the same bound on the probability of six crossings of an annulus,
we can also conclude that the scaling limits of the left and right
boundaries of $K^n_{T^n_1}$ are almost surely simple (continuous) curves,
as in the proof of Lemma~\ref{sub-conv}.

It is also possible to conclude that the intersection of the scaling limit
of the left and right boundaries of $K^n_{T^n_1}$ with the boundary of $D$
almost surely does not contain arcs of positive length.
In fact, if that were the case, it would be possible to find a subdomain
$D'$ with three counterclockwise points $z_1,z_2,z_3$ on its boundary such
that the probability that an exploration path started at $z_1$ and stopped
when it first hits the arc $\overline{z_2 z_3}$ of $\partial D'$ has a
positive probability, in the scaling limit, of hitting at $z_2$ or $z_3$,
contradicting Cardy's formula (which, by Theorem~\ref{strong-cardy}, holds
for all subsequential scaling limits).
This means that the scaling limit of $K^n_{T^n_1}$ almost surely
satisfies the condition in~(\ref{hulls}) and is therefore a hull.
It says as well that, almost surely, the scaling limit $\tilde\gamma$ of
$\gamma^{\delta_n}_n$ does not ``stick" to the boundary of $\tilde G_1$,
which implies that also $\tilde K_{\tilde T_1}$ satisfies the condition
in~(\ref{hulls}) and is therefore a hull.

It also implies that $D \setminus \tilde K_{\tilde T_1}$ and $\tilde G_2$
are admissible domains since the part of the boundary of either
$D \setminus \tilde K_{\tilde T_1}$ or $\tilde G_2$ that belongs to the
boundary of $\tilde K_{\tilde T_1}$ can be split up, by removing the single
point $\tilde\gamma(\tilde T_1)$, into two pieces which are, by an application
of the proof of Lemma~\ref{sub-conv}, simple continuous curves, while the
remaining part of the boundary of either $D \setminus \tilde K_{\tilde T_1}$
or $\tilde G_2$ is a Jordan arc whose interior does not touch the hull
$\tilde K_{\tilde T_1}$.
(Notice however, that they need not be Jordan domains because
$\tilde K_{\tilde T_1}$ has cut-points with positive probability --
see Figure~\ref{cut-point-fig}).
This will be important later, when we need to apply
Lemma~\ref{hull} (and therefore Theorem~\ref{strong-cardy}),
Corollary~\ref{cor-unif-conv} and Lemmas~\ref{double-crossing}-\ref{touching}
again.

Then, since hulls are characterized by their ``envelope" (see Lemma~\ref{hull}
and the discussion preceding it), the joint convergence in distribution of
$\{ \partial A^n_{1,1},\partial A^n_{1,2} \}$ to
$\{ \partial\tilde A_{1,1},\partial\tilde A_{1,2} \}$ would be enough to
conclude that $K^n_{T^n_1}$ converges to $\tilde K_{\tilde T_1}$ as $n \to \infty$,
and in fact that $(\gamma^{\delta_n}_n,K^n_{T^n_1})$ converges in distribution
to $(\tilde\gamma,K_{\tilde T_1})$ (and this will be valid also for the
original subsequence $k$ and not just for the further subsequence $k_n$).
In order to get that, as explained before, we can use the convergence in
distribution of $\gamma^{\delta_n}_n$ to $\tilde\gamma$ and apply almost
the same arguments as used in the proof of Lemma~\ref{boundaries}.
The only difference is that, in proving claim (C), we cannot use the
bound on the probability of three crossings of an annulus centered at
a boundary point because we are not necessarily dealing with a convex
domain.
To replace that bound we use once again
Lemmas~\ref{double-crossing}-\ref{touching}
below (a simpler version concerning Jordan domains would again suffice
here, but not when we iterate the argument -- see below).

We can then conclude that $K^n_{T^n_1}$ converges in distribution to
$\tilde K_{\tilde T_1}$, which in turn implies the joint convergence
in distribution of $(K^k_{T^k_1},\gamma^{\delta_k}_k(T^k_1))$ to
$(\tilde K_{\tilde T_1},\tilde\gamma(\tilde T_1))$ and concludes the
first step of the argument.

%
%

We next need to prove that
$((K^k_{T^k_1},\gamma^{\delta_k}_k(T^k_1)),(K^k_{T^k_2},\gamma^{\delta_k}_k(T^k_2)))$
converges in distribution to
$((\tilde K_{\tilde T_1},\tilde\gamma(\tilde T_1)),(\tilde K_{\tilde T_2},\tilde\gamma(\tilde T_2)))$.
Since we have already proved the convergence of $(K^k_{T^k_1},\gamma^{\delta_k}_k(T^k_1))$ to
$(\tilde K_{\tilde T_1},\tilde\gamma(\tilde T_1))$, we claim that all we really need to prove
is the convergence of $(K^k_{T^k_2} \setminus K^k_{T^k_1},\gamma^{\delta_k}_k(T^k_2))$
to $(\tilde K_{\tilde T_2} \setminus \tilde K_{\tilde T_1},\tilde\gamma(\tilde T_2))$.
To see this, notice that $K^k_{T^k_2} \setminus K^k_{T^k_1}$ is distributed like
the hull of a percolation exploration path inside $D_k \setminus K^k_{T^k_1}$.
Besides, the convergence in distribution of $(K^k_{T^k_1},\gamma^{\delta_k}_k(T^k_1))$
to $(\tilde K_{\tilde T_1},\tilde\gamma(\tilde T_1))$ implies that we can find versions
of $(\gamma^{\delta_k}_k,K^k_{T^k_1})$ and $(\tilde\gamma,\tilde K_{\tilde T_1})$
on some probability space $(\Omega',{\cal B}',{\mathbb P}')$ such that
$\gamma^{\delta_k}_k(\omega')$ converges to $\tilde\gamma(\omega')$ and
$(K^k_{T^k_1},\gamma^{\delta_k}_k(T^k_1))$ converges to
$(\tilde K_{\tilde T_1},\tilde\gamma(\tilde T_1))$ for all $\omega' \in \Omega'$.
These two observations imply that, if we work with the coupled versions of
$(\gamma^{\delta_k}_k,K^k_{T^k_1})$ and $(\tilde\gamma,\tilde K_{\tilde T_1})$,
we are in the same situation as before, but with $D_k$ (resp., $D$) replaced by
$D_k \setminus K^k_{T^k_1}$ (resp., $D \setminus \tilde K_{\tilde T_1}$) and $a_k$
(resp., $a$) by $\gamma^{\delta_k}_k(T^k_1)$ (resp., $\tilde\gamma(\tilde T_1)$).
As already remarked, $D \setminus \tilde K_{\tilde T_1}$ and $\tilde G_2$ are
admissible domains, which allows us to use
Theorem~\ref{strong-cardy} (and
therefore Lemma~\ref{hull}), Corollary~\ref{cor-unif-conv} and
Lemmas~\ref{double-crossing}-\ref{touching}.


Then the conclusion that
$((K^k_{T^k_1},\gamma^{\delta_k}_k(T^k_1)),(K^k_{T^k_2},\gamma^{\delta_k}_k(T^k_2)))$
converges in distribution to
$((\tilde K_{\tilde T_1},\tilde\gamma(\tilde T_1)),(\tilde K_{\tilde T_2},\tilde\gamma(\tilde T_2)))$
follows from the same arguments as before, again using Corollary~\ref{cor-unif-conv}
and Lemmas~\ref{double-crossing}-\ref{touching}.
In order to get claim (C), in places where the exploration path comes close
to the boundary of the past hull we can use the bound on the probability of
six crossings of an annulus in the plane (as already seen in the case $k=3$
of the proof of Theorem~\ref{thm-convergence}), while in places where it comes
close to the remaining portion of the boundary (i.e.,
$\partial D$ or the Jordan arc
$\overline{cd}$ in Figure~\ref{cut-point-fig}) we can use
Lemmas~\ref{double-crossing}-\ref{touching}.

\begin{figure}[!ht]
\begin{center}
\includegraphics[width=8cm]{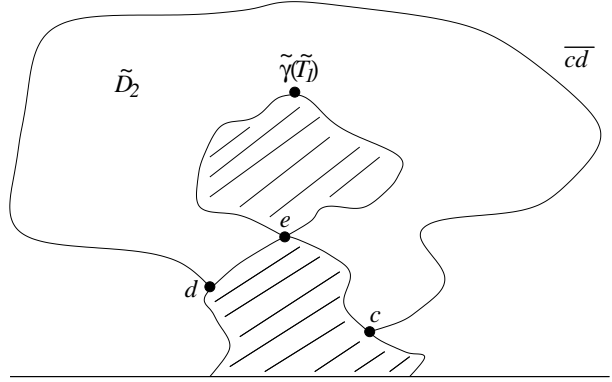}
\caption{Schematic figure representing a hull (shaded) with a cut-point $e$,
resulting in a non-Jordan, but admissible, $\tilde D_2$.}
\label{cut-point-fig}
\end{center}
\end{figure}

We can now iterate those same arguments $j$ times, for any $j>1$.
It is in fact easy to see by induction that the domains
$D \setminus \tilde K_{\tilde T_j}$ and $\tilde G_j$
that appear in the successive steps are admissible for all $j$.
Therefore we can keep using
Theorem~\ref{strong-cardy} (and therefore Lemma~\ref{hull}),
Corollary~\ref{cor-unif-conv} and Lemmas~\ref{double-crossing}-\ref{touching}.
If we keep track at each step of the previous ones, in the spirit
of Theorem~\ref{thm-convergence}, this provides the \emph{joint}
convergence of all the curves and fillings involved at each step
and concludes the proof of Theorem~\ref{spatial-markov}. \fbox{}

\begin{lemma} \label{double-crossing}
Let $\{ (D_k,a_k,c_k,d_k) \}$ be a sequence of domains admissible with
respect to $(a_k,c_k,d_k)$ and let $\gamma^{\delta}_k$ be the percolation
exploration path in $D_k$ started at $a_k$ and stopped when it first hits
the counterclockwise arc
$J'_k=\overline{c'_k d'_k} \subset J_k=\overline{c_k d_k}$ of $\partial D_k$.
Assume that, as $k \to \infty$, $(D_k,a_k,c_k,c'_k,d_k,d'_k)$ converges to
$(D,a,c,c',d,d')$, where $D$ is a domain admissible with respect to $(a,c,d)$
and $J=\overline{c'd'} \subset J=\overline{cd}$.
Let ${\cal E}^{\delta}_k(J_k;\varepsilon,\varepsilon') =
\{ \bigcup_{v \in J_k \setminus J'_k} {\cal B}^{\delta}_k(v;\varepsilon,\varepsilon') \}
\cup \{ \bigcup_{v \in J'_k} {\cal A}^{\delta}_k(v;\varepsilon,\varepsilon') \}$,
where ${\cal A}^{\delta}_k(v;\varepsilon,\varepsilon')$ is the event that
$\gamma^{\delta}_k$ contains a segment that stays within $B(v,\varepsilon)$ and
has a double crossing of the annulus $B(v,\varepsilon) \setminus B(v,\varepsilon')$
without that segment touching $\partial D^{\delta}_k$, and
${\cal B}^{\delta}_k(v;\varepsilon,\varepsilon')$ is the event that
$\gamma^{\delta}_k$ enters $B(v,\varepsilon')$, but is stopped outside
$B(v,\varepsilon)$ and does not touch $\partial D_k \cap B(v,\varepsilon)$.
Then, for any $\varepsilon>0$, 
\begin{equation} \label{double}
\lim_{\varepsilon' \to 0} \limsup_{\substack{k \to \infty \\ \delta \to 0}}
{\mathbb P}({\cal E}^{\delta}_k(J_k;\varepsilon,\varepsilon')) = 0.
\end{equation}
Essentially, this means that as $\delta \to 0$ (and $k \to \infty$),
it becomes increasingly unlikely that the exploration path ever comes
close to $J_k$ without quickly touching $\partial D^{\delta}_k$ nearby.
\end{lemma}

Lemma~\ref{double-crossing} easily implies the following result,
used in Theorem~\ref{spatial-markov} to show that
if $\tilde\gamma$ touches $\partial D$, then $\gamma^{\delta_n}_n$ touches
$\partial D^{\delta_n}_n$ nearby, for $n$ large enough.
\begin{lemma} \label{touching}
With the notation and assumptions of Lemma~\ref{double-crossing},
\begin{equation} \label{touch}
\lim_{\varepsilon \to 0} \, \lim_{\varepsilon' \to 0}
\limsup_{\substack{k \to \infty \\ \delta \to 0}}
{\mathbb P}(\bigcup_{v \in J_k \setminus J'_k}{\cal A}^{\delta}_k(v;\varepsilon,\varepsilon')) = 0.
\end{equation}
\end{lemma}

\noindent {\bf Proof.} First of all notice that for $v$ in $J_k \setminus J'_k$
but not in $B(c',\varepsilon)$ and not in $B(d',\varepsilon)$, the events
${\cal A}^{\delta}_k(v;\varepsilon,\varepsilon')$ and
${\cal B}^{\delta}_k(v;\varepsilon,\varepsilon')$ are exactly the same because
the exploration path is, by definition, stopped on $J'_k$.
Therefore, we only have to prove that the event corresponding to the union over
$v \in \{ J_k \setminus J'_k \} \cap \{ B(c',\varepsilon) \cup B(d',\varepsilon) \}$
of ${\cal A}^{\delta}_k(v;\varepsilon,\varepsilon')$ has probability
going to zero as $\varepsilon \to 0$.
We already know from Lemma~\ref{double-crossing} that
${\cal B}^{\delta}_k(v;\varepsilon,\varepsilon')$ happens with small
probability for those points.
This is, however, not sufficient because the exploration path could enter
$B(v,\varepsilon')$, then exit $B(v,\varepsilon)$, and then re-enter it
and touch $J'_k$ inside $B(v,\varepsilon)$, which is not an event in
${\cal B}^{\delta}_k(v;\varepsilon,\varepsilon')$.
But such an event would imply that $\gamma^{\delta}_k$ first
touches $J'_k$ inside one of the two balls of radius $\varepsilon$
centered at $c'_k$ and $d'_k$, and by an application of Cardy's formula
the probability that the latter happens goes to zero as
$\varepsilon \to 0$. \fbox{} \\

The proof of Lemma~\ref{double-crossing} is partly based on relating
the failure of~(\ref{double}) to the occurrence with strictly positive
probability of certain continuum limit ``mushroom" events (see
Lemma~\ref{mushroom}) that we will show must have zero probability
because otherwise there would be a contradiction to Lemma~\ref{equal},
which itself is a consequence of the continuity of Cardy's formula
with respect to the domain boundary.
In both of the next two lemmas, we denote by $\mu$ any subsequence
limit of the probability measures for the collection of all colored
(blue and yellow) $\cal T$-paths on \emph{all} of ${\mathbb R}^2$,
in the Aizenman-Burchard sense (see Remark~\ref{ab}).
We recall that in our notation, $D$ represents an open domain
and $\overline{z_1 z_2}$, $\overline{z_3 z_4}$ represent closed
segments of its boundary.
In Lemma~\ref{equal} below, we restrict attention to a Jordan domain $D$
since that case suffices for the use of Lemma~\ref{equal} in the proof
of Lemma~\ref{double-crossing}.

\begin{lemma} \label{equal}
For $(D,z_1,z_2,z_3,z_4)$, with $D$ a Jordan domain, consider the following
crossing events, ${\cal C}^*_i = {\cal C}^*_i(D,z_1,z_2,z_3,z_4)$, where $*$
denotes either blue or yellow and $i=1,2,3$:
\begin{equation} \nonumber
{\cal C}^*_1 = \{ \exists \text{ a $*$ path in the closure } \overline D
\text{ from } \overline{z_1 z_2} \text{ to } \overline{z_3 z_4} \},
\end{equation}
\begin{equation} \nonumber
{\cal C}^*_2 = \{ \exists \text{ a $*$ path in } D \text{ from the interior of }
\overline{z_1 z_2} \text{ to the interior of } \overline{z_3 z_4} \},
\end{equation}
\begin{equation} \nonumber
{\cal C}^*_3 = \{ \exists \text{ a $*$ path starting and ending \emph{outside} }
\overline D \text{ whose restriction to } D \text{ is as in } {\cal C}^*_2 \}.
\end{equation}
Then $\mu({\cal C}^*_1)=\mu({\cal C}^*_2)=\mu({\cal C}^*_3)=\Phi_D(z_1,z_2;z_3,z_4)$.
\end{lemma}

\noindent {\bf Proof.} The proof is similar to that of Theorem~\ref{strong-cardy},
but easier because $D$ is here a Jordan domain.
Indeed, it is enough to construct a new Jordan domain $\tilde D(\varepsilon)$
(with appropriately selected points $\tilde z_1(\varepsilon),\tilde z_2(\varepsilon),
\tilde z_3(\varepsilon),\tilde z_4(\varepsilon)$  on the boundary and corresponding
events $\tilde{\cal C}^*_i$) such that the occurrence of $\tilde{\cal C}^*_1$
in $\tilde D(\varepsilon)$ implies the occurrence of ${\cal C}^*_3$ in $D$ and
with $(\tilde D(\varepsilon),\tilde z_1(\varepsilon),\tilde z_2(\varepsilon),
\tilde z_3(\varepsilon),\tilde z_4(\varepsilon)) \to (D,z_1,z_2,z_3,z_4)$ as
$\varepsilon \to 0$.
The continuity of Cardy's formula (Lemma~\ref{cont-cardy} in Appendix~\ref{rado})
does the rest. \fbox{}

\begin{figure}[!ht]
\begin{center}
\includegraphics[width=8cm]{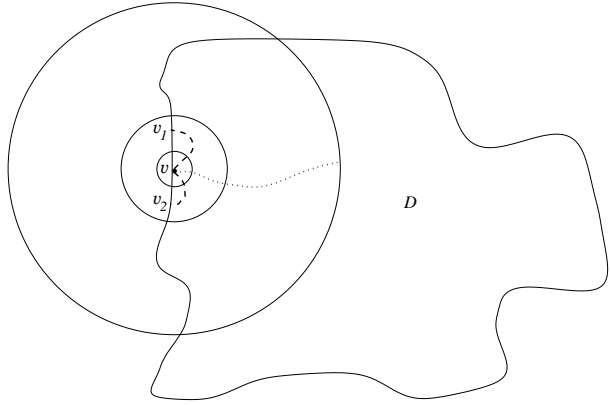}
\caption{A yellow ``mushroom" event. The dashed path is blue
and the dotted path is yellow.
The three circles centered at $v$ in the figure have radii
$\varepsilon/8$, $\varepsilon/3$, and $\varepsilon$ respectively.}
\label{fig1-lemma7-4}
\end{center}
\end{figure}

\begin{lemma} \label{mushroom}
For $(D,a,c,d)$ as in Lemma~\ref{double-crossing}, $v \in J \equiv \overline{cd}$,
and $\varepsilon>0$, we define $U^{yellow}(D,\varepsilon,v)$, the yellow
``mushroom" event (at $v$), to be the event that there is a yellow path
in $\overline D$ from $v$ to $\partial B(v,\varepsilon)$ and a blue path
in $\overline D$, between some pair of distinct points $v_1,v_2$ in
$\partial D \cap \{ B(v,\varepsilon/3) \setminus B(v,\varepsilon/8) \}$,
that passes through $v$ and such that this blue path is \emph{between}
$\partial D$ and the yellow path (see Figure~\ref{fig1-lemma7-4}).
We similarly define $U^{blue}(D,\varepsilon,v)$ with the colors interchanged
and $U^*(D,\varepsilon,J) = \cup_{v \in J} U^*(D,\varepsilon,v)$ where $*$
denotes blue or yellow.
Then for any deterministic domain $D$ and any
$0 < \varepsilon < \min \{ |a-c|, |a-d| \}$, $\mu(U^*(D,\varepsilon,J))=0$.
\end{lemma}


\begin{figure}[!ht]
\begin{center}
\includegraphics[width=8cm]{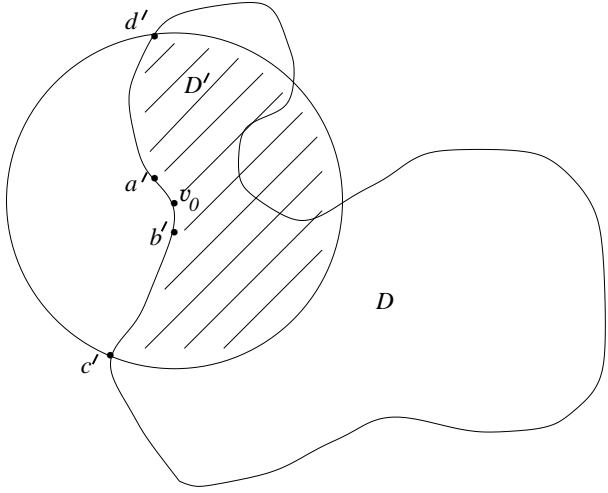}
\caption{Construction of the domain $D'$ (shaded)
used in the proof of Lemma~\ref{mushroom}.}
\label{fig2-lemma7-4}
\end{center}
\end{figure}

\noindent {\bf Proof.} If $\mu(U^*(D,\varepsilon,J))>0$ for some
$\varepsilon>0$, then there is some segment $\overline{a'b'} \subset J$
of $\partial D$ of diameter not larger than $\varepsilon/10$ such that
\begin{equation}
\mu(\cup_{v \in \overline{a'b'}} U^*(D,\varepsilon,v)) > 0.
\end{equation}
Choose any point $v_0 \in \overline{a'b'}$ and consider the new domain $D'$
whose boundary consists of the correctly chosen (as we explain below) segment
of the circle $\partial B(v_0,\varepsilon/2)$ between the two points $c',d'$
where $\partial D$ first hits $\partial B(v_0,\varepsilon/2)$ on
either side of $v_0$, together with the segment from $d'$ to $c'$ of
$\partial D$ (see Figure~\ref{fig2-lemma7-4}).
The correct circle segment between $c'$ and $d'$ is the (counter) clockwise
one if $v_0$ is between $c'$ and $d'$ along $\partial D$ when $\partial D$
is oriented (counter) clockwise.
It is also not hard to see that since $\varepsilon < \min \{ |a-c|, |a-d| \}$,
$D'$ is a Jordan domain, so that Lemma~\ref{equal} can be applied.
In the new domain $D'$, $\overline{a'b'}$ is the same curve segment as it
was in the old domain $D$, but $\overline{c'd'}$ is now a segment of the
circle $\partial B(v_0,\varepsilon/2)$.
It should be clear that
\begin{equation}
\cup_{v \in \overline{a'b'}} U^*(D,\varepsilon,v) \subset
{\cal C}^*_1(D',a',b',c',d') \setminus {\cal C}^*_3(D',a',b',c',d')
\end{equation}
which yields a contradiction of Lemma~\ref{equal} if
$\mu(U^*(D,\varepsilon,J))>0$. \fbox{} \\

\noindent {\bf Proof of Lemma~\ref{double-crossing}.} We first note that
since the probability in~(\ref{double}) is nonincreasing in $\varepsilon$,
we may assume that $\varepsilon < \min \{|a-c|,|a-d| \}$, as requested by
Lemma~\ref{mushroom}.

Let us first consider the simpler case of
${\cal B}^{\delta_k}(v;\varepsilon,\varepsilon')$ in which $v \in J'_k$.
We follow the exploration process until time $T$, when it first
touches $\partial B(v,\varepsilon')$ for some $v \in J'_k$, and
consider the annulus $B(v,\varepsilon) \setminus B(v,\varepsilon')$.
Let $\pi_Y$ be the leftmost yellow $\cal T$-path and $\pi_B$
the rightmost blue $\cal T$-path in $\Gamma(\gamma^{\delta}_k)$
at time $T$ that cross $B(v,\varepsilon) \setminus B(v,\varepsilon')$.
$\pi_Y$ and $\pi_B$ split the annulus
$B(v,\varepsilon) \setminus B(v,\varepsilon')$ into three sectors
that, for simplicity, we will call the central sector, containing
the crossing segment of the exploration path, the yellow (left) sector,
with $\pi_Y$ as part of its boundary, and the blue (right) sector,
the remaining one, with $\pi_B$ as part of its boundary.

We then look for a yellow ``lateral" crossing within the yellow sector
from $\pi_Y$ to $\partial D_k$ and a blue lateral crossing within
the blue sector from $\pi_B$ to $\partial D_k$.
Notice that the yellow sector may contain ``excursions"
of the exploration path coming off $\partial B(v,\varepsilon)$,
producing nested yellow and blue excursions off
$\partial B(v,\varepsilon)$, and the same for the blue sector.
But for topological reasons, those excursions are such that
for every group of nested excursions, the outermost one is
always yellow in the yellow sector and blue in the blue sector.
Therefore, by standard percolation theory arguments, the conditional
probability (conditioned on $\Gamma(\gamma^{\delta}_k)$ at time $T$)
to find a yellow lateral crossing of the yellow sector from $\pi_Y$
to $\partial D_k$ is bounded below by the probability to find a
yellow circuit in an annulus with inner radius $\varepsilon'$
and outer radius $\varepsilon$.
An analogous statement holds for the conditional probability
(conditioned on $\Gamma(\gamma^{\delta}_k)$ at time $T$ and
also on the entire percolation configuration in the yellow
sector) to find a blue lateral crossing of the blue sector
from $\pi_B$ to $\partial D_k$.
Thus for any fixed $\varepsilon>0$, by an application of the
Russo-Seymour-Welsh lemma~\cite{russo,sewe}, the conditional
probability to find both a yellow lateral crossing within the
yellow sector from $\pi_Y$ to $\partial D_k$ and a blue lateral
crossing within the blue sector from $\pi_B$ to $\partial D_k$
goes to one as $\varepsilon'\to 0$.

But if such yellow and blue crossings are present, the exploration
path is forced to touch $J'_k$ before exiting $B(v,\varepsilon)$,
and if that happens, the exploration process is stopped, so that
it will never exit $B(v,\varepsilon)$ and the union over $v \in J'_k$
of ${\cal B}^{\delta}_k(v;\varepsilon,\varepsilon')$ cannot occur.
This concludes the proof of this case.

Let us now consider the remaining case in which $v \notin J'_k$.
The basic idea of the proof is then that by straightforward weak
convergence and related coupling arguments, the failure of~(\ref{double})
would imply that \emph{some} subsequence limit $\mu$ would satisfy
$\mu(U^{yellow}(D,\varepsilon,J) \cup U^{blue}(D,\varepsilon,J))>0$,
which would contradict Lemma~\ref{mushroom}.
This is essentially because the close approach of an exploration
path on the $\delta$-lattice to $J_k \setminus J'_k$ without quickly
touching nearby yields one two-sided colored $\cal T$-path (the
``perimeter" of the portion of the hull of the exploration path
seen from a boundary point of close approach) and a one-sided
$\cal T$-path of the other color belonging to the percolation
cluster not seen from the boundary point (i.e., shielded by the
two-sided path).
Both the two-sided path and the one-sided one are subsets of
$\Gamma(\gamma^{\delta}_k)$.

\begin{figure}[!ht]
\begin{center}
\includegraphics[width=8cm]{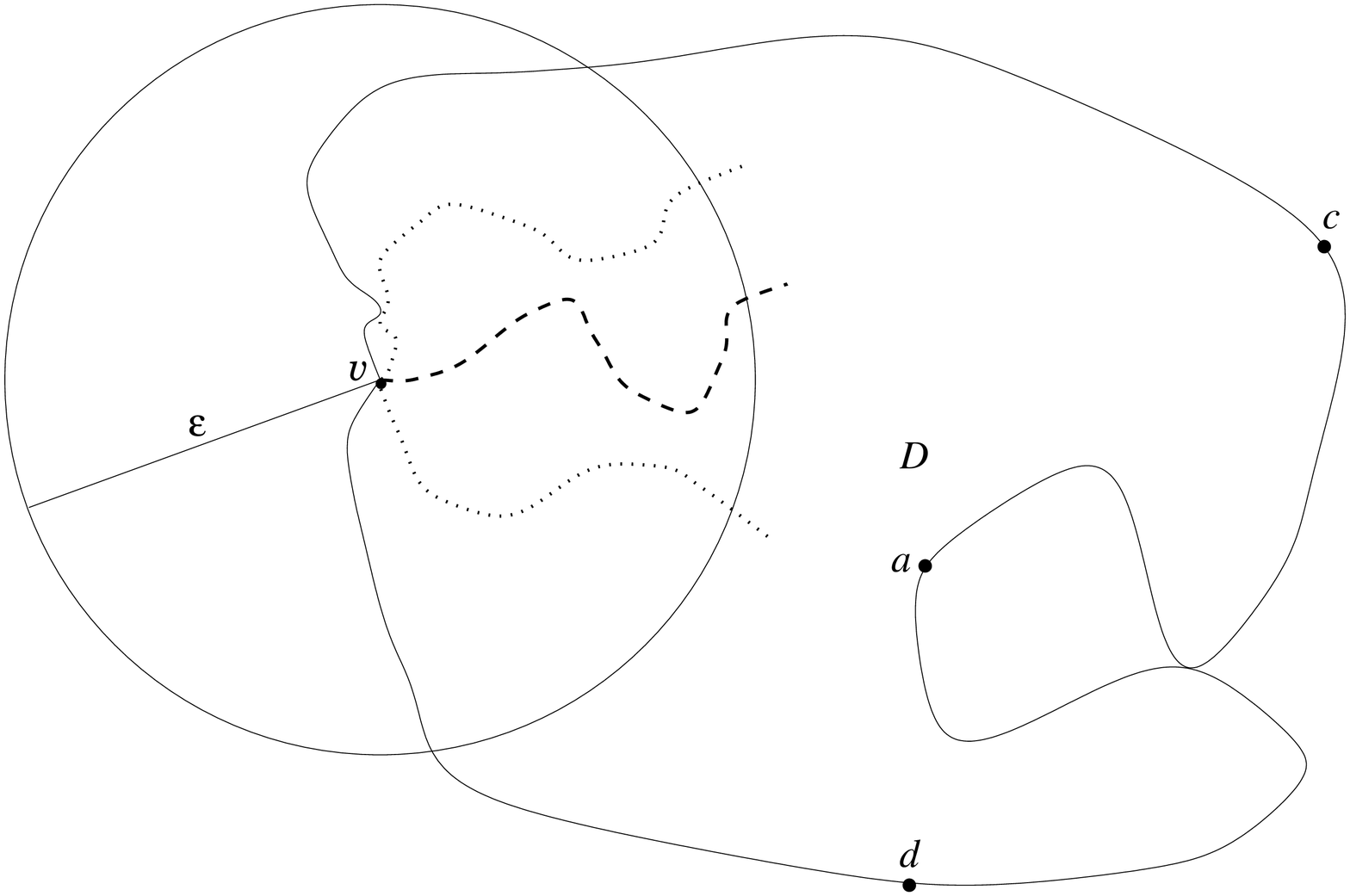}
\caption{The event consisting of a yellow double crossing
and a blue crossing used as a first step for obtaining a
blue mushroom event in the proof of Lemma~\ref{double-crossing}.
The dashed crossing is blue and the dotted crossing is yellow.}
\label{fig1-lemma7-2}
\end{center}
\end{figure}

Assume by contradiction that~(\ref{double}) is false, so that close
encounters without touching happen with bounded away from zero probability.
Consider for concreteness an exploration path $\gamma^{\delta}_k$
that has a close approach to a point $v$ in the counterclockwise
arc $\overline{d'_k d_k}$.
The exploration path may have multiple close approaches to $v$ with
differing colors of the perimeter as seen from $v$, but for topological
reasons, the last time the exploration path comes close to $v$, it
must do so in such a way as to produce a yellow $\cal T$-path $\pi_Y$
(seen from $v$) that crosses $B(v,\varepsilon) \setminus B(v,\varepsilon')$
twice, and a blue path $\pi_B$ that crosses it once
(see Figure~\ref{fig1-lemma7-2}).
This is so because the exploration process that produced
$\gamma^{\delta}_k$ ended somewhere on $J'_k$ (and outside
$B(v,\varepsilon)$), which is to the right of (i.e., clockwise to) $v$.

The presence of $\pi_Y$ implies that there are a yellow leftmost
$\cal T$-path $\pi_L$ and a yellow rightmost $\cal T$-path $\pi_R$
(looking at $v$ from inside $D_k$) crossing the annulus
$B(v,\varepsilon) \setminus B(v,\varepsilon')$.
The paths $\pi_L$ and $\pi_R$ split the annulus
$B(v,\varepsilon) \setminus B(v,\varepsilon')$ into three sectors,
that we will call the central sector, containing $\pi_B$, the left
sector, with $\pi_L$ as part of its boundary, and the right sector,
with $\pi_R$ as part of its boundary.
Again for topological reasons, all other monochromatic crossings
of the annulus are contained in the central sector, including at
least one blue path $\pi_B$.
As in the previous case, the left and right sectors can contain
nested monochromatic excursions off $\partial B(v,\varepsilon)$,
but this time for every group of excursions, the outermost one is
yellow in both sectors.

Now consider the annulus $B(v,\varepsilon/3) \setminus B(v,\varepsilon/8)$.
We look for a yellow lateral crossing within the left sector
from $\pi_L$ to $\partial D_k$ and a yellow lateral crossing
within the right sector from $\pi_R$ to $\partial D_k$.
Since the outermost excursions in both sectors are yellow,
the conditional probability to find a yellow lateral crossing
within the left sector from $\pi_L$ to $\partial D_k$ is bounded
below by the probability to find a yellow circuit in an annulus with
inner radius $\varepsilon'$ and outer radius $\varepsilon$,
and an analogous statement holds for the conditional probability
to find a yellow lateral crossing within the right sector from
$\pi_R$ to $\partial D_k$.
Thus for any fixed $\varepsilon>0$, by an application of the
Russo-Seymour-Welsh lemma~\cite{russo,sewe}, the conditional
probability to find both yellow lateral crossings goes to one
as $\varepsilon'\to 0$.
But the presence of such yellow crossings would produce a (blue)
mushroom event, leading to a contradiction with Lemma~\ref{mushroom}. \fbox{} \\

We are finally ready to prove the main result of this section which
implies statement (S) at the beginning of Section~\ref{main-results}.

\begin{corollary} \label{conv-to-sle-1}
Consider a sequence $\{(D_k,a_k,b_k)\}$ of Jordan domains with two
distinct selected points $a_k,b_k$ on their boundaries $\partial D_k$.
Assume that $(D_k,a_k,b_k) \to (D,a,b)$, where $D$ is a Jordan domain
with two distinct selected points on its boundary $\partial D$.
Denote by $\gamma^{\delta}_k \equiv \gamma^{\delta}_{D_k,a_k,b_k}$ the
percolation exploration path inside (the $\delta$-approximation
$D^{\delta}_k$ of) $D_k$ from (the e-vertex closest to) $a_k$ to
(the e-vertex closest to) $b_k$.
Then for any sequence $\delta_k \downarrow 0$, as $k \to \infty$,
$\gamma^{\delta_k}_k$ converges in distribution to the trace $\gamma_{D,a,b}$
of chordal $SLE_6$ inside $D$ from $a$ to $b$.
\end{corollary}

\noindent {\bf Proof.} It follows from~\cite{ab} that $\gamma^{\delta_k}_k$
converges in distribution along subsequence limits $k_n$.
Since we have proved that the filling of any such subsequence
limit $\tilde\gamma$ satisfies the spatial Markov property
(Theorem~\ref{spatial-markov}) and the exit distribution of
$\tilde\gamma$ is determined by Cardy's formula (Theorem~\ref{strong-cardy}),
we can deduce from Theorem~\ref{characterization} that the limit is unique
and that the law of $\gamma^{\delta_k}_k$ converges, as $k \to \infty$,
to the law of the trace $\gamma_{D,a,b}$ of chordal $SLE_6$ inside $D$
from $a$ to $b$. \fbox{}

\refstepcounter{section}
\section*{Appendix \thesection: Sequences of Conformal Maps} \label{rado}

In this appendix, we give some results about sequences of conformal maps
that are used in various places throughout the paper.
For more details, the interested reader should consult~\cite{pommerenke}.


\begin{definition} \emph{(see Section~1.4 of~\cite{pommerenke})} \label{kernel-conv}
Let $w_0 \in {\mathbb C}$ be given and let $\{ G_n \}$ be
domains with $w_0 \in G_n \subset {\mathbb C}$.
We say that $G_n \to G$ as $n \to \infty$ with respect to $w_0$
in the sense of {\bf kernel convergence} if
\begin{enumerate}
\item either $G = \{ w_0 \}$, or else $G$ is a domain $\neq {\mathbb C}$
with $w_0 \in G$ such that some neighborhood of every $w \in G$ lies
in $G_n$ for large $n$; and
\item for $w \in \partial G$ there exist $w_n \in \partial G_n$ such
that $w_n \to w$ as $n \to \infty$.
\end{enumerate}
\end{definition}

It is clear from the definition that every subsequence limit also converges
to $G$ and it is also easy to see that the limit is uniquely determined.
With this definition we can now state Carath\'eodory's kernel theorem~\cite{caratheodory}.

\begin{theorem} \emph{(see Theorem~1.8 of~\cite{pommerenke})} \label{kernel-thm}
Let $f_n$ map ${\mathbb D}$ conformally onto $G_n$ with $f_n(0)=w_0$
and $f'_n(0)>0$.
If $G = \{ w_0 \}$, let $f(z) \equiv w_0$; otherwise let $f$ map
$\mathbb D$ conformally onto $G$ with $f(0)=w_0$ and $f'(0)>0$.
Then, as $n \to \infty$, $f_n \to f$ locally uniformly in ${\mathbb D}$
if and only if $G_n \to G$ with respect to $w_0$.
\end{theorem}

The next result, Rad\'o's theorem~\cite{rado}, deals with sequences
of Jordan domains and is used in the main body of the paper.
In this case the conformal maps have a continuous extension to
${\mathbb D} \cup \partial {\mathbb D}$.
\begin{theorem} \emph{(see Theorem~2.11 of~\cite{pommerenke})} \label{rado-thm}
For $n=1,2,\ldots$ , let $J_n$ and $J$ be Jordan curves parametrized
respectively by $\phi_n(t)$ and $\phi(t)$, $t \in [0,1]$, and let
$f_n$ and $f$ be conformal maps from $\mathbb D$ onto the inner
domains of $J_n$ and $J$ such that $f_n(0)=f(0)$ and $f_n'(0)>0$,
$f'(0)>0$ for all $n$.
If $\phi_n \to \phi$ as $n \to \infty$ uniformly in $[0,1]$
then $f_n \to f$ as $n \to \infty$ uniformly in $\overline{\mathbb D}$.
\end{theorem}

The type of convergence of sequences of Jordan domains
$\{ G_n \}$ to a Jordan domain $G$ encountered in the main body of the
paper (i.e., in the sense that $\partial G_n$ converges, as $n \to \infty$,
to $\partial G$ in the uniform metric~(\ref{distance}) on continuous curves)
is clearly sufficient to apply Theorem~\ref{rado-thm}.
%
%
In Appendix~\ref{convergence-to-sle}, however, we have to deal with
domains that are not Jordan, and therefore we cannot use Rad\'o's theorem.
The tools needed to deal with those situations are described below.

\begin{definition} \emph{(see Section~2.2 of~\cite{pommerenke})} \label{locally-connected}
The closed set $A \subset {\mathbb C}$ is called {\bf locally connected}
if for every $\varepsilon>0$ there is $\delta>0$ such that, for any two
points $a,b \in A$ with $|a-b|<\delta$, we can find a continuum $B$ with
diameter smaller than $\varepsilon$ and with $a,b \in B \subset A$.
\end{definition}

In the definition above, a {\bf continuum} denotes a compact
connected set with more than one point.
We remark that every continuous curve (with more than one point) is a
locally connected continuum (the converse is also true: every locally
connected continuum is a curve).
The concept of local connectedness gives a topological answer to the problem
of global extension of a conformal map to the domain boundary, as follows.
\begin{theorem} \emph{(see Continuity Theorem in Section~2.1 of~\cite{pommerenke})}
\label{cont-thm}
Let $f$ map the unit disk $\mathbb D$ conformally onto $G \subset {\mathbb C} \cup \{\infty\}$.
Then the function $f$ has a continuous extension to ${\mathbb D} \cup \partial {\mathbb D}$
if and only if $\partial G$ is locally connected.
\end{theorem}

When $f$ has a continuous extension to ${\mathbb D} \cup \partial {\mathbb D}$,
we do not distinguish between $f$ and its extension.
This is always the case for the conformal maps considered in this paper.
The problem wether this extension is injective on $\overline{\mathbb D}$
has also a topological answer, as follows.
\begin{theorem} \emph{(see Carathh\'eodory Theorem in Section~2.1 of~\cite{pommerenke})}
\label{cara-thm}
In the notation of Theorem~\ref{cont-thm}, the function $f$ has a continuous
and injective extension if and only if $\partial G$ is a Jordan curve.
\end{theorem}


When considering sequences of domains whose boundaries are locally connected
the following definition is useful.
\begin{definition} \emph{(see Section~2.2 of~\cite{pommerenke})}
\label{unif-locally-connected}
The closed sets $A_n \subset {\mathbb C}$ are {\bf uniformly locally connected}
if, for every $\varepsilon>0$, there exists $\delta>0$ independent of $n$ such
that any two points $a_n,b_n \in A_n$ with $|a_n-b_n|<\delta$ can be joined by
continua $B_n \subset A_n$ of diameter smaller than $\varepsilon$.
\end{definition}

The convergence of domains used in this paper (i.e., $G_n \to G$ if
$\partial G_n \to \partial G$ in the uniform metric~(\ref{distance})
on continuous curves) clearly implies kernel convergence, which
immediately allows us to use Theorem~\ref{kernel-thm}.
However, we need uniform convergence in $\overline{\mathbb D}$.
This is guaranteed by Rad\'o's theorem in the case of Jordan domains;
in the non-Jordan case, sufficient conditions to have uniform convergence
are stated in the next theorem.

\begin{theorem} \emph{(see Corollary~2.4 of~\cite{pommerenke})}
\label{unif-conv}
Let $\{ G_n \}$ be a sequence of bounded domains such that,
for some $0<r<R<\infty$, $B(0,r) \subset G_n \subset B(0,R)$
for all $n$ and such that $\{ {\mathbb C} \setminus G_n \}$
is uniformly locally connected.
Let $f_n$ map $\mathbb D$ conformally onto $G_n$ with $f_n(0)=0$.
If $f_n(z) \to f(z)$ as $n \to \infty$ for each $z \in {\mathbb D}$,
then the convergence is uniform in $\overline{\mathbb D}$.
\end{theorem}

In order to use Theorem~\ref{unif-conv} in Appendix~\ref{convergence-to-sle}
we need the following lemma.
The definitions of admissible domain and the related notion of
convergence are given just before Theorem~\ref{strong-cardy} in
Appendix~\ref{convergence-to-sle}.

\begin{lemma} \label{lemma-unif-loc-conn}
Let $\{ (G_n,a_n,c_n,d_n) \}$ be a sequence of domains admissible
with respect to $(a_n,c_n,d_n)$ and assume that, as $n \to \infty$,
$(G_n,a_n,c_n,d_n) \to (G,a,c,d)$, where $G$ is a domain admissible
with respect to $(a,c,d)$.
Then the sequence of closed sets $\{ {\mathbb C} \setminus G_n \}$
is uniformly locally connected.
\end{lemma}

\noindent{\bf Proof.} In order to prove the lemma, we claim that it
suffices to focus on pairs of points on the boundaries, i.e., to show that:
$(*)$ for every $\varepsilon>0$, there exists $\delta=\delta(\varepsilon)>0$
independent of $n$ such that any two points $u_n,v_n \in \partial G_n$ with
$|u_n-v_n| < \delta$ can be joined by a continuum of diameter $<\varepsilon$
contained in the complement ${\mathbb C} \setminus G_n$ of $G_n$.

To verify our claim, let us assume $(*)$ for the moment, and consider
two points $u_n,v_n \in {\mathbb C} \setminus G_n$ (but not necessarily
in $\partial G_n$) with $|u_n-v_n| < \delta'$, where
$\delta' = \min \{ \frac{1}{3}\delta(\frac{\varepsilon}{3}), \frac{\varepsilon}{3} \}$.
If (at least) one of the two points, say $u_n$, is at distance greater
than $\delta'$ from $\partial G_n$, then we can connect $u_n$ and $v_n$
using the closed ball of radius $\delta'$ centered at $u_n$, since
$v_n \in \overline{B(u_n,\delta')} \subset {\mathbb C} \setminus G_n$.
If both points are at distance smaller than $\delta'$ from $\partial G_n$,
we can connect each point to a closest point on $\partial G_n$ by a
straight segment of length smaller than $\delta'$.
Those two points on $\partial G_n$ can then be connected to each
other by a continuum $B_n$ of diameter $<\varepsilon/3$ contained
in ${\mathbb C} \setminus G_n$, and the union of $B_n$ with
the two straight segments gives a continuum of diameter $<\varepsilon$
connecting $u_n$ with $v_n$ and contained in ${\mathbb C} \setminus G_n$.

We now prove $(*)$.
Since $\partial G_n \to \partial G$ in the uniform metric~(\ref{distance})
on continuous curves, for every $\varepsilon>0$ there exists
$n_0=n_0(\varepsilon)$ such that for all $n \geq n_0$,
$\text{d}(\partial G_n, \partial G) < \varepsilon$.
The admissibility of $G$ implies that we can split its boundary
into three Jordan arcs, $J_1 = \overline{da}$, $J_2 =\overline{ac}$,
$J_3 = \overline{cd}$, such that $J_3$ does not touch the interior
of either $J_1$ or $J_2$.
We can do the same with $\partial G_n$, letting $J_{1,n} = \overline{d_n a_n}$,
$J_{2,n} =\overline{a_n c_n} $ and $J_{3,n} = \overline{c_n d_n}$.
Let $\phi_{i,n}(t)$ and $\phi_i(t)$, $t \in [0,1]$, $i=1,2,3$ be
parametrizations of $J_{i,n}$ and $J_i$ respectively, with
$\sup_{t \in [0,1]} |\phi_{i,n}(t)-\phi_i(t)| < \varepsilon$ for
$n \geq n_0$ and $i=1,2,3$.

Let us assume, by contradiction, that $(*)$ is false.
Then there are indices $k$ (actually $n_k$, but we abuse notation a
bit) and points $u_k, v_k \in \partial G_k$ with $|u_k - v_k| \to 0$
(as $k \to \infty$) that cannot be joined by a continuum of diameter
$<\varepsilon$ contained in ${\mathbb C} \setminus G_k$.
By compactness considerations, we may assume that $u_k \to u$ and
$v_k \to v$ as $k \to \infty$, with $u=v$.
Suppose that $u_k$ and $v_k$ belong to the interior of the same Jordan
arc $J_{i,k}$ for all $k$ large enough.
Let $u_k = \phi_{i,k}(\tau_k)$, $v_k = \phi_{i,k}(\tau'_k)$, $u = \phi_i(\tau)$
and $v = \phi_i(\tau')$.
It follows that $\tau_k \to \tau$ and $\tau'_k \to \tau'$, and since $J_i$
is a Jordan arc, $\tau=\tau'$.
For $k$ large enough, the function $\phi_{i,k}$ maps the closed segment of
$[0,1]$ between $\tau$ and $\tau'$ onto a continuum in $J_{i,k}$ containing
$u_k$ and $v_k$ whose diameter tends to zero as $k \to \infty$, leading to
a contradiction with our assumption.

Similar reasoning gives a contradiction if $u_k$ and $v_k$ both belong
to $J_{1,k} \cup J_{3,k}$ or both belong to $J_{2,k} \cup J_{3,k}$ for
all $k$ large enough, since the concatenation of $J_{1,k}$ with $J_{3,k}$
or of $J_{2,k}$ with $J_{3,k}$ is still a Jordan arc.
The above reasoning applies except when $u(=v)$ is on both $J_1$ and $J_2$.
When $u=v=a$, one can paste together small Jordan arcs on $J_{1,k}$ and
$J_{2,k}$ to get a suitable continuum leading to a contradiction.
The sole remaining case is when for all $k$ large enough, $u_k$ belongs
to the interior of $J_{1,k}$ and $v_k$ belongs to the interior of $J_{2,k}$.

(Notice that we are ignoring the ``degenerate" case in which $c=d$
coincides with the ``last" [from $a$] double-point on $\partial G$,
and $J_3$ is a simple loop.
In that case $u_k$ and $v_k$ could converge to $u=v=c=d \in J_1 \cap J_2$
and $u_k$ or $v_k$ could still belong to $J_{3,k}$ for arbitrarily large $k$'s.
However, in that case one can find two distinct points on $J_3$, $c'$
and $d'$, such that $D$ is admissible with respect to $(a,c',d')$, and
points $c'_k$ and $d'_k$ on $J_{3,k}$ converging to $c'$ and $d'$
respectively, and define accordingly new Jordan arcs, $J'_1, J'_2, J'_3$
and $J'_{1,k}, J'_{2,k}, J'_{3,k}$, so that $u_k \in J'_{1,k}$ and
$v_k \in J'_{2,k}$ for $k$ large enough.
We assume that this has been done if necessary, and for simplicity
of notation drop the primes.)

In this case let $[u_k v_k]$ denote the closed straight line segment
in the plane between $u_k$ and $v_k$.
Imagine that $[u_k v_k]$ is oriented from $u_k$ to $v_k$ and let $v'_k$
be the first point of $J_{2,k}$ intersected by $[u_k v_k]$ and $u'_k$
be the previous intersection of $[u_k v_k]$ with $\partial G_k$.
Clearly, $u'_k \notin J_{2,k}$.
For $k$ large enough, $u'_k$ cannot belong to $J_{3,k}$ either, or
otherwise in the limit $k \to \infty$, $J_3$ would touch the interior
of $J_1$ and $J_2$.
We deduce that for all $k$ large enough, $u'_k \in J_{1,k}$.
Since $J_{1,k}$ and $J_{2,k}$ are continuous curves and therefore
locally connected, $u_k$ and $u'_k$ belong to a continuum $B_{1,k}$
contained in $J_{1,k}$ whose diameter goes to zero as $k \to \infty$,
and the same for $v_k$ and $v'_k$ (with $B_{1,k}$ and $J_{1,k}$
replaced by $B_{2,k}$ and $J_{2,k}$).

Since the interior of $[u'_k v'_k]$ does not intersect any
portion of $\partial G_k$, it is either contained in $G_k$
or in its complement ${\mathbb C} \setminus G_k$.
If $[u'_k v'_k] \subset {\mathbb C} \setminus G_k$, we have a
contradiction since the union of $[u'_k v'_k]$ with $B_{1,k}$
and $B_{2,k}$ is contained in ${\mathbb C} \setminus G_k$ and
is a continuum containing $u_k$ and $v_k$ whose diameter goes
to zero as $k \to \infty$.

If the interior of $[u'_k v'_k]$ is contained in $G_k$, let us
consider a conformal map $f_k$ from $\mathbb D$ onto $G_k$.
Since $\partial G_k$ is locally connected, the conformal map
$f_k$ extends continuously to the boundary of the unit disc.
Let $u'_k = f_k(u^*_k)$, $v'_k = f_k(v^*_k)$, $a_k = f_k(a^*_k)$,
$c_k = f_k(c^*_k)$ and $d_k = f_k(d^*_k)$.
The points $c^*_k, d^*_k, u^*_k, a^*_k, v^*_k$ are in counterclockwise
order on $\partial {\mathbb D}$, so that any curve in $\mathbb D$
from $a^*_k$ to the counterclockwise arc $\overline{c^*_k d^*_k}$
must cross the curve from $u^*_k$ to $v^*_k$ whose image under $f_k$
is $[u'_k v'_k]$.
This implies that any curve in $G_k$ going from $a_k$ to the
counterclockwise arc $\overline{c_k d_k}$ of $\partial G_k$ must
cross the (interior of the) line segment $[u'_k v'_k]$.
Then, in the limit $k \to \infty$, any curve in $G$ from $a$ to
the counterclockwise arc $\overline{cd}$ must contain the limit
point $u=\lim_{k \to \infty}u'_k=\lim_{k \to \infty}v'_k=v$.
On the other hand, except for its starting and ending point, any
such curve is completely contained in $G$, which implies that either
$u=v=a$ or else that (in the limit $k \to \infty$) the counterclockwise
arc $\overline{cd}$ is the single point at $u=v=c=d$.
We have already dealt with the former case.
In the latter case, one can paste together small Jordan arcs
from $u'_k$ to $d_k$, from $d_k$ to $c_k$, and from $c_k$ to $v'_k$,
and take the union with $B_{1,k}$ and $B_{2,k}$ (defined above) to
get a suitable continuum in ${\mathbb C} \setminus G_k$
containing $u_k$ and $v_k$, leading to a contradiction.
This concludes the proof. \fbox{} \\

%

Theorem~\ref{unif-conv}, together with Theorem~\ref{kernel-thm} and
Lemma~\ref{lemma-unif-loc-conn}, implies the following result, which
is used in Appendix~\ref{convergence-to-sle}.

\begin{corollary} \label{cor-unif-conv}
With the notation and assumptions of Lemma~\ref{lemma-unif-loc-conn}
(and also assuming that $G_n$ and $G$ contain the origin), let $f_n$ map
${\mathbb D}$ conformally onto $G_n$ with $f_n(0)=0$ and $f'_n(0)>0$,
and $f$ map $\mathbb D$ conformally onto $G$ with $f(0)=0$ and $f'(0)>0$.
Then, as $n \to \infty$, $f_n \to f$ uniformly in $\overline{\mathbb D}$.
\end{corollary}

\noindent {\bf Proof.} As already remarked, the convergence of $\partial G_n$
to $\partial G$ in the uniform metric~(\ref{distance}) on continuous curves
(which is part of the definition of $(G_n,a_n,c_n,d_n) \to (G,a,c,d)$)
easily implies that the conditions in Carath\'eodory's kernel theorem
(Theorem~\ref{kernel-thm}) are satisfied and therefore that $f_n$
converges to $f$ locally uniformly in $\mathbb D$, as $n \to \infty$.
By an application of Lemma~\ref{lemma-unif-loc-conn}, the sequence
$\{ {\mathbb C} \setminus D_n \}$ is uniformly locally connected, so that
we can apply Theorem~\ref{unif-conv} to conclude that, as $n \to \infty$,
$f_n$ converges to $f$ uniformly in $\overline{\mathbb D}$. \fbox{} \\

We conclude this appendix with a simple lemma, used in the proof of
Theorem~\ref{strong-cardy}, about the continuity of Cardy's formula
with respect to the shape of the domain and the positions of the four
points on the boundary.

\begin{lemma} \label{cont-cardy}
For $\{ (D_n,a_n,b_n,c_n,d_n) \}$ and $(D,a,b,c,d)$ as in
Theorem~\ref{strong-cardy}, let $\Phi_n$ denote Cardy's formula
(see~(\ref{cardy-formula})) for a crossing inside $D_n$ from the
counterclockwise segment $\overline{a_n c_n}$ of $\partial D_n$ to
the counterclockwise segment $\overline{b_n d_n}$ of $\partial D_n$
and $\Phi$ the corresponding Cardy's formula for the limiting domain $D$.
Then, as $n \to \infty$, $\Phi_n \to \Phi$.
\end{lemma}

\noindent {\bf Proof.} Let $f_n$ be the conformal map that takes
$\mathbb D$ onto $D_n$ with $f_n(0)=0$ and $f'_n(0)>0$, and let $f$
denote the conformal map from $\mathbb D$ onto $D$ with $f(0)=0$
and $f'(0)>0$; let $z_1 = f^{-1}(a)$, $z_2 = f^{-1}(c)$, $z_3 = f^{-1}(b)$,
$z_4 = f^{-1}(d)$, $z_1^n = f_n^{-1}(a_n)$, $z_2^n = f_n^{-1}(c_n)$,
$z_3^n = f_n^{-1}(b_n)$, and $z_4^n = f_n^{-1}(d_n)$.
We can apply Corollary~\ref{cor-unif-conv} to conclude that, as $n \to \infty$,
$f_n$ converges to $f$ uniformly in $\overline{\mathbb D}$.
This, in turn, implies that, as $n \to \infty$, $z_1^n \to z_1$,
$z_2^n \to z_2$, $z_3^n \to z_3$, and $z_4^n \to z_4$.

Cardy's formula for a crossing inside $D_n$ from the counterclockwise
segment $\overline{a_n c_n}$ of $\partial D_n$ to the counterclockwise
segment $\overline{b_n d_n}$ of $\partial D_n$ is given by
\begin{equation} \label{cardy1}
\Phi_n =
\frac{\Gamma(2/3)}{\Gamma(4/3) \Gamma(1/3)} \eta^{1/3}_n
{}_2F_1(1/3,2/3;4/3;\eta_n),
\end{equation}
where
\begin{equation}
\eta_n =
\frac{(z_1^n-z_2^n)(z_3^n-z_4^n)}{(z_1^n-z_3^n)(z_2^n-z_4^n)}.
\end{equation}
Because of the continuity of $\eta_n$ in $z_1^n$, $z_2^n$, $z_3^n$,
$z_4^n$, and the continuity of Cardy's formula~(\ref{cardy1}) in $\eta_n$,
the convergence of $z_1^n \to z_1$, $z_2^n \to z_2$, $z_3^n \to z_3$ and
$z_4^n \to z_4$ immediately implies the convergence of $\Phi_n$ to $\Phi$. \fbox{}

\bigskip
\bigskip

\noindent {\bf Acknowledgements.} We are grateful to
Greg Lawler, Oded Schramm and Wendelin Werner for various
interesting and useful conversations and to Stas Smirnov
for communications about a paper in preparation.
We note in particular that a discussion with Oded Schramm
at the November 2004 Northeast Probability Seminar at the
CUNY Graduate Center, about dependence of exploration paths
with respect to small changes of domain boundaries, pointed
us in a direction that eventually led to
Lemmas~\ref{double-crossing}-\ref{mushroom}.
F.~C. thanks Wendelin Werner for an invitation to Universit\'e
Paris-Sud 11, and Vincent Beffara and Luiz Renato Fontes for
many helpful discussions.
We are especially grateful to Vincent Beffara for pointing out a
gap in a preliminary version of Appendix~\ref{convergence-to-sle}.
We thank Lai-Sang Young for comments about the presentation
of our results.
F.~C. acknowledges the kind hospitality of the Courant Institute
where part of this work was completed.

\end{document}